\documentclass[11pt]{article}
\usepackage{mathrsfs}
\usepackage{}

\usepackage{bbding}
\usepackage{pifont}
\usepackage{amsfonts}
\usepackage{geometry}
\usepackage{amsmath}
\usepackage{amsthm}
\usepackage{amssymb}
\usepackage{mathptmx}
\usepackage{graphicx}
\usepackage{bbding}
\usepackage{tikz}
\usepackage{amscd}

\def\nd{\noindent}

\def\<{\leq}
\def\>{\geq}

\newtheorem{thm}{Theorem}[section]
\newtheorem{lem}{Lemma}[section]
\newtheorem{prop}{Proposition}[section]
\newtheorem{cor}{Corollary}[section]
\newtheorem{defi}{Definition}[section]
\newtheorem{rem}{Remark}[section]

\begin{document}
\title{ Many Three Dimensional Objects Inspired From Finite Groups }
\author{Zhi Chen}
\date{}
\maketitle

\begin{abstract}
Starting from considering deeper relationship between conjugacy classes and irreducible representations of  a finite group $G$, we find some quite simple $R-$matrice defined by using finite groups. This construction produces many sets (or topological spaces) admitting braid group actions. We introduce conceptions "extended $R-$matrix" and "generalized extended $R-$matrix" generalizing Turaev's enhanced $R-$matrix, which can still give invariants of oriented links. With these new frames, we show that above $R-$matrix, together with certain commuting pairs (essentially conjugacy classes of commuting pairs ) of $G$ can give integer invariants of oriented links. We construct some group dominating these integer invariants and prove these groups are link invariant by themselves. Using the language of the (colored) tangle category, we extended above invariant to invariant of  links and ribbon links colored by commuting pairs. We show given a oriented link diagram $L$, a suitable weighted sum of above invariant on all kinds of  coloring of $L$ (by conjugacy classes of $G$) is invariant under both two types of Kirby moves, thus giving a invariant for closed three manifolds. We define a group dominating those invariants, and prove this group is a three manifold invariant by itself.

\end{abstract}

\section{introduction}

For a finite group $G$, denote the set of  irreducible complex linear representations of $G$ as $\mathscr{R}(G) $, and  the set of conjugacy classes of $G$ as $ \mathscr{C}(G) $. It is well known that the number of  $\mathscr{R}(G)$ equals the number of   $\mathscr{C}(G)$   \cite{Se}. For some finite groups, for example, the permutation groups $S_n$, and finite groups in $SL(2,\mathbb{C})$, there even exist natural bijections between the sets $\mathscr{R}(G)$ and $\mathscr{C}(G)$.  There are also interesting relations between irreducible representations and conjugacy classes in the Langlands duality. Our original motivation is to search for possible further relations between irreducible representations and conjugacy classes. In the proof of some kind of "Modular properties" of the central algebra $\mathbb{C}G$ we found the following maps which are key objects of this paper.
\begin{align*}
&\psi _{g_1 ,g_2} :G\times G \rightarrow G\times G : (x,y )\mapsto (yx^{-1}g_1 ^{-1}x, x  ) ,\\
&\psi ^{-} _{g_1 ,g_2} :G\times G \rightarrow G\times G : (x,y )\mapsto (y, xy^{-1} g_2 y ) .
\end{align*}

The complex linear space spanned by the set $G\times G$ is naturally identified with $\mathbb{C}G\otimes \mathbb{C}G$. So $\psi _{g_1 ,g_2}$ ($\psi^{-} _{g_1 ,g_2}  $) induces a linear automorphism of $\mathbb{C}G\otimes \mathbb{C}G$, which is still denoted as $\psi_{g_1 ,g_2}$ ($\psi^{-} _{g_1 ,g_2} $).  We observed these
maps satisfy braid relations on the space $ \mathbb{C}G\otimes \mathbb{C}G \otimes \mathbb{C}G $ (Theorem \ref{thm:basicbraidrelation})  .
\begin{align*}
&( \psi_{g_2 ,g_3} \otimes id_{\mathbb{C}G} )( id_{\mathbb{C}G}\otimes \psi_{g_1 ,g_3}  )(  \psi_{g_1 ,g_2} \otimes id_{\mathbb{C}G}  ) \\
=&(id_{\mathbb{C}G}\otimes \psi_{g_1 ,g_2}    )( \psi_{g_1 ,g_3} \otimes id_{\mathbb{C}G} )( id_{\mathbb{C}G}\otimes \psi_{g_2 ,g_3}    ) .
\end{align*}

Besides, there are
$$\psi^{-} _{g_2 ,g_1} \psi_{g_1 ,g_2} = id_{\mathbb{C}G \otimes \mathbb{C}G}  ,\  and\  \psi_{g_2 ,g_1} \psi_{g_1 ,g_2} \neq  id_{\mathbb{C}G \otimes \mathbb{C}G} . $$

So these maps indeed define nontrivial actions of $G-$colored braids on tensor products of $\mathbb{C}G$. Especially, if we take $g_1 = g_2 = g$, then the map
$\psi_{g,g}$ is a nontrivial $R-$matrix.

  Now, on the side of irreducible representations, let $R_{G}$ be the complex linear space spanned by $\mathscr{R}(G)$. There is a natural commutative product on $R_G$ which make it into the so called representation ring. We emphasis the viewpoint that $R_{G}$ is a based algebra ( Definition \ref{defi:basedalgebra} ), that is, an algebra with a prescribed basis. In this case, the basis is the set of irreducible representations. As a based algebra, $R_G$ satisfies the "modular property" described in Theorem \ref{thm:repmodular}.

 On the side of conjugacy classes, let $\mathscr{C}_{G}$ be the complex linear space spanned by $\mathscr{C}(G) $. We defined a product on $\mathscr{C}_{G}$ which make it into a commutative algebra also (Theorem \ref{thm:centerproduct} ). It is a based  algebra whose prescribed basis is the set of conjugacy classes. Our first observation is that $\mathscr{C}_{G}$ also satisfies similar modular properties (Theorem \ref{thm:centermodular}).

 Let $g_1 , g_2 ,...,g_k$ be a series of elements in $G$. In above construction we introduced the follow set:
$$\mathscr{D}_{g_1 , g_2 ,...,g_k } =  \{ ( a_1 ,a_2 ,..., a_k   ) \in G^k | g_{1 } a_1 g_{2 } a_2 ...g_{k } a_k =e ; a_1 a_2 ...a_k =e \}.$$

 A part of the modular properties of $\mathscr{C}_{G}$ is (2) of Theorem \ref{thm:centermodular},   that is,  for any permutation $\sigma \in S_{k}$,
 $| \mathscr{D}_{g_1 , g_2 ,...,g_k }  | =|  \mathscr{D}_{g_{\sigma(1)} , g_{\sigma(2)} ,...,g_{\sigma(k)} }  | $. For a proof we found a bijection
 \begin{align*}
  \phi_{v,v+1} :\mathscr{D}_{g_1 ,...g_v , g_{v+1},...,g_k } &\rightarrow  \mathscr{D}_{g_1 ,...g_{v+1} , g_{v},...,g_k }  \\
     (a_1 , a_2 , ..., a_k ) &\mapsto  ( a_1 , ..., a_{v-2} , a_{v-1} g_{i_v} a_v , a_v ^{-1} g_{i_v } ^{-1} , a_v a_{v+1} , a_{v+2}, ..., a_k    ) .
 \end{align*}

  In Theorem \ref{thm:satisfybraid},we observe the maps $\phi_{v,v+1}$ satisfy nontrivial braid relations. So through the maps $\phi_{v,v+1}$, finite groups are connected to 3-topology. For above construction of the sets $\mathscr{D}_{g_1 ,..., g_k}$ and  the maps $\phi_{v,v+1}$, $G$ need not be a finite group. Actually $G$ can be any group. If we take a topological group $G$ then the sets $\mathscr{D}_{g_1 ,g_2 ,...,g_k } $ are topological spaces.  So  this construction in fact produce many sets admitting braid group actions. But it is not the emphases to study those braid group actions.

  It is natural to ask if those braid group actions could lead to Link invariants and 3-manifold invariants. We can do nothing with an braid group action in that form as $\phi_{v,v+1}$.   In Theorem \ref{thm:identifybraidaction}, we prove that above braid group actions are essentially equivalent to braid group actions induced by the following simpler maps introduced in Definition \ref{defi:refinedbraidaction}:
\begin{align*}
 \psi_{i;g_1,...,g_n }: \mathscr{G}_{g_1 ,..., g_i ,g_{i+1} ,..., g_n } &\rightarrow   \mathscr{G}_{g_1 ,..., g_{i+1},g_{i} ,..., g_n } \\
(b_1, ..., b_i ,b_{i+1} ,..., b_n ) &\rightarrow (b_1 ,...,b_{i+1} b_i ^{-1} g_i ^{-1} b_i ,b_i ,..., b_n  ) .
\end{align*}

Where we set $ \mathscr{G}_{g_1 ,..., g_i ,g_{i+1} ,..., g_n }=\underbrace{ G\times G\times ...\times G}_{n}$ for any sequence $g_1 ,g_2 ,..., g_n \in G$, whose $i-$th component can be understood as "a set $G$ colored by an element $g_i$". And the map $\psi_{i:g_1 ,..., g_n }$ are actually associated to certain "$G-$colored braid".

For above operator $\psi_{i;g_1,...,g_n }$ if we let $n=2$ and $g_1 = g_2 =g$, then the linearized map $\psi _{1,g,g} :\mathbb{C} G \otimes \mathbb{C} G \rightarrow  \mathbb{C} G \otimes \mathbb{C} G $ ( we simply denote it as $\psi_{g,g}$ )  is a R-matrix, something closer to link invariants.

In \cite{Tu1}, Turaev introduced the notion "enhanced R-matrices" ,and use them to construct link invariants from those R-matrices from quantum groups.  After some attempts we realized the R-matrix $\psi _{g,g}$ possibly can not "be enhanced" to give a enhanced R-matrix (thus give a link invariant). Instead we find the notion enhanced R-matrix can be relaxed to "extended R-matrix " (Definition \ref{defi:extented}) $(I , f , c)$ ,where $I \in End(V\otimes V)$ is a R-matrix, and $f,c\in End(V)$, they satisfy some conditions together. In Theorem \ref{thm:extendedgivelink} we prove any extended R-matrix can still give a link invariant. The key notion for this construction is the "modified R-matrix " introduced in Lemma \ref{lem:modifiedRmatrix}.

 In Definition \ref{defi:extendedpair} , we introduced the following conception " extended pair " for a group:

\begin{defi}
\begin{enumerate}
\item[(1)] We call $( g,b )\in G\times G $ as a extended pair for $G$, if  $gb=bg$ and $b^2 =1$. We denote the set of extended pairs for $G$ as $\mathscr{E}G $.
\item[(2)] We call $( g,b )\in G\times G $ as a commuting pair for $G$, if  $gb=bg$. We denote the set of commuting pairs for $G$ as $\mathscr{C}G $.
\end{enumerate}
\end{defi}

Denote the set of extended pairs for $G$ as $\mathscr{E}G$. Notice each $g\in G$ can give a extended pair $(g,e)$, where $e$ is the unit. In this way $G$ can be looked as a subset of $\mathscr{E}G$.

 In Theorem \ref{thm:extendableGRmatrix},we prove that the R-matrix $\psi_{g,g}$ along with any extended pair $(g,b)$ can give a extended R-matrix in some way. Thus,through Theorem \ref{thm:extendedgivelink},  for any finite group $G$, and extended pair $(g,b)$, we have a integer link invariant $\Lambda_{G;(g,b)} (-) $. We proved if a extended pair $(g_1 ,b_1 )$ is "conjugated" to another one $(g_1 , b_1)$, that is, there exist $h\in G$ such that
$g_2 = hg_1 h^{-1}  , b_2 = h b_1 h^{-1}  $, then $\Lambda _{G; (g_1 , b_1) }= \Lambda_{G; (g_2 , b_2 )}$.

Above construction of link invariants can be extended in several aspects. Still let $G$ be a finite group. A oriented (framed) $\mathscr{E}G$-tangle is a oriented (framed) tangle with every component associated with a element of $\mathscr{E}G$. Denote the category of $\mathscr{E}G$-colored tangles ( $\mathscr{E}G$-colored framed tangles  )as $\mathcal{T}_{\mathscr{E}G}$   ( $ \mathbb{T}_{\mathscr{E}G} $ ) . By using Turaev's construction for tangle invariants, in Theorem \ref{thm:extended colored link invariant} , we prove that certain  modification of the map $\psi_{g_1 ,g_2}$ can give  a functor  $F: \mathcal{T}_{\mathscr{E}G} \rightarrow \mathcal{V}  $  in Theorem \ref{thm:extended colored link invariant} for $\mathscr{E}G$-tangles with some additional datum,   where $\mathcal{V} $ is the category of complex linear spaces.   This tangle invariant naturally contain the link invariant $\Lambda_{G;(g,b)}$   as special cases.

 Similarly, in Theorem  \ref{thm:extended colored framed link invariant} we  define a functor  $\mathbb{F}: \mathbb{T}_{\mathscr{E}G} \rightarrow \mathcal{V} $, by using the unmodified map  $\psi_{g_1 ,g_2} $ and some similar additional datum. If restrict this functor to the set of  ribbon links, we obtain many invariants for ribbon links. We also prove that if change the colour $(g,b)$  of a string (or a component ) of a $\mathscr{E}G- $colored tangle to another extended pair $(hgh^{-1} , hbh^{-1}  )$ conjugated to it, then the resulted invariant is essentially the same. We we can think that for the functors $ F$ and $\mathbb{F}$, the colored links are actually colored by conjugacy classes of extended pairs.

 In Theorem \ref{thm:3mfdintegerinvariant} we prove that given  ribbon $L=L_1 \sqcup ...\sqcup L_m $ with $m$ components, a suitable weighted  sum of  $ \mathbb{F}( L_{d_1 , ..., d_m } ) $ is invariant under both Kirby I and Kirby II move, where $L_{d_1 ,..., d_m} $ means the $\bar{G}-$ colored link obtained  by associating the $d_i -$th conjugacy class of $G$ to the component $L_i$. Thus this sum is an invariant $\mathscr{F}(M_{L}) $ of the 3 manifold $M_{L}$ obtained by doing surgery along the ribbon link $L$.  This formula parallel the WRT 3-mfd invariant very nicely, where the important conception " irreducible representation of quantum group " and " quantum dimensional of a irreducible representation " are replace by "conjugacy class of $G$" and "the number of elements in a conjugacy class" respectively.

 We have constructed many integer invariants for links and 3 manifolds. We observe that there are several groups behind these integer invariants, dominating them actually with some additional datum. In Definition \ref{defi:groupinvariantbraidversion}, we define a group $\mathcal{G}_{L,b}$
 and a group $\mathbb{G}_{L,b}$. The distinct place of the group $\mathcal{G}_{L,b}$  is that the value of the functor $F$ on a extended pair $(g,e)$ equal the number of morphisms from $\mathcal{G}_{L,b}$ to the given finite group $G$ satisfying certain relations ( Proposition \ref{prop:groupinandnumberin}). Similar result holds between the group  $\mathbb{G}_{L,b} $ and the functor $\mathbb{F}$.

 In Theorem \ref{thm:groupinvariant1} we prove that the group $\mathcal{G}_{L,b} $ has no relation with the braid $b$ and is a link invariant itself. In Proposition \ref{prop:graphpresentation} we give another presentation of the group $\mathcal{G}_{L,b}$ based on a link diagram of $L$. In proposition \ref{prop:framedgroupinvariant} and Proposition \ref{prop:simplifiedpresentation} we give a presentation of the group $\mathbb{G}_{L,b}$ based on a given link diagram of $L$.

 In section 9, we construct another  link group invariant $\bar{\mathcal{G} }_{L}$ that dominating those integer invariants $\Lambda_{G; (g,b)}$ by using extended pairs of $G$, and have $\mathcal{G}_{L,b} $ as a quotient group. In section 10, we prove that the group $\mathbb{G}_{L}$ ($\mathbb{G}_{L,b} $  ) is invariant under the second type of Kirby moves, and change simply under the first type of Kirby moves. From the group $\mathbb{G}_{L}$ we derive a group $\hat{ \mathbb{G}}_{L}$ which is invariant under both types of Kirby moves, so be a group invariant of the 3 manifold $M_{L}$.

 In section 10, we introduce "generalized extended $R-$matrix " (Definition \ref{defi:generalizedextended}  ) and its special case "special generalized extended $R-$matrix " (Definition \ref{defi:specialgeneralizedextended}  ), which have the extended $R-$matrix as special case. We show in Theorem \ref{thm:linkinvariantfromgeneralizedextended} that the generalized extended $R-$matrix still give link invariants, in a slightly more complicated way. On this more general setting, we show for every finite group $G$, and for any commuting pair $(g,b)$ of $G$ (not just extended pairs ), we can build a generalized extended $R-$matrix, so give an invariant for oriented links. In Definition \ref{defi:groupinvariantbraidversion3}, we define another group invariant of links that dominating above integer invariants from commuting pairs, and stronger than the group invariant we defined in Section 9.

 We will discuss the relationship between above group invariants and known group invariants of links and 3 manifolds in later works.

 At last we present several necessary conceptions including the R-matrix due to Drinfeld, and the enhanced R-matrix due to Turaev.

\begin{defi}
The braid group $B_n$ is the group with generating set $\{ \sigma_1 ,\sigma_2 ,...,\sigma_{n-1} \} $  submitting to following relations.

(1) $\sigma _i \sigma_j = \sigma_j \sigma_i $ for $|i-j| \geq 2$,

(2) $\sigma_i \sigma_{i+1} \sigma_i = \sigma_{i+1} \sigma_i \sigma_{i+1} $.
\end{defi}

The following conception of $R$-matrix is closely related with braid groups.

\begin{defi}
Let $V$ be a finite dimensional vector space. A $R$-matrix on $V\otimes V$ is a invertible operator $R \in End( V\otimes V) $ such that the following Yang-Baxter equation is satisfied. Where both sides belongs to $End (V\otimes V\otimes V )$.

$( R \otimes  id  ) (id \otimes R ) (R\otimes id )=(id \otimes R ) (R\otimes id ) (id \otimes R )  $.
\end{defi}

If $R $ is a $R$-matrix on $V\otimes V$, suppose $R=\sum _{l} x_l \otimes y_l \in End(V) \otimes End(V)$. In
$End( V^{ \otimes n } )$ and for $1\leq i<j\leq n$, we set $R_{ij} =\sum_{l} id \otimes ...\otimes x_l \otimes ...\otimes y_l \otimes ...\otimes id $ where the terms "$x_l$" and "$y_l$" are in  the $i-$th and the $j-$th position respectively, other components are all $id$'s. Then it is evident that the map $\sigma_i \mapsto R_{i,i+1} $ extends to a $B_n$ representation on $End( V^{\otimes n} )$.




In \cite{Tu1}, Turaev introduced the following conception of enhanced $R-$matrice.

\begin{defi}
\label{defi:enhanced}
An enhanced $R-$matrix consists of the datum: $I \in End(V\otimes V )$, $f\in End(V) $, $\lambda ,\mu \in \mathcal {C}$, satisfying the following relations.

(1) $I$ is a $R-$matrix.

(2) $f  $ is  invertible and $(f\otimes f )\circ I=I\circ (f\otimes f )  $.

(3) $tr _2 ( I \circ (f\otimes f)) = \lambda \mu f$; $tr _2 (I^{-1} \circ (f\otimes f)  )= \lambda ^{-1} \mu f$.
\end{defi}

It is known that all $R-$matrice obtained from quantum groups as above can be "enhanced" by certain $f$ to because enhanced $R-$matrice.  An enhanced $R-$matrix leads to a link invariant as follows.

For any positive integer $n$, the $R-$matrix $I$ induces a $B_n $ representation $\rho_n$ on $V^{\otimes n} $.   Suppose $\beta \in B_n$ is any braid. Denote the closure of $\beta $(an oriented link) as $\hat{\beta}$, and denote the braid index of $\beta$ (the number of positive generators in $\beta$ minus the number of negative generators in $\beta$  ) as $\alpha (\beta)$ Then set

$$ \rho (\beta ) = \lambda ^{ \alpha(\beta) } \mu ^{-n} tr (\rho _n (\beta) \circ f^{ \otimes n} ) . $$

It isn't hard to show above $\rho(\beta)$ is invariant under both two Markov moves, thus defines a link invariant.

In above definition, if we set $\bar{f}= \frac{1}{\mu} f$, then we have

$$ tr_2 (I\circ (\bar{f} \otimes \bar{f} ) ) = \lambda \bar{f}  ;\ tr_2 (I^{-1} \circ (\bar{f} \otimes \bar{f} ) ) = \lambda ^{-1} \bar{f} .$$

So the data $(I , \bar{f} , \lambda , \mu =1)$ is still an enhanced $R-$matrix, and it is easy to see the resulted link invariant is the same. So in definition 1.3 we can set $\mu =1$ without any loss.

Since partial trace will be used in several occasions, we give its definition and main properties here for later use.

\begin{defi}
Suppose $V$ is a linear space with a basis $\{ v_1 ,..., v_n \}$. Let $f\in End(V^{\otimes n})$ such that
$f ( v_{i_1 } \otimes ...\otimes v_{i_n } ) = \sum_{j_1 ,..., j_n  }  \alpha _{i_1 ,..., i_n} ^{j_1 ,..., j_n } v_{j_1 ,..., j_n }.$  Then $tr_{n} (f) \in End(V^{\otimes (n-1)})$ is defined as
$ tr_{n} (v_{i_1} \otimes ...\otimes v_{i_{n-1}}) = \sum_{j_1 ,..., j_{n-1}} (\sum _{i=1} ^{n} \alpha _{i_1 ,..., i_{n-1} , i} ^{j_1 ,..., j_{n-1} ,i} ) v_{j_1 ,..., j_{n-1}}$. The definition is well defined,has no relation with the choice of basis of $V$.
\end{defi}

\begin{lem}
\label{lem:partialtrace}
Let $V$ be a linear space, $f_1 ,f_2 \in End(V^{\otimes (n-1)}) $ and $g\in End(V^{\otimes n})$, then

$ tr_n ( (f_2 \otimes id_V ) \circ g \circ (f_1 \otimes id_V ) ) = f_2 \circ tr_2 (g) \circ f_1 . $
\end{lem}

\section{Center and the representation ring}
  Let $G$ be a finite group, let $\{ [g_i ] \} _{i=0,...,N}$ be the set of conjugacy classes in $G$. Where $g_i \in G$, and $[g_i ]$ denotes the conjugacy class of $g_i$. And we suppose $g_0 =e$.  We set $\mathscr{C}_G = \mathbb{C}< [g_i ] >_{i=1,...,N}$ to be the linear space with $\{ [g_i ] \} _{i=0,...,N}$ as a basis. We will use the notation $\mathbb{C}<v_i >_{i=0,1,...,N}$ to denote the complex linear space with a basis $\{  v_i \}_{i=0,1,..., N}$ in many occasions.

\begin{thm}
\label{thm:centerproduct}
On $\mathscr{C}_G $ we define the following product
$  [g_i ] \centerdot  [g_j ] =  \sum _{h\in G} [ g_i h g_j h^{-1} ], $
then $ \mathscr{C}_G$ become a commutative and associative algebra.

\end{thm}

\begin{pf}
For commutativity,
$$[g_i ] \centerdot  [g_j] =  \sum _{h\in G} [ g_i h g_j h^{-1} ]=   \sum _{h\in G} [ g_j h^{-1} g_i h ]=  \sum _{h\in G} [ g_j h g_i h^{-1} ] = [g_j ][g_i] .$$

For associativity,
\begin{align*}
([g_i] [g_j] )[g_k] &= \sum _{h_1 , h_2 \in G  } [ g_i h_1 g_j h_1 ^{-1} h_2 g_3 h_2 ^{-1}   ]
 =  \sum _{h_1 , h_2 \in G  } [  g_i h_1  g_j h_1 ^{-1} h_2 g_3  (h_1 ^{-1} h_2 )^{-1}  h_1 ^{-1}  ] \\
 & =  \sum _{h_1 , h_3 \in G  } [ g_i h_1 g_j h_3 g_3 h_3 ^{-1}  h_1 ^{-1}  ]
 = [g_i ] ([g_j ] [g_k ]).
 \end{align*}
\end{pf}

\begin{defi}
\label{defi:basedalgebra}

A based algebra is a algebra $A$ together with a prescribed basis of $A$.
\end{defi}

Later we will show $\mathscr{C}_G$ is isomorphic to the center of the group algebra $ \mathbb{C} G $. So $\mathscr{C}_G$ is not something new. The main viewpoint of this paper is to study $\mathscr{C}_G$ as a based algebra, with the natural basis formed by conjugacy classes. Before that we recall knowledge about representative rings £¬ which are typical examples of based algebras.

Suppose the set of irreducible complex representations of $G$ is $\{ \rho _i \} _{i=0,...,N}$, where we let $\rho _0$ be the trivial representation.  The representation ring $R_G$ is $\mathbb{C} < \rho _i  >_{i=0,...,N}$ as a linear space. The product in $R_G$ is defined by tensor product of representations. With the natural basis $\{ \rho _i \} _{i=0,1,...,N}$ $R_G$ become a based algebra. We define two sets of constants as follows.
\begin{enumerate}
\item[(1)] $\alpha _{i,j} ^{k} : $  $ \rho _i \otimes \rho _j = \sum_{k} \alpha _{i,j} ^{k} \rho _k $;

\item[(2)] $A_{i_1 ,i_2 ,..., i_k } :$ $A_{i_1 ,i_2 ,...,i_k } = \dim Hom_{G} (\rho _0 , \rho _{i_1 } \otimes \rho _{i_2 } \otimes ... \otimes \rho _{i_k } ) $.
\end{enumerate}
There is a natural involution on the set $\{ 0,1,2,...,N \}$: $i\mapsto i^* $ such that $\rho _{i^*  }$ is the dual of $\rho _i $. Evidently $0^* =0 $.  The following theorem summarize some features of these numbers.

\begin{thm}
\label{thm:repmodular}
\begin{enumerate}

\item[(1)] $ \alpha _{i,j} ^{k} = A_{i,j,k^* } $;

\item[(2)] $A_{i_1 , i_2 ,..., i_k } = A_{i_{\sigma (1)} , i_{\sigma (2)} ,..., i_{\sigma (k)}} $ for any $\sigma \in S_k $(the permutation group ).

\item [(3)] $\rho _{i_1 } \rho _{i_2 } ... \rho _{i_l } = \sum _{k} A_{k^* , i_1 , i_2 ,..., i_l  } \rho _k  $;

\item[(4)] $A_{i,j} = 0 $ if $j \neq i^* $ and $A_{i,j} = 1 $ if $j = i^* $.

\item [ (5)] $A_{i_1 , i_2 ,..., i_k } = \sum _{j=0} ^{N} A_{i_1 ,..., i_v ,j } A_{j^* , i_{v+1} ,..., i_k }$ for any $1< v < k$.
\end{enumerate}

\end{thm}

\begin{pf}
By definition,
$$\alpha _{i,j} ^{k}= \dim Hom_{G} ( \rho_k , \rho_i \otimes \rho_j ) = \dim Hom_{G} (\rho_0 , \rho _i \otimes \rho _j \otimes \rho _{k^* }) = A_{i,j, k^* }.$$
 So we have the statement (1). The statement (2) is true because $R_G$ is a commutative ring.

 The statement (3) follows from
$$ \dim Hom _G ( \rho _k ,  \rho _{i_1 } \otimes \rho _{i_2 } \otimes ... \otimes \rho _{i_l } ) = \dim Hom _G ( \rho _0 , \rho _k \otimes \rho _{i_1 } \otimes \rho _{i_2 } \otimes ... \otimes \rho _{i_l } ) = A_{k^* , i_1 , i_2 ,..., i_l } .$$

The statement (4) is because
 $$ A_{i,j} =\dim Hom_G ( \rho _0 , \rho _i \otimes \rho _j  )= \dim Hom _G ( \rho _{i^* } , \rho _j  ).$$

For (5),
 \begin{align*}
(\rho _{i_1 } \otimes \rho _{i_2 } \otimes ...\otimes \rho _{i_v }) \otimes ( \rho _{v+1 } \otimes ...  \otimes \rho _{i_k })  &= (  \oplus _{j_1 }  A_{i_1 , i_2 ,..., i_v , j_1 ^*  } \rho _{j_1  }   ) \otimes ( \oplus _{j_2 } A_{j_2 ^* , i_{v+1} , i_{v+2 } , ..., i_k  } \rho _{j_2 }  ) \\
   &= \oplus _{j_1 , j_2  }  A_{i_1 , i_2 ,..., i_v ,j_1 } A_{j_2 , i_{v+1} ,..., i_k }  \rho _{j_1 } \otimes \rho _{j_2 }.
 \end{align*}

so

\begin{align*}
A_{ i_1 , i_2 , ..., i_k  } &= \dim Hom_G ( \rho _0 , \rho _{i_1 } \otimes \rho _{i_2 } \otimes ...\otimes \rho _{i_k } ) = \sum _{j_1 , j_2 } A_{i_1 , i_2 , ..., i_v , j_1 } A_{j_2 ,i_{v+1 } , i_{v+2 } , ..., i_k } A_{j_1 , j_2 } \\
 &= \sum _{j } A_{i_1 , i_2 , ..., i_v , j } A_{j^* ,i_{v+1 } , i_{v+2 } , ..., i_k }  .
\end{align*}

The last equality sign is by (4).

\end{pf}

\section{ $\mathscr{C}_G $ as a based algebra }

First, for the based algebra $\mathscr{C}_G $, there is also a natural involution on the set of indices $\{ 0,1,2,..., N \} $: $i \mapsto i ^{   \star  }   $ such that $[g_i ^{-1}] = [g_{i^{   \star    } }]$. We define a kind of set important for this paper as follows. For convenience, we denote the set
$\underbrace{G\times G\times ...\times G}_{n} $ simply as $G^{n}$.

\begin{defi}
\label{defi:keyset}
Let $G$ be a group. Suppose $\{ g_1 , g_2 ,..., g_n \} \subset G$. We set

$\mathscr{D}_{g_1 , g_2 ,...,g_n } =  \{ ( a_1 ,a_2 ,..., a_n   ) \in G^n | g_{1 } a_1 g_{2 } a_2 ...g_{n } a_n =e ; a_1 a_2 ...a_n =e \} $.

\end{defi}

Then we set some constants and other  structures needed later as follows.

\begin{enumerate}
\item[(1)] $\bar {B} _{i_1, i_2 , ..., i_k }$: For $0\leq i_1 ,i_2 ,..., i_k \leq N $, Choose any representative $g_{i_{\alpha }} $ from $[ g_{i_{\alpha }}]$.  Let

  $\bar {B} _{i_1 , i_2 ,..., i_k } =| \mathscr{D}_{ g_{ i_1 } ,g_{ i_2 } ,..., g_{ i_k }  }  |$;

\item[(2)] $\bar {\beta } _{i,j} ^k $: these constants are defined by  $[g_i ] [g_j ] = \sum _{k} \bar { \beta } _{i,j} ^k [g_k ]$;

\item[(3)] Let $C_{g_k }$ be the centralizer of $g_k $, that is , the subgroup consists of elements commutes with $g_k $. And we set $v_k = |C_{g_k } |$. It is easy to see $v_k $ doesn't depend on the choice of representatives $g_k $ in $[g_k ]$ ,so they are well defined . We set $w_k = |G| / v_k =|[g_k ]| $.
\end{enumerate}

\begin{thm}
\label{thm:centermodular}
\begin{enumerate}
\item[(1)] The constants $\bar {B} _{i_1 , i_2 ,..., i_k }$ don't depend on the choice of representatives $g_{i_{\alpha }}\in [ g_{i_{\alpha }}]$, so they are well defined.

\item[(2)] $\bar {B} _{i_1 , i_2 ,..., i_k } = \bar { B} _{i_{\sigma (1)} , i_{\sigma (2)} ,..., i_{\sigma (k)} }$ for any $\sigma \in S_k $.

\item[(3)] $\bar {B} _{i,j } =0 $ if $j \neq i^{\star } $ and $\bar {B} _{i, i^{\star }} =v_i $.

\item[(4)] $\bar { \beta } _{i,j} ^{k} = \frac{1}{v_k } \bar {B} _{i,j, k^{\star }}$.

\item[(5)] $\bar {B} _{i_1 , i_2 ,..., i_k } = \sum _{j=0} ^N  \frac{1}{v_j }  \bar {B} _{i_1 , i_2 , ..., i_{v-1} ,j } \bar {B} _{j^{\star } , i_{v} , i_{v+1} , ..., i_k } $.
\end{enumerate}

\end{thm}

\begin{pf}
Statement (1): For some fixed $v : 1\leq v\leq k$,   If choose another representative $g_{ i_v } ^{'} = h g_{ i_v } h^{-1 } \in [ g_{ i_v } ] $, then there is a bijective map
\begin{align*}
 \phi : \mathscr{D}_{g_{i_1} , ..., g_{i_k }} & \rightarrow \mathscr{D}_{ g_{i_1 } ,..., g_{i_v} ^{'} ,...,g_{i_k }  }\\
 by\  (a_1 , a_2 ,..., a_k ) & \mapsto ( a_1 , a_2 ,..., a_{v-1} h^{-1} , h a_v ,..., a_k   ).
\end{align*}

 So the number $| \mathscr{D}_{g_{i_1} , ..., g_{i_k }} |$ is unchanged if we change one representative, then also unchanged if we change all representatives.

 Statement (2): We only need to prove for $1\leq v \leq k-1 $, $\bar{B}_{i_1 ,..., i_k } = \bar{B}_{i_1 ,..., i_{v-1} , i_{v+1} i_v , i_{v+2} , ..., i_k }$. It is proved by construction of the following bijective map $\phi _{v,v+1} :$
\begin{align*}
  \mathscr{D}_{ g_{i_1 } , ..., g_{i_k }} & \rightarrow \mathscr{D}_{ g_{i_1 } , ..., g_{i_{v+1}} , g_{i_v} ,..., g_{i_k }}  \\
(a_1 , a_2 , ..., a_k ) & \mapsto ( a_1 , ..., a_{v-2} , a_{v-1} g_{i_v} a_v , a_v ^{-1} g_{i_v } ^{-1} , a_v a_{v+1} , a_{v+2}, ..., a_k    ) .
\end{align*}

 Statement (3) is evident.

 Now we prove (4). By Definition 1.1,
 $$\bar{\beta }_{i,j} ^{k} =| \{ h\in G |  g_i h g_j h^{-1} = a g_k a^{-1} , for\ some\ a\in G \}   |.$$
 Since  $g_i h g_j h^{-1} = a g_k a^{-1}$ is equivalent to $g_i h g_j h^{-1} a g_k ^{-1} a^{-1} $, and notice  $ h \cdot (h^{-1} a) \cdot a^{-1} =e  $,  we consider the following  map $\phi$ :
\begin{align*}
  \mathscr{D}_{g_i , g_j  , g_{k^{\star }} } & \rightarrow    \{ h\in G |  g_i h g_j h^{-1} = a g_k a^{-1} , for\ some\ a\in G \}  \\
 (a_1 , a_2 , a_3 ) & \mapsto a_1   .
\end{align*}

It is easy to see $\phi$ is surjective, and for any $h$ in the left side, $\phi ^{-1} (h) $ consists of $v_k =|C_{g_k }|$ elements. So we have $ \bar{B} _{i,j , k^{\star }} = v_k \bar{\beta }_{i,j} ^{k} $ and (4) is proved.

For (5), we observe that if
$ g_{i_1} x_1 g_{i_2 } x_2 ...g_{i_{v-1} } x_{v-1} g_j  y_1 =e $ and $x_1 x_2 ...x_{v-1} y_1 =e $;

 $g_{j}^{-1} y_2 g_{i_{v}} x_v g_{i_{v+1}} ...g_{i_k } x_k =e $ and $y_2 x_v x_{v+1} ...x_k =e$, then

$g_{i_1 } x_1 g_{i_2 } x_2 ...g_{i_{v-1}} x_{v-1} y_2 g_{i_{v } } x_v ...g_{i_k } x_k y_1 =e $ and $x_1 x_2 ...x_{v-1} y_2 x_v ...x_k y_1 =e$.

It enable us to define the following map $\Phi$ :
\begin{align*}
 \mathscr{D}_{i_1 , i_2 , ..., i_{v-1} ,j } \times \mathscr{D}_{j^{\star } , i_v ,...,i_k } & \rightarrow
  \mathscr{D}_{i_1 , i_2 ,..., i_k } \\
[  (x_1 , x_2 , ..., x_{v-1} , y_1 ) , (y_2 , x_{v}, x_{v+1} ,..., x_k )] & \mapsto (x_1 , x_2 , ..., x_{v-2} , x_{v-1} y_2 , x_v , ..., x_k y_1 ) .
\end{align*}

Suppose $(x_1 , x_2 , ..., x_{v-1} , x_v ,..., x_k ) \in  \mathscr{D}_{i_1 , i_2 ,..., i_k } $,
Suppose $[ x_1 x_2 ... x_{v-1} g_{i_v } x_v ...g_{i_k } x_k ] = [g_j ]  $, choose any $y_1 $ such that
$ x_1 x_2 ... x_{v-1} g_{i_v } x_v ...g_{i_k } x_k = y_1 ^{-1} g_j y_1    $. And we set $y_2 = y_1 x_1 x_2 ...x_{v-1}$.
Then \begin{align*}
 ( x_1 , ..., x_{v-2}, x_{v-2}^{-1} x_{v-3}^{-1} ...x_1 ^{-1} y_1 ^{-1}, y_1  ) & \in \mathscr{D}_{i_1 , i_2 ,..., i_{v-1} ,j }  , \\
 (y_2 , x_v , x_{v+1} ,..., x_k y_1 ^{-1}) & \in \mathscr{D}_{j^{\star} , i_{v+1} ,..., i_k }, \  and \\
 \Phi ([   ( x_1 , ..., x_{v-2}, x_{v-2}^{-1} x_{v-3}^{-1} ...x_1 ^{-1} y_1 ^{-1}, y_1 ),& (y_2 , x_v , x_{v+1} ,..., x_k y_1 ^{-1}) ])  = (x_1 , x_2 , ..., x_{v-1} , x_v ,..., x_k).
\end{align*}

Then the map $\Phi$ is surjective, and it isn't hard to see that every element in

 $\Phi ^{-1} (  (x_1 , x_2 , ..., x_{v-1} , x_v ,..., x_k)  )$ is constructed in this way by reversing the process. In above process of constructing an element in  $\Phi ^{-1} (  (x_1 , x_2 , ..., x_{v-1} , x_v ,..., x_k)  )$ , we have a freedom of the choice of $y_1 $
Since there are altogether $|C_{g_j ^{-1} } | =| C_{g_j }  | = v_j $ kinds of way to choose $y_1 $, we have:

$| \Phi ^{-1} (  (x_1 , x_2 , ..., x_{v-1} , x_v ,..., x_k)  )  | = v_j $ ( $j $ is determined by $(x_1 , x_2 , ..., x_{v-1} , x_v ,..., x_k)$ )

So we have
$\bar{B} _{i_1 , i_2 ,..., i_k } =  \sum _{j=0} ^N  \frac{1}{v_j }  \bar {B} _{i_1 , i_2 , ..., i_{v-1} ,j } \bar {B} _{j^{\star } , i_{v} , i_{v+1} , ..., i_k } .$

\end{pf}

$\\$

The theorem 2.1 have a very similar form with theorem 1.1. To obtain a complete similar form we do the following "normalization". \begin{enumerate}
\item[(1)] Set $T_i = \frac{1}{\sqrt{v_i}  } [g_i ] $;
\item[(2)] Set $B_{i_1 ,i_2 ,..., i_k } = \frac{1}{\sqrt{v_{i_1} v_{i_2} ...v_{i_k }}   } \bar{B}_{i_1 ,i_2 ,..., i_k} $;
\item[(3)] Set $ \beta _{i,j} ^{k} = \frac { \sqrt{v_k }}{ \sqrt{v_i v_j } } \bar{\beta } _{i,j} ^{k} $.
\end{enumerate}

Then for the new based algebra $(\mathscr{C}_G , {T_0, T_1 , ..., T_N}  )$ we have

\begin{cor}
\begin{enumerate}
\item[(1)] $ T_i T_j = \sum _{k} \beta _{i,j} ^k T_k $ .
\item[(2)] $ B_{i_1 , i_2 ,..., i_k } =  B_{i_{\sigma (1)} , i_{\sigma (2)} ,..., i_{\sigma (k)} }$ for any $\sigma \in S_k $.
\item[(3)] $ B_{i,j } =0 $  if $ j \neq i^{\star } $ and $ B_{i, i^{\star }} =1$ .
\item[(4)] $ \beta  _{i,j} ^{k} =  B_{i,j, k^{\star }} $.
\item[(5)] $ B_{i_1 , i_2 ,..., i_k } = \sum _{j=0} ^N    B_{i_1 , i_2 , ..., i_{v-1} ,j }  B_{j^{\star } , i_{v} , i_{v+1} , ..., i_k } $.
\end{enumerate}
\end{cor}

We introduce the following definition:

\begin{defi}
A based algebra $( \mathscr{C} , \{v_i \} _{i=0} ^N  )$  ($N$ can be taken to be infinity) is called a modular algebra if
there is a set of constants $\{  A_{i_1 , ..., i_k }  | k\geq 1 ,  0\leq i_1 , ..., i_k \leq N   \} $ and there is a involution in the indices set $\{ 0,1,2,..., N\} $: $i\mapsto i^* $
such that: \begin{enumerate}
\item[(1)] The constants are symmetric: $ A_{i_1 , ..., i_k } =A_{i_{\sigma (1)} , i_{\sigma (2)} ,..., i_{\sigma (k)} }$ for any $\sigma \in S_k $.

\item[(2)] If $\beta _{i,j} ^k $ are the structure constants for $\mathscr{C}$, that is , $v_i v_j = \sum _{k} \beta _{i,j} ^k v_k $ , then $\beta _{i,j} ^{k} = A_{k^* ,i,j }$.

\item[(3)] For fixed $v$,$A_{i_1 , i_2 ,..., i_k } = \sum _{j=0} ^N    A_{i_1 , i_2 , ..., i_{v-1} ,j }  A_{j^{\star } , i_{v} , i_{v+1} , ..., i_k } $.
\end{enumerate}

\end{defi}

For a modular algebra $ ( \mathscr{C} , {v_i }_{i=0} ^N  ) $, if all the constants $ A_{i_1 , ..., i_k } $ are great than or equal $0$, we call it as a positive modular algebra. If all the constants  $A_{i_1 , ..., i_k } $ are integers, we call it as a integral modular algebra. Then, the representation ring $R_G $ of a finite group $G$  with the natural basis is a positive integral modular algebra. Yet the above algebra $\{ \mathscr{C}_G , \{  T_i \} _{i=0} ^N    \}$ is only a positive modular algebra. If we take $G$ as a compact Lie group and let $R_G$ be its representation ring, then it isn't hard to see $R_G$ with the natural basis is a  infinite dimensional positive integral modular algebra .

\section{ Further properties of $\mathscr{C}_G  $  }

 In this section we study some functional properties of the based algebra $\mathscr{C}_G  $.  We denote the group algebra of $G$ as $\mathbb{C} G$, which is a linear space with the set $G$ as a basis, whose multiplication rule  follow from the multiplication rule in $G$. Denote the center of $\mathbb{C} G$ as $C_G $.

 For any conjugacy class $[g_i ]$ in $G$, we set $c_i = \sum _{g\in [g_i ]} g  \in \mathbb{C} G$. It is well known that $c_i $ belongs to $C_G$ , and $\{  c_i \} _{i=0,...,N}$ is a basis of $C_G $.

 \begin{thm}
  There is a natural isomorphism from $\mathscr{C}_G  $  to  $C_G $.
 \end{thm}

 \begin{pf}  Suppose $c_i c_j = \sum _{k} \gamma _{i,j} ^{k}  c_k $. Set the following three set:
 \begin{enumerate}
 \item[(1)] $A= \{  ( \bar{g}_i , \bar{g}_j  ) \in G^2  |  \bar{g}_i \in [g_i ] , \bar{g}_j \in [g_j ] , \bar{g}_i \bar{g}_j \in [g_k ]    \} $;

 \item[(2)] $B= \{  (h_1 , h_2 , h_3 ) | h_1 g_i h_1 ^{-1} h_2 g_j h_2 ^{-1} h_3 g_k h_3 ^{-1} =e    \}$;

 \item[(3)] The set $C= \mathscr{D}_{i,j ,k^{\star }}$, which contains $\bar{B}_{i,j,k^{\star }}$ elements by our notation.
 \end{enumerate}

 First we see $\gamma _{i,j} ^{k} = |A|/ w_k = |A|v_k /|G|$.

 There is a natural map $\Phi : B\rightarrow C$, $(h_1 , h_2 ,h_3 ) \mapsto (h_1 ^{-1} h_2 , h_2 ^{-1} h_3 , h_3 ^{-1} h_1 )$. Evidently the map $\Phi$ is surjective, and $(h_1 , h_2 , h_3 ) $ maps to the same image with $(h^{'} _1 , h^{'} _2 , h^{'} _3 )$ if and only if there is a $g \in G$ such that $  (h^{'} _1 , h^{'} _2 , h^{'} _3 ) = (gh_1 ,gh_2 gh_3 )$. So the inverse image of any element in $C$ under $\Phi $ consists of $|G|$ elements. So we have
 $|B|= |G||C|$.

There is another natural map $\Psi : B\rightarrow A$, $(h_1 ,h_2 ,h_3 ) \mapsto (h_1 g_i h_1 ^{-1} , h_2 g_j h_2 ^{-1} )$.  It is easy to see $\Psi$ is surjective , and the inverse image of any element in $A$ under $\Psi$ consists of $v_i v_j v_k $ elements. So we have: $|B| = v_i v_j v_k |A| $.

So, $|A| = \frac{|G|}{ v_i v_j v_k  } |C|$. So we have $\gamma _{i,j} ^{k^{\star}} = \frac{1}{ v_i v_j } |C| = \frac{1}{v_i v_j } \bar{B} _{i,j, k^{\star }}$. Which implies    $c_i c_j = \sum _{k}  \frac{1}{ v_i v_j } \bar{B}_{i,j,k^{\star }} c_k $.  Compare this identity with (4) of theorem 2.1,  we see

$ [g_i ] \mapsto v_i c_i  $,$0\leq i\leq N $ extends to an isomorphism  from $\mathscr{C}_G $ to $C_G $. For later use, we denote this isomorphism as $J$.

 \end{pf}

 $\\$ Above isomorphism $J$ implies the following corollary.

 \begin{cor}
 There is a one to one correspondence between the set of irreducible representations of $G$ and the maximal ideals of the algebra $\mathscr{C}_G$.

 \end{cor}

  For any conjugacy class $[g_i ]$ of $G$ ,now  we define a function  $\Lambda _i$ on the set of irreducible representations of $G$ as follows.   For any irreducible representation $\rho _j $ , the value
  $\Lambda _i (  \rho _j  )$ is determined as follows.  Since $J( [g_i ]  ) \in C_G $ , then it acts on the representation $\rho _j $ as  multiplication of some constant $\mu  _i ^j $ by Schur lemma.  We set $\Lambda _i ( \rho _j )$ as this $\mu _i ^j $.

  Denote the set of complex functions on the set of irreducible representations of $G$ as $C(G, \mathbb{C} )$, which is a commutative algebra by multiplication of functions.

  \begin{lem}

  The map $[g_i ] \mapsto \Lambda _i $ $(i=0,1,...,N)$ extends to an isomorphism $ \mathscr{C}_G \rightarrow C(G, \mathbb{C}) $.

  \end{lem}

  \begin{pf}
  First there is a natural isomorphism (between algebras ) $K:  C_G \rightarrow C(G, \mathbb{C}) $,such that for any $c\in C_G$,the value of $K(c) $ on $ \rho _j $ is the constant $\lambda $ such that $K(c)$ acts on $\rho _j $ as $\lambda $ . Now We have a isomorphism $ K\circ  J : \mathscr{C} _G \rightarrow C(G, \mathbb{C})$, and by definition $\Lambda _i = K\circ J ( [g_i ] )$ so the lemma follows.
  \end{pf}
  $\\$

  We denote the character of $\rho _i $ as $ \chi _i  $, the value of $\chi _i $ on any element in $[ g_j  ]$ as $\lambda _i ^j$. Then $( \lambda _i ^j  ) _{(N+1) \times (N+1)}$ is the character matrix of $G$. Denote the dimension of $\rho _i $ as $n_i $.

  Since $J(  [g_j ] ) = v_j c_j $, we have $v_j w_j \lambda _i ^j = n_i \mu _i ^j $, which implies
  $ \Lambda _j (  \rho _i  ) = \mu  _i ^j = \frac {v_j w_j }{ n_i  } \lambda _i ^j = \frac{|G|}{n_i } \lambda _i ^j . $
  We obtain  the following identity satisfied by irreducible characters of $G$ as corollary of Lemma 3.1 .

  \begin{cor}
  Suppose $\rho _k $ is a irreducible representation of $G$ with dimension $n_k $ and character $\chi _k$,  $ [g_i ]$ and $[ g_j  ]$ be two conjugacy classes of $G$. Then

  $$  \chi _k ( g_i ) \chi _k ( g_j ) =\frac {n_k }{ |G| } \sum _{l} \bar{ \beta} _{i,j} ^l  \chi _k ( g_l  ) = \frac{ n_k }{ |G| } \sum _{l} \bar{ B } _{i,j, l^{\star }} \frac {1 }{v_l  } \chi _k (g_l ) $$

  \end{cor}

  \begin{pf}

  Compute the value $K\circ J ( [g_i ] [g_j ]  ) ( \rho _k   )$. Since $K\circ J$ is a morphism then it equals
  $$ K\circ J ([g_i ]) (\rho _k ) K\circ J ([g_j ]) (\rho _k ) = \frac{|G| ^2}{ n_k ^2 } \chi _k (g_i ) \chi _k ( g_j ) $$.
   On the other hand, since $ [ g_i ] [ g_j ] =\sum _{l} \bar{\beta }_{i,j} ^l  [g_l ] $, then it also equals
  $ \sum _{l}  \bar{\beta }_{i,j} ^l  \frac{|G|}{ n_k } \chi _k (g_l ) $. And the corollary follows.

  \end{pf}

  \section{Braid group actions}

 It could be interesting to pursuit deeper similarities between conjugacy classes and irreducible representations. But from now on we enter three dimensional topology starting from the following key observation.

 Recall in proof of $(2)$ of Theorem \ref{thm:centermodular}, we introduced certain bijective maps $\phi _{v,v+1}$. Roughly speaking, the following lemma says that these maps $\phi _{v,v+1}$ satisfy nontrivial braid relations. We put it in more details as follows. Choose $n$ elements $g_1 , g_2 ,..., g_n $ in $G$, let

 $$  \mathscr{D} _{ \{
 g_1  , g_2  ,...,g_n  \} } = \sqcup  _{s\in  S_n } \mathscr{D} _{ g_{s(1)} , g_{s(2)},..., g_{s(n)} } ,$$
 where $S_n $ denote the $n$-th permutation group.

 \begin{defi}
 For $1\leq v\leq n-1$, we define two maps $\phi _{v,v+1} $ and $\psi _{v,v+1}  $ from $\mathscr{D} _{\{g_1  , g_2  ,...,g_n  \} }$ to itself, which send elements of $\mathscr{D} _{ g_{s(1)} , g_{s(2)} ,..., g_{s(n)} } $to $\mathscr{D} _{ g_{s(1)}  ,..., g_{s(v+1)} ,g_{s(v)},..., g_{s(n)} }$ as follows.
 \begin{align*}
  \phi _{v,v+1}(v>1) :\mathscr{D}_{ g_{s(1) } , ..., g_{s(n) }} & \rightarrow \mathscr{D}_{ g_{s(1) } , ..., g_{s(v+1)} , g_{s(v)} ,..., g_{s(n) }}  \\
(a_1 , a_2 , ..., a_n ) & \mapsto ( a_1 , ..., a_{v-2} , a_{v-1} g_{s(v)}  a_v , a_{v}^{-1} g_{s(v)}^{-1} , a_v  a_{v+1} , a_{v+2}, ..., a_n    ) .\\
  \psi _{v,v+1} (v>1):\mathscr{D}_{ g_{s(1) } , ..., g_{s(n) }} & \rightarrow \mathscr{D}_{ g_{s(1) } , ..., g_{s(v+1)} , g_{s(v)} ,..., g_{s(n) }}  \\
(a_1 , a_2 , ..., a_n ) & \mapsto ( a_1 , ..., a_{v-2} , a_{v-1}  a_v , g_{s(v+1)}^{-1} a_{v}^{-1} , a_v g_{s(v+1)} a_{v+1} , a_{v+2}, ..., a_n    ) .\\
\phi _{1,2} :\mathscr{D}_{ g_{s(1) } , ..., g_{s(n) }} & \rightarrow \mathscr{D}_{ g_{s(2) } ,g_{s(1)} ,..., g_{s(n) }}  \\
(a_1 , a_2 , ..., a_n ) & \mapsto ( a_1 ^{-1} g_{s(1)} ^{-1} , a_1 a_2 , a_3 , ..., a_{n-1}, a_n g_{s(1)} a_1    ) .\\
\psi _{1£¬2} :\mathscr{D}_{ g_{s(1) } , ..., g_{s(n) }} & \rightarrow \mathscr{D}_{ g_{s(2) } ,g_{s(1)}, ..., g_{s(n) }}  \\
(a_1 , a_2 , ..., a_n ) & \mapsto (g_{s(2)} ^{-1} a_1 ^{-1} , a_1 g_{s(2)} a_2 , a_3 ,  ..., a_{n-1}, a_n a_1   ) .\\
\end{align*}
\end{defi}

\begin{thm}
 \label{thm:satisfybraid}
 (1)  $\psi _{v,v+1} \phi _{v,v+1} =id $ and $\phi _{v,v+1} \psi _{v,v+1} =id $, so $\phi _{v,v+1} $ is a bijection. From now on we denote the map $\psi _{v,v+1}$ as $\phi _{v,v+1}^{-1}$.

(2)  $\phi _{u, u+1} \phi _{v,v+1} = \phi _{v,v+1} \phi _{u,u+1}$ if $|u-v| \geq 2$.

(3)  $ \phi _{v,v+1} \phi _{v+1 , v+2 } \phi _{v, v+1} = \phi _{v+1, v+2 } \phi _{v, v+1} \phi _{v+1, v+2} $.

 \end{thm}

 \begin{pf} These identities are proved by direct computation. we only need to take care of the change of lower index.
 (1) Action of $\psi_{v,v+1} \phi_{v,v+1}$ on an element $(a_1 , ..., a_n ) \in \mathscr{D}_{g_{s(1)} , ..., g_{s(n)}}$ is described as follows.
\begin{align*}
 &(a_1 , ..., a_n ) \mapsto  ( a_1 , ..., a_{v-2} , a_{v-1} g_{s(v)}  a_v , a_{v}^{-1} g_{s(v)}^{-1} , a_v  a_{v+1} , a_{v+2}, ..., a_n     ) \\
&\mapsto ( a_1, ..., (a_{v-1} g_{s(v)} a_v) (a_{v}^{-1} g_{s(v)}^{-1}), g_{s(v)}^{-1} (a_{v}^{-1} g_{s(v)}^{-1})^{-1} , (a_{v}^{-1} g_{s(v)}^{-1})g_{s(v)} (a_v a_{v+1}) , a_{v+2}, ...,a_{n}   )\\
&=(a_1 ,..., a_n ).
 \end{align*}
 So we have $\psi _{v,v+1} \phi _{v,v+1} =id $. The other identity  $\phi _{v,v+1} \psi _{v,v+1} =id $ can be proved similarly.

 (2) We only need to prove the cases $|u-v|=2$ .  The cases $|u-v|>2$ are evident. So, suppose $u= v+2$ and $v \geq 2 $,  action of the left side on an element $(a_1 ,..., a_n ) \in \mathscr{D} _{g_{s(1)} , ..., g_{s(n)} }$ is described as follows.
\begin{align*}
&(a_1 ,..., a_n ) \mapsto  ( a_1 , ..., a_{v-2} , a_{v-1} g_{s(v)}  a_v , a_{v}^{-1} g_{s(v)}^{-1} , a_v  a_{v+1} , a_{v+2}, ..., a_n     ) \\
&\mapsto (a_1 , ...,  a_{v-2} , a_{v-1} g_{s(v)}  a_v , a_{v}^{-1} g_{s(v)}^{-1} , (a_v a_{v+1}) g_{s(v+2)} a_{v+2} ,
 a_{v+2} ^{-1} g_{s(v+2)} ^{-1} , a_{v+2} a_{v+3} , a_{v+4} ,..., a_{n} ).
 \end{align*}
 Action of the right side of (2) on the same element is described as follows.
\begin{align*}
 &(a_1 , ..., a_n ) \mapsto ( a_1, ..., a_v , a_{v+1} g_{s(v+2)} a_{v+2} , a_{v+2} ^{-1} g_{s(v+2)} ^{-1} , a_{v+2} a_{v+3} , a_{v+4} , ...,a_{n })\\
&\mapsto (a_1 ,..., a_{v-2} , a_{v-1} g_{s(v)} a_{v} , a_{v}^{-1} g_{s(v)}^{-1} , a_{v} ( a_{v+1} g_{s(v+2)} a_{v+2} ) ,  a_{v+2} ^{-1} g_{s(v+2)} ^{-1} , a_{v+2} a_{v+3} , a_{v+4} , ...,a_{n }).
 \end{align*}
 So the actions are the same. The case $v=1, u=v+2$ can be proved similarly.

 (3) First we consider the cases $v\geq 2$. Action of the left side of (3) on an element $(a_1 , ..., a_n ) \in \mathscr{D}_{g_{s(1)} , ..., g_{s(n)}}$ is described as follows.
 \begin{align*}
 &(a_1 , ..., a_n ) \mapsto  ( a_1 , ..., a_{v-2} , a_{v-1} g_{s(v)}  a_v , a_{v}^{-1} g_{s(v)}^{-1} , a_v  a_{v+1} , a_{v+2}, ..., a_n     ) \\
 &\mapsto (a_1 , ..., a_{v-2} , a_{v-1} g_{s(v)} a_v , (a_v ^{-1} g_{s(v)} ^{-1} ) g_{s(v)} (a_{v} a_{v+1}) , (a_v a_{v+1})^{-1} g_{s(v)}^{-1} , (a_{v} a_{v+1} )a_{v+2} ,...,a_n  ) \\
&\mapsto (a_1 ,..., (a_{v-1} g_{s(v)} a_v ) g_{s(v+1)} (a_{v+1}) , a_{v+1}^{-1} g_{s(v+1)} ^{-1} ,
 a_{v+1} (a_v a_{v+1})^{-1} g_{s(v)} ^{-1} , a_v a_{v+1} a_{v+2} , ..., a_n )\\
 &=(  a_1 ,..., a_{v-1} g_{s(v)} a_v  g_{s(v+1)} a_{v+1} , a_{v+1}^{-1} g_{s(v+1)} ^{-1} ,
 a_v ^{-1} g_{s(v)} ^{-1} , a_v a_{v+1} a_{v+2} , ..., a_n ).
 \end{align*}
 Action of the right side on the same element is described as follows.
\begin{align*}
 &(a_1 , ..., a_n )\mapsto ( a_1 ,..., a_{v-1} , a_v g_{s(v+1)} a_{v+1} , a_{v+1} ^{-1} g_{s(v+1)} ^{-1} , a_{v+1} a_{v+2} ,..., a_n )  \\
&\mapsto (a_1 , ..., a_{v-1} g_{s(v)} (a_v g_{s(v+1)} a_{v+1}) , (a_v g_{s(v+1)} a_{v+1}) ^{-1} g_{s(v)}^{-1}  \\
 &\ \ \ \ \ ,(a_v g_{s(v+1)} a_{v+1} ) a_{v+1}^{-1} g_{s(v+1)}^{-1} , a_{v+1} a_{v+2} ,..., a_n ) \\
&= ( a_1 , ..., a_{v-1} g_{s(v)} a_v g_{s(v+1)} a_{v+1} , (a_v g_{s(v+1)} a_{v+1}) ^{-1} g_{s(v)}^{-1} , a_v  , a_{v+1} a_{v+2} ,..., a_n   )\\
 &\mapsto (a_1 , ..., a_{v-1} g_{s(v)} a_v g_{s(v+1)} a_{v+1} ,a_{v+1} ^{-1} g_{s(v+1)}^{-1} a_v ^{-1} g_{s(v)}^{-1} (g_{s(v)}) a_v , a_v ^{-1} g_{s(v)} ^{-1} , a_v a_{v+1} a_{v+2} ,..., a_n  ).
 \end{align*}
 Comparing the results shows the actions are the same. The case for $v=1$ can be proved similarly.

\end{pf}

\begin{rem}
(1) Action of $\phi _{v,v+1} \phi _{v , v+1 } $ on the element $(a_1 ,..., a_n ) \in \mathscr{D} _{g_{s(1)} , ..., g_{s(n)} }$ is described as follows.
\begin{align*}
&(a_1 ,..., a_n ) \mapsto  ( a_1 , ..., a_{v-2} , a_{v-1} g_{s(v)}  a_v , a_{v}^{-1} g_{s(v)}^{-1} , a_v  a_{v+1} , a_{v+2}, ..., a_n     ) \\
&\mapsto (a_1 ,..., a_{v-2} , (a_{v-1} g_{s(v)} a_v ) g_{s(v+1)} a_v ^{-1} g_{s(v)} ^{-1} , g_{s(v)} a_v g_{s(v+1)} ^{-1} , a_v ^{-1} g_{s(v)} ^{-1} a_v a_{v+1} ,...,a_n ) .
\end{align*}
 We see usually the result isn't $(a_1 , ..., a_n )$, which means the braid actions are nontrivial ,or equivalently
$\phi _{v,v+1} ^2 \neq id $.

(2) If we take $g_1 = g_2 =...=g_n =g $, then we have action of $B_n$ on the set $\mathscr{D}_{g,...,g}$.

(3) Above proof has nothing to do with the fact that $G$ is a finite group. If we take a Lie group $G$, we can still define sets $\mathscr{D}_{\{ g_1 ,..., g_n \} }$,and the maps $\phi_{v,v+1}$ and the lemma still holds. So we have many new actions of braid groups on manifolds by smooth homeomorphisms.

\end{rem}

Since many representations of braid groups are related with some link invariants, for example, the braid group representations from R-matrix and quantum groups, we want to see if there are some link invariants related with above braid group actions. The map $\phi_{v,v+1}$ looks like a R-matrix. But corresponding to the fact that a R-matrix acts on two components in a tensor product of quantum group representations, the map $\phi_{v,v+1}$ change three position in the element $(a_1 ,..., a_n )$. We found the following essentially equivalent but more clear braid group actions by "breaking $a_i$ into $b_i b_{i+1} ^{-1}$.".

\begin{defi}
\label{defi:refinedbraidaction}
Let $G$ be any group, choose elements $g_1 ,..., g_n $ in $G$, we set
$\mathscr{G}_{g_1, ..., g_n} = (G)_{g_1} \times (G)_{g_2 } \times ... \times (G)_{g_n } =G^n $,  and define maps $\psi_{i;g_1 ,..., g_n } $ and $\psi _{i; g_1 ,..., g_n} ^{-1}$ $(1\leq i\leq n-1)$ as follows.
\begin{align*}
\psi_{i;g_1,...,g_n }: \mathscr{G}_{g_1 ,..., g_i ,g_{i+1} ,..., g_n } &\rightarrow   \mathscr{G}_{g_1 ,..., g_{i+1},g_{i} ,..., g_n } \\
(b_1, ..., b_i ,b_{i+1} ,..., b_n ) &\rightarrow (b_1 ,...,b_{i+1} b_i ^{-1} g_i ^{-1} b_i ,b_i ,..., b_n  )\\
\psi_{i;g_1 ,..., g_n } ^{-1} : \mathscr{G}_{g_1, ..., g_i , g_{i+1}, ...,g_n } &\rightarrow \mathscr{G}_{g_1, ..., g_{i+1} ,g_{i} ,...,g_n } \\
(b_1 ,..., b_i ,b_{i+1} ,..., b_n )&\mapsto (b_1 ,...,b_{i+1} , b_i b_{i+1}^{-1} g_{i+1} b_{i+1} ,..., b_n )
\end{align*}

\end{defi}

\begin{rem} (1)For any $g_1 , ..., g_n $ ,the set $\mathscr{G}_{g_1 ,..., g_n }$ is the same set $ G\times G\times ...\times G$. The symbols $g_1 ,..., g_n $ are used to define the maps $\psi _i$. Notice for different series $\{ g_1 ,..., g_n \} $, the maps $\sigma_i $ are different. Another notation $(G)_{g_1} \times (G)_{g_2 } \times ... \times (G)_{g_n } $ for $\mathscr{G}_{g_1 ,..., g_n }$ is used to form a parallel with the tensor product $V_{\lambda _1} \otimes ...\otimes V_{\lambda _n }$ of quantum group representations, where $\lambda_i $ indicates irreducible representations, just dual to $g_i $ (or its conjugacy class $[g_i ]$) in our "far conjugacy class- irreducible representation duality".

(2)For simplicity from now on we denote a sequence like $(b_1 ,..., b_i ,b_{i+1} ,..., b_n )$ just as $(..., b_i , b_{i+1} ,...)$.
\end{rem}

\begin{thm}
\label{thm:basicbraidrelation}
\begin{enumerate}
\item[ (1)] The map $\psi_{i,g_1 ,..., g_n } $ defined above is a bijection with $\psi_{i ; g_1 ,..., g_{i+1} ,g_i ,..., g_n  } ^{-1}$ as inverse.

\item[(2)] $\psi _{i; g_1 ,..., g_{j+1} ,g_j ,..., g_n } \psi _{j; g_1 ,..., g_n } = \psi _{j; g_1 ,...,g_{i+1} ,g_i ,..., g_n }  \psi_{i; g_1 ,..., g_n } $ for $|i-j|\geq 2$.

\item[(3)] $\psi _{i; g_1 ,..., g_{i+1} , g_{i+2} ,g_i ,...,g_n } \psi _{i+1; g_1, ..., g_{i+1} ,g_i ,g_{i+2} ,..., g_n } \psi _{i; g_1 ,..., g_n  }$

$= \psi_{i+1; g_1, ..., g_{i+2} ,g_i ,g_{i+1} ,...,g_n } \psi _{i;g_1 ,..., g_i ,g_{i+2} ,g_{i+1} ,..., g_n  } \psi _{i+1; g_1 ,..., g_n }$.
\end{enumerate}
\end{thm}

\begin{pf} We only prove (3), other two statements can be proved similarly by direct computation. The action of the left side on an element $(b_1 ,..., b_n ) \in \mathscr{G}_{g_1 ,..., g_n }$ is as follows.
\begin{align*}
&(...,b_i ,b_{i+1},b_{i+2}, ... ) \mapsto (...,b_{i+1} b_i ^{-1} g_i ^{-1} b_i ,b_i  ,b_{i+2},... )  \\
&\mapsto (...,b_{i+1} b_i ^{-1} g_i ^{-1} b_i , b_{i+2} b_i ^{-1} g_i ^{-1} b_i ,b_i ,... ) \\
&\mapsto (...,b_{i+2} b_i ^{-1} g_i ^{-1} b_i (b_{i+1} b_i ^{-1} g_i ^{-1} b_i )^{-1} g_{i+1} ^{-1} (b_{i+1} b_i ^{-1} g_i ^{-1} b_i ), b_{i+1} b_i ^{-1} g_i ^{-1} b_i, b_i ,... )\\
&= (...,b_{i+2} b_{i+1} ^{-1} g_{i+1}^{-1} b_{i+1} b_i ^{-1} g_i ^{-1} b_i , b_{i+1} b_i ^{-1} g_i ^{-1} b_i, b_i ,... ).
\end{align*}

The action of the right side of (3) on the same element is as follows.
\begin{align*}
&(...,b_i , b_{i+1} , b_{i+2} ,...) \mapsto (..., b_i , b_{i+2 } b_{i+1} ^{-1} g_{i+1} ^{-1} b_{i+1} , b_{i+1} ,...)\\
&\mapsto ( ...,b_{i+2 } b_{i+1} ^{-1} g_{i+1} ^{-1} b_{i+1} b_i ^{-1} g_i ^{-1} b_i ,b_i , b_{i+1} ,... ) \\
&\mapsto (..., b_{i+2 } b_{i+1} ^{-1} g_{i+1} ^{-1} b_{i+1} b_i ^{-1} g_i ^{-1} b_i , b_{i+1} b_i ^{-1} g_i ^{-1} b_i , b_i ,... ).
\end{align*}
Comparison of the results give (3).

\end{pf}

Now we relate above two kinds of braid group actions.
\begin{defi}
For a sequence of elements $g_1 ,..., g_n $ of $G$, define a subset of $G^{n}$
$$\bar{\mathscr{D}  }_{g_1 ,..., g_n }= \{ (b_1 ,...,b_n )\in G^{n} |
b_1 ^{-1} g_1 b_1 b_2 ^{-1} g_2 b_2 ... b_n ^{-1} g_n b_n =e \} .$$

\end{defi}
We let the group $G$ acts on $\bar{\mathscr{D}}_{g_1,..., g_n }$ from the right side as follows. For any $(b_1 ,..., b_n ) \in \mathscr{D}_{g_1 ,..., g_n } $ and any $g\in G$, let $ (b_1 ,..., b_n ) \cdot g = (b_1 g , ..., b_n  g ) $.  Then we have

\begin{lem}
\begin{enumerate}
\item[(1)] The map $\psi _{i, g_1 ,..., g_n }$ restricts to an bijection from $\bar{\mathscr{D}  }_{g_1 ,..., g_n }$ to
$\bar{\mathscr{D}  }_{g_1 ,...,g_{i+1} ,g_i ,..., g_n }$.

\item[(2)]The set $\bar{\mathscr{D}  }_{g_1 ,..., g_n }$ is stable by the action of $G$. Denote the quotient set as
$ \bar{\mathscr{D}  }_{g_1 ,..., g_n } /G $, in which the equivalent class containing $(b_1 ,..., b_n )$ will be denoted as $[ (b_1 ,..., b_n ) ]$.

\item[(3)] The maps $\psi _{i; g_1 ,..., g_n }$ commutes with the action of $G$.
\end{enumerate}
\end{lem}
\begin{pf}
(1) An element $(b_1 ,..., b_n ) \in \bar{\mathscr{D}  }_{g_1 ,..., g_n }$ is mapped to $(..., b_{i+1} b_i ^{-1} g_i ^{-1} b_i , b_i ,... )$ by  $\psi _{i; g_1 ,..., g_n }$. The following identities show this element is in$\bar{\mathscr{D}  }_{g_1 ,...,g_{i+1} ,g_i ,..., g_n }$.
\begin{align*}
 &b_1 ^{-1} g_1 b_1 ... (  b_{i+1} b_i ^{-1} g_i ^{-1} b_i )^{-1} g_{i+1}  (  b_{i+1} b_i ^{-1} g_i ^{-1} b_i ) b_i ^{-1} g_i  b_i ... b_n ^{-1} g_n ^{-1} b_n   \\
&= b_1 ^{-1} g_1 b_1 ... b_i ^{-1} g_i b_i b_{i+1}^{-1} g_{i+1} b_{i+1} ... b_n ^{-1} g_n  b_n =e .
\end{align*}

 Then we can show  that the inverse map $ \psi _{i; g_1 ,..., g_{i+1} ,g_i ,..., g_n } ^{-1} $ of $\psi _{i; g_1 ,..., g_n }$ send every element of $\bar{\mathscr{D}  }_{g_1 ,...,g_{i+1} ,g_i ,..., g_n }$ to $\bar{\mathscr{D}  }_{g_1  ,..., g_n }$ similarly ,which implies the claimed bijectivity.

(2) For any $( b_1 ,..., b_n ) \in \bar{\mathscr{D}  }_{g_1 ,..., g_n } $ and any $g\in G$, we have
\begin{align*}
 &(b_1 g ) ^{-1} g_1 (b_1 g ) (b_2 g)^{-1} g_2 (b_2 g) ... (b_n g)^{-1} g_n (b_n g) \\
 &= g^{-1} b_1 ^{-1} g_1 b_1 b_2 ^{-1} g_2 b_2 ...b_n ^{-1} g_n b_n g= g^{-1} e g=e .
\end{align*}

so we have $(b_1 ,..., b_n   )\cdot g \in \bar{\mathscr{D}  }_{g_1 ,..., g_n } . $

(3) We only need to show for any $ ( b_1 ,..., b_n ) \in \bar{ \mathscr{D}}_{g_1 ,..., g_n } $ ,any $g\in G$ and any $1\leq i\leq n-1$,  $\psi _{i, g_1 ,..., g_n  } ( (b_1 ,..., b_n    )\cdot g   )= (\psi _{i, g_1 ,..., g_n  } ( (b_1 ,..., b_n    )  ))\cdot g  $.  Now the left side equals
\begin{align*}
&(  b_1 g,..., (b_{i+1} g ) (b_i g)^{-1} g_i ^{-1} (b_i g) , b_i g,..., b_n g  )\\
&=(b_1 g,..., b_{i+1} b_i ^{-1} g_i ^{-1} b_i g , b_i g,..., b_n g  )= (b_1 ,..., b_{i+1} b_i ^{-1} g_i ^{-1} b_i , b_i ,..., b_n   ) \cdot g
\end{align*}
equals the right side.
\end{pf}\\

As a corollary of above lemma, the map  $\psi _{i; g_1 ,..., g_n }$ induce a bijective map from $\bar{\mathscr{D}  }_{g_1 ,..., g_n }/G$ to $\bar{\mathscr{D}  }_{g_1 ,...,g_{i+1} ,g_i ,..., g_n }/G $.  We denote it as $\bar{\psi } _{i; g_1 ,..., g_n }$. It is evident that these maps $\bar{\psi } _{i; g_1 ,..., g_n }$ satisfy the same relations in Theorem \ref{thm:basicbraidrelation}.

\begin{thm}
\label{thm:identifybraidaction}
(1) The following map $H_{g_1 ,..., g_n } $ is well defined ,and is bijective.
\begin{align*}
 \bar{\mathscr{D}  }_{g_1 ,..., g_n }/G &\rightarrow  \mathscr{D}_{g_1 ,..., g_n }\\
 [(b_1 ,..., b_n )] &\mapsto (b_1 b_2 ^{-1}  ,..., b_{n-1} b_n ^{-1} , b_n b_1 ^{-1}  ).
\end{align*}
(2) Maps $H_{g_1 ,..., g_n }$ commute with braid group actions in the obvious sense. So through these bijections
$H_{g_1 ,..., g_n }$, the braid group actions on sets $ \mathscr{D}_{g_1 ,..., g_n } $ and $ \bar{\mathscr{D}  }_{g_1 ,..., g_n }/G $ are identified.
\end{thm}
\begin{pf}
(1) The map is well defined because any other representative $(b_1 g ,..., b_n g  )$ is mapped to $( b_1 g (b_2 g)^{-1} ,..., b_n g (b_1 g)^{-1} ) $ ,which is the same element. The map is surjective because for any $(a_1 ,..., a_n )\in
\mathscr{D}_{g_1 ,..., g_n }$, we have
 \begin{align*}
 &( e, a_1 ^{-1} , (a_1 a_2)^{-1} , (a_1 a_2 a_3)^{-1} ,..., (a_1 ... a_{n-1} )^{-1} ) \in \bar{\mathscr{D}  }_{g_1 ,..., g_n }, and \\
&H_{g_1 ,..., g_n } ([ ( e, a_1 ^{-1} , (a_1 a_2)^{-1} , (a_1 a_2 a_3)^{-1} ,..., (a_1 ... a_{n-1} )^{-1} ) ]  )=
 (a_1 ,..., a_n ) .
 \end{align*}

 The map is injective because, if $[(b_1 ,..., b_n ) ]$ and $[(d_1 ,..., d_n ) ]$ are mapped to the same element, which means $b_1 b_2 ^{-1} = d_1 d_2 ^{-1}$, ... ,$b_n b_1 ^{-1} =d_n d_1 ^{-1}  $, which implies $d_1 ^{-1} b_1 = d_2 ^{-1} b_2 =...= d_n ^{-1} b_n  $. So we have $[(b_1 ,..., b_n ) ] =[(d_1 ,..., d_n )]$.

 (2) We only need to show $\phi _{i; g_1 ,..., g_n } H_{g_1 ,..., g_n } = H_{g_1 ,..., g_{i+1} ,g_i ,..., g_n } \bar{\psi }_{i; g_1 ,..., g_n } $. In the cases $i\geq 2 $, for any $(b_1 ,..., b_n ) \in \bar{\mathscr{D}  }_{g_1 ,..., g_n } $, the action of the left side on $[(b_1 ,..., b_n )]$ is as follows.
\begin{align*}
 &[(b_1 ,..., b_n )] \mapsto (b_1 b_2 ^{-1} ,..., b_n b_1 ^{-1} )  \\
&\mapsto (b_1 b_2 ^{-1} , ...,( b_{i-1} b_i ^{-1} ) g_i (b_i b_{i+1} ^{-1} ), (b_i b_{i+1}^{-1} )^{-1} g_i ^{-1} , (b_i b_{i+1}^ {-1}) b_{i+1} b_{i+2} ^{-1} ,..., b_n b_1 ^{-1}   )
\end{align*}

 Action of the right side map  on the same element is as follows.
\begin{align*}
 &[(b_1 ,..., b_n )] \mapsto [(b_1 ,..., b_{i+1} b_i ^{-1} g_i ^{-1} b_i , b_i , ..., b_n  ) ] \\
&\mapsto ( b_1 b_2 ^{-1} ,..., b_{i-1} ( b_{i+1} b_i ^{-1} g_i ^{-1} b_i  )^{-1} , b_{i+1} b_i ^{-1} g_i ^{-1} b_i (b_i )^{-1} ,b_i b_{i+2}^{-1} ,..., b_n b_1 ^{-1}  ).
\end{align*}
A comparison of the resulted elements finish the proof.
 \end{pf}\\

 Next we show the braid group actions on the sets $\bar{\mathscr{D}}_{g_1 ,..., g_n } /G  $ essentially only depend on the conjugacy classes $[g_1] ,..., [g_n ]$.   Suppose we replace $g_i $ in the sequence $g_1 ,..., g_n $ with a conjugate $h g_i h^{-1}$. Then we have

 \begin{prop}
 \label{prop:identifytwobraidactions}
 The following map
 \begin{align*}
 H_{i;h} : \bar{\mathscr{D}}_{g_1 ,..., g_n } &\rightarrow \bar{\mathscr{D}}_{g_1 ,..., h g_i h^{-1} ,..., g_n }\\
 (b_1 ,..., b_n ) &\mapsto (b_1 ,..., h b_i ,..., b_n  ).
 \end{align*}
 is a well defined bijective map ,and commutes with braid group actions in the obvious sense.
\end{prop}
  \begin{pf}
 Welldefinedness of the map is easy to see. It is bijective because we can easily show the map $[(b_1 ,..., b_n ) \mapsto (b_1 ,..., h^{-1}b_i ,..., b_n  )   ]$ is its inverse.

 For the second statement, we only need to show that $H_{i+1; h} \psi _{i, g_1 ,..., g_n } =
 \psi _{i, g_1 ,..., hg_i h^{-1} ,..., g_n } H_{i;h } $. Here $H_{i+1 ;h}$ is a map from $\bar{\mathscr{D}}_{g_1 ,...,g_{i+1} , g_i,..., g_n }$ to $ \bar{\mathscr{D}}_{g_1 ,...,g_{i+1} , hg_i h^{-1},..., g_n } $.

 Now, action of the left side map on an element $(b_1 ,..., b_n  )$ is as follows.
$$(b_1 ,..., b_n) \mapsto (..., b_{i+1} b_i ^{-1} g_i ^{-1} b_i ,b_i ,... ) \mapsto (..., b_{i+1} b_i ^{-1} g_i ^{-1} b_i , h b_i ,... ) .$$
 Action of the right side map on the same element is as follows.
$$(b_1 ,..., b_n ) \mapsto (...,h b_i , b_{i+1} ,...  ) \mapsto
  ( ..., b_{i+1} (h b_i )^{-1} (h g_i h^{-1} )^{-1} (h b_i) , h b_i ,...  ) .$$
A comparison of the results finish the proof.
\end{pf}\\

In the same way we can prove the braid actions on the set $\mathscr{D}_{\{ g_1 ,..., g_n  \}}$ also essentially depends on the conjugacy classes $[g_1] ,..., [g_n]$. By using Proposition \ref{prop:identifytwobraidactions} repeatedly, we see if replace every element in the sequence $g_1 ,..., g_n $ with a conjugate, forming a new sequence $h_1 g_1 h_1 ^{-1} ,..., h_n g_n h_n ^{-1} $ ,then the braid group action on sets $\{ \bar{\mathscr{D}}_{g_{s(1)} ,..., g_{s(n)} }  \} _{s\in S_{n}} $ is isomorphic to the braid group action on $\{ \bar{\mathscr{D}}_{h_{s(1)}  g_{s(1)} h_{s(1)} ^{-1} ,..., h_{s(n)} g_{s(n)} h_{s(n)} ^{-1} }  \} _{s\in S_{n}} $. So ,to emphasize the dependence on conjugacy classes, in some occasions we write the set $\bar{\mathscr{D}}_{g_1 ,..., g_n }$ as $\bar{\mathscr{D}}_{[g_1 ] ,..., [g_n ]} $ and talk about the braid group actions on $\bar{\mathscr{D}}_{[g_1 ] ,..., [g_n ]} $ .

Summing up results of last sections, for the aim to understand those new braid group actions
, we present the following "far duality picture ". Roughly speaking, one side of the duality is around conjugacy classes of finite groups,  yet another side is around irreducible representations of a compact Lie group (or the associated Lie algebra, or ,the associated quantum group ). Both sides show some mysteries parallel properties, for example ,the braid group actions , and relations to links and 3- manifolds  to be introduced in following sections.

\begin{tabular}{|l|l|l|}
\multicolumn{2}{c}{ Far duality picture}\\
\hline
   Conjugacy\ classes\ of\  $G$ & irreducible\ representations\ of $\mathcal {G} $ \\
\hline
 the\ center\ $C_{G}$ & the representation\ ring\ $R_{\mathcal{G}}$ \\
\hline
 a\ trivial\ class\ $[e]$ & a\ trivial\ representation\ $V_{\lambda_0 }$\\
\hline
 involution\ $[g] \mapsto [g^{-1}]$ & involution\ $V_{\lambda }\mapsto V_{\lambda ^{*}} $\\
\hline
  properties\ in\ Theorem \ref{thm:centermodular}   & properties\ in\ Theorem\ \ref{thm:repmodular} \\
\hline
the\ sets\ $\mathscr{D}_{[g_1 ],..., [g_n ]}$ and $\bar{\mathscr{D}}_{[g_1] ,..., [g_n ]}$ &
the spaces $Hom_{\mathcal {G}} (V_{\lambda_{0}}, V_{\lambda_{1}} \otimes ...\otimes V_{\lambda_{n}}  )$ \\
\hline
braid\ group\ actions\ on $\mathscr{D}_{[g_1 ],..., [g_n ]}$ &  braid\ group\ actions\ on\\
 by\ $\phi_i $.    &  $Hom_{\mathcal {G}} (V_{\lambda_{0} } , V_{\lambda_{1}} \otimes ...\otimes V_{\lambda_{n}})$\\

   & by\ quantum\ group\ and\ R-matrix.\\
\hline
?   & quantum\ link\ invariants. \\
\hline
?   & WRT\ invariants\ for\ 3-manifolds.\\
\hline
?   & categorification\ of\ quantum\ invariants.\\
\hline
\end{tabular}\\
\\

The correspondence between $\mathscr{D}_{[g_1] ,..., [g_n ]} ( \bar{\mathscr{D}}_{[g_1] ,..., [g_n ]}  ) $ and the space $Hom_{G} (V_{\lambda_0 } , V_{\lambda_{1}} \otimes ...\otimes V_{\lambda_n }  )$ can be strengthened by the fact
that the role played by the  numbers $A_{i_1 ,..., i_k } = \dim _{G} (\rho_0 , \rho_1 \otimes ...\otimes \rho_k  )$ in theorem1.1 is very similar with the role played by the numbers $\bar{B}_{i_1 ,..., i_k }= |\mathscr{D}_{[g_{i_1}] ,..., [g_{i_k }]} |$ in theorem 2.1. Yet in this correspondence, the compact Lie group $\mathcal{G}$ is substituted with a finite group $G$. That is why we call above duality a "conceptual" but not exact duality. Even though, it can show us the prospects before the new braid group actions.

If $g_1 = g_2 =...=g_n=g $,then we have action of the braid group $B_n $ on the set $\underbrace{ (G)_{g} \times ...\times (G)_{g}}_{n}=G^{n}
$. Now the action of a braid generator $\sigma_i $ only affects the $i-$th and $(i+1)-$th factors in $G^n$, just as the $R-$matrix. And when $g=e$, then the action degenerate to a action of the permutation group, just as when the quantum parameter $q$ goes to 1, the $R-$matrix goes to simply a transportation.
From these observations, we could recognize the map $\psi _{i;g,g}$ on $G\times G$ as a kind of "dual $R-$matrix" growing out of finite groups. Now, since the $R-$matrix can give quantum invariants of links and 3-manifolds, we could make the following hypothesises:
\begin{itemize}
\item[(1)] The new braid group actions can also give invariants for links and 3-manifolds.
\item[(2)]The theoretical frame for those  invariants for links and 3-manifolds could be similar with the theoretical frame for quantum invariants. So, the theory for quantum invariants can serve as a guide for us to look for the new invariants for links and 3-manifolds.
\end{itemize}
In the following sections we will certify these hypothesis by constructing those new invariants.

\begin{rem}
Since these braid group actions connect finite group with three dimensional topology, you may ask the relationship between them with the famous Dijkgraaf-Witten invariant \cite{DW}. They are indeed connected, roughly as follows. According to the settings of  Section 3.2 of \cite{BK} by B.Bakalov, J.A.Kirillov , consider the quantum double $D(G)$ of $\mathbb{C}G$, the group algebra of a finite group. For every conjugacy class $[g]$ of $G$, there is a module $V_{[g]}$ of $D(G)$, whose representation space is naturally identified with $\mathbb{C}G$. Then the map resulted from the universal $R-$matrix of $D(G)$:
$$\mathbb{C}G\otimes \mathbb{C}G \cong V_{[g_1]} \otimes V_{[g_2 ]} \rightarrow  V_{[g_2]} \otimes V_{[g_1 ]} \cong \mathbb{C}G\otimes \mathbb{C}G $$
can be identified with the key map $\psi_{2, g_1 , g_2}$.  Even though, we think the map $\psi_{2, g_1 , g_2}$ has its own meaning, because of the related braid group action $\phi_{v,v+1}$ (Theorem \ref{thm:satisfybraid}), because of the maps $\psi_{2, g_1 , g_2}$ and $\phi_{v,v+1}$ can be defined for any group but not only finite groups, and because many later constructions in this paper and forthcoming ones built on these maps.
\end{rem}

\section{Link invariants from new braid group actions }

Let $G$ be a finite group and let $g_1 ,..., g_n $ be a sequence of elements in $G$. In Definition \ref{defi:refinedbraidaction} we introduced certain bijective map $\psi _{i; g_1 ,..., g_n }$ from $G^n$ to itself. In the following we
let $g_1 = g_2 =...=g_n =g \in G$, and denote the map $\psi _{i;g ,..., g }$ simply as $\psi_{i;g}$. Now we linearize everything, and $\psi_{i; g}$ induces an automorphism of $\mathbb{C} G \otimes \mathbb{C} G \otimes ...\otimes \mathbb{C}G $ which will be still denoted as $\psi_{i;g}$. From Theorem \ref{thm:basicbraidrelation} we see $\psi_{g} = \psi _{1, g} $  is a $R-$matrix on $\mathbb{C} G \otimes \mathbb{C} G $. In the last part of section 5 we proved that if substitute $g$ with a conjugate $hgh^{-1}$ then the resulted braid actions are essentially the same, so if we want to emphasize the dependence on the conjugacy class we can write the map $\psi_{ g}$ as $\psi_{[g]}$.

Now according to the remarks after Definition \ref{defi:enhanced}, if we can enhance the $R-$matrix $\psi_{g}$ to be an enhanced $R-$matrix, then we obtain a link invariant. But after some investigation, we found this $R-$matrix possibly can't be enhanced. Instead of that, if we relax the conditions in Definition \ref{defi:enhanced} to give a definition of the following "extended $R-$matrix "  then, firstly every extended $R-$matrix still give a link invariant, secondly  the $R-$matrix $\psi_{g}$ can be "extended" to an extended $R-$matrix.

\begin{defi}
\label{defi:centraltor}
Suppose $I\in End (V\otimes V)$ is a $R-$matrix. An  element $h \in End(V)$ is called central with respect to $I$ if

$(id_{V} \otimes h )\cdot I = I\circ (h\otimes id_{V} )  $ and $(h\otimes id_{V}  )\cdot I = I \cdot (id\otimes h) .  $

\end{defi}

\begin{defi}
\label{defi:extented}
Suppose $V$ is a finite dimensional linear space.  Let $I\in End (V\otimes V)$, $f , c \in End(V)$ such that $c$ is invertible.  The data $(I , f , c)$ is called an extended $R-$matrix if the following conditions are satisfied.
\begin{enumerate}
\item[(1)] $I$ is a $R-$matrix.

\item[(2)] $(f\otimes f)\cdot I = I \cdot (f\otimes f).$

\item[(3)] $c$ is central with respect to $I$.

\item[(4)] $tr_2 ( ( id_{V}  \otimes f)\cdot I )=   c $, $tr_2 ( (id_{V}\otimes f) \cdot I^{-1} ) = c^{-1}  .$

\item[(5)] $f\cdot c=c\cdot f$.
\end{enumerate}
\end{defi}

 \begin{figure}[htbp]

  \centering
  \includegraphics[height=4.5cm]{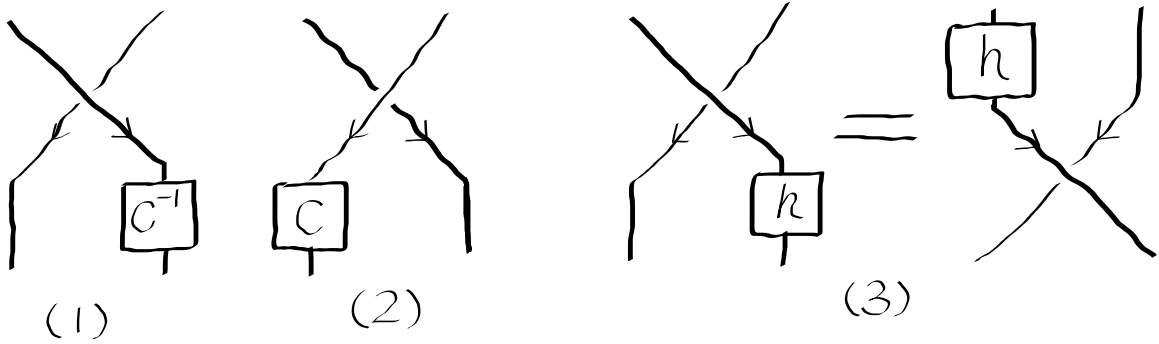}
  \caption{Modified R-matrix}
\end{figure}
\begin{rem}
\label{rem:extendedRmatrix}
(1)The extended $R-$matrix is a generalization of the enhanced $R-$matrix, since in Definition \ref{defi:enhanced} if we set $\mu=1$ as in the remark after Definition \ref{defi:enhanced}, then the data $(I, f, \lambda ( id_V )  )$ satisfies all conditions in Definition \ref{defi:extented}.

(2)If in an extended $R-$matrix $(I, f, c)$, the third component $c=id_V$, we call the triple $(I,f,c)$ as a special extended $R-$matrix.
\end{rem}
By using an extended $R-$matrix we can also construct a link invariant. The first step is to modify the $R-$matrix $I$.

\begin{lem}
\label{lem:modifiedRmatrix}
Let $(I, f, c)$ be an extended $R-$matrix on a linear space $V$. If we set $\bar{I}= (id_V \otimes c^{-1}  )\circ I $ ,then $\bar{I}$ is still a $R-$matrix and $\bar{I}^{-1} =  (c \otimes id_V )\circ I^{-1}$. Besides, $\bar{I} \circ (f\otimes f) = (f\otimes f) \circ \bar{I} . $ We call $\bar{I}$ as the $R-$matrix modified by $c$.
\end{lem}
\begin{pf}
The statements can be easily proved by using following graphs.
\end{pf}\\

 \begin{figure}[htbp]

  \centering
  \includegraphics[height=3.8cm]{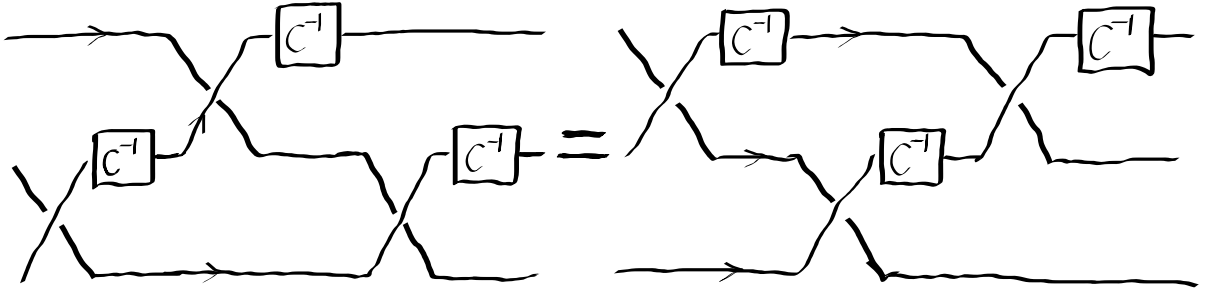}
  \caption{The Yang-Baxter relation }
\end{figure}

Since the element $\bar{I}$ is still a $R-$matrix, it induces a braid group representation $\bar{\rho} _n : B_n \rightarrow End(V^{\otimes n} ) $.

\begin{thm}
 \label{thm:extendedgivelink}
 Let $\beta \in B_n $, set $\bar {\rho} (\beta ) = tr ( f^{\otimes n} \cdot \bar{\rho}_n (\beta)  ) $. Then $\bar{\rho }(\beta) $ is invariant by both Markov moves, thus define a link invariant.
\end{thm}

\begin{pf}
For invariance under the first type of Markov move, let $\gamma \beta \gamma ^{-1}$ be any conjugate of $\beta$. Then
$$\bar{\rho }( \gamma \beta \gamma^{-1} ) = tr (f^{\otimes n} \cdot \bar{\rho}_n ( \gamma \beta \gamma ^{-1} )   )= tr ( \bar{\rho}_n (\gamma ) f^{\otimes n}  \bar{\rho}_n (\beta )  \bar{\rho} _n (\gamma^{-1})  )
= tr (f^{\otimes n}\bar{\rho}_n (\beta)  )= \bar{\rho} (\beta ) .$$
For the second type of Markov moves,
First, we have
\begin{align*}
 &tr_{n+1} ( f^{\otimes (n+1)} ( id_{V} ^{\otimes n} \otimes c^{-1}  ) (id_{V} ^{\otimes (n-1)} \otimes I  ) \bar{\rho}_{n} (\beta)  ) \\
& = tr_{n+1} ( f^{\otimes (n+1)} ( id_{V} ^{\otimes (n-1)} \otimes I ) ( id_{V} ^{ \otimes (n-1)} \otimes c^{-1} \otimes id_V  ) \bar{\rho}_n (\beta)   ) \\
&= ( f^{\otimes (n-1)} \otimes (f\cdot c)   ) ( id_{V} ^{\otimes (n-1)} \otimes c^{-1} ) \bar{\rho}_{n} (\beta ) = f^{\otimes n} \bar{\rho}_{n} (\beta).
\end{align*}
 So we have
\begin{align*}
  &\bar{\rho} ( \sigma_{n} \beta  ) = tr ( f^{\otimes (n+1)} ( id_{V} ^{\otimes n} \otimes c^{-1}  ) (id_{V} ^{\otimes (n-1)} \otimes I  ) \bar{\rho}_{n} (\beta)  ) \\
&= tr( tr_{n+1} (   f^{\otimes (n+1)} ( id_{V} ^{\otimes n} \otimes c^{-1}  ) (id_{V} ^{\otimes (n-1)} \otimes I  ) \bar{\rho}_{n} (\beta)   )  )
 =tr(   f^{\otimes n} \bar{\rho}_{n} (\beta) ) = \bar{ \rho } ( \beta ).
 \end{align*}

 Similarly we can prove  $ \bar{\rho} ( \sigma_{n}^{-1} \beta  ) =  \bar{\rho } (\beta ) $.

So the proof is finished.
\end{pf}\\

\begin{lem}
Let $(I, f, c)$ be an extended $R-$matrix on a linear space $V$. Suppose $\bar{I}$ is the modified $R-$matrix as in Lemma \ref{lem:modifiedRmatrix}. Then the triple $( \bar{I} , f, id_{V}  )$ is a special extended $R-$matrix, as defined in Remark \ref{rem:extendedRmatrix}.
\end{lem}
This lemma can be proved by direct computations. It is easy to see the link invariant resulted from $(\bar{I}  ,f, id_{V} )$ is the same with the link invariant resulted from $(I, f,c )$.

 \begin{figure}[htbp]

  \centering
  \includegraphics[height=4.5cm]{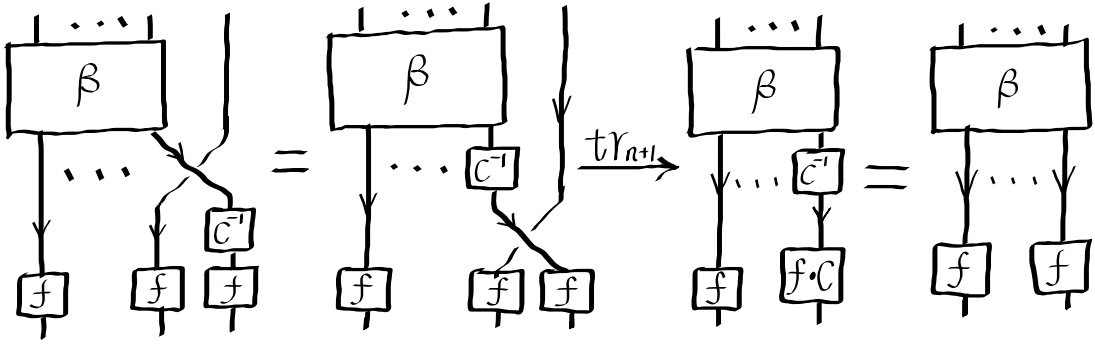}
  \caption{The second Markov move }
\end{figure}

\begin{defi}
\label{defi:extendedpair}
We call $( g,b )\in G\times G $ as a extended pair for $G$, if  $gb=bg$ and $b^2 =1$. We denote the set of extended pairs for $G$ as $\mathscr{E}G $.
\end{defi}

Easy examples are $\{ ( g, e ) \} _{g\in G}$ , where $e$ is the unit. In this way we can view $G$ as a subset of $\mathscr{E}G$ naturally.   We call two extended pair $(g_1 ,b_1)$ and $(g_2 , b_2 ) $ as conjugate to each other, if there is $h\in G$ such that $g_2 = h g_1 h^{-1}  $ and $b_2 = h b_1 h^{-1}  $. Conjugation is evidently a equivalence relation in the set of extended pairs, and we call the relevant equivalence classes as conjugacy classes (of extended pairs   ). Denote the conjugacy class containing  a extended pair $(g,b )$ as $[(g,b) ] $.

\begin{thm}
\label{thm:extendableGRmatrix}
Let $(g, b) $ be a extended pair for $G$ . Let $\psi _g $ be the $R-$matrix as before.  Denote $c=( bg^{-1} ) _{L} \in End (\mathbb{C}G )$ to be the left action of $b g^{-1}$ on $\mathbb{C} G$: $ x\mapsto b g^{-1} x $ for $x\in G$. Then $c^{-1} = ( b^{-1} g  )_{L} = (bg )_{L} $. Let $f= b_{L}\in End (\mathbb{C}G ) $. Then the data $( \psi_g ,  f  ,c )$ is an extended $R-$matrix.
\end{thm}

\begin{pf} We need to show the combination $( \psi_g ,  f  ,c )$ satisfies (1)$ \sim $(5) of Definition \ref{defi:extented}. (1) and (5) are evident. For (2),

let $x,y \in G$ are any elements, then
\begin{align*}
& (( f\otimes f  )\cdot \psi_{g}) ( x\otimes y ) = (f\otimes f) ( yx^{-1} g^{-1} x \otimes x  )= byx^{-1} g^{-1} x \otimes bx  , and\\
&(\psi_{g} \cdot (f\otimes f ) )(x\otimes y )= \psi_{g} ( bx\otimes by  )= by(bx)^{-1} g^{-1} (bx) \otimes bx =  by x^{-1} g^{-1} x \otimes bx  .
\end{align*}
Where the last equality sign follows from $bg=gb$.  For (3),
\begin{align*}
&(\psi_{g} \cdot (c\otimes id_{\mathbb{C} G} )   ) (x\otimes y) =\psi_g ( bg^{-1} x\otimes y  ) =y(bg^{-1} x  )^{-1} g^{-1} bg^{-1} x\otimes bg^{-1} x  \\
&= yx^{-1} g^{-1} x \otimes bg^{-1} x   ,\\
&(  id_{\mathbb{C} G} \otimes c ) \cdot \psi_{g} ( x\otimes y )= (id_{\mathbb{C} G} \otimes c ) ( yx^{-1} g^{-1} x \otimes x ) =yx^{-1} g^{-1} x \otimes bg^{-1} x  .
\end{align*}

For (4), first we have $ (f\otimes f)\cdot \psi_{g} (x\otimes y )= byx^{-1} g^{-1} x \otimes bx $.  It is easy to see that to get the partial trace $tr_2 $ of this map, we only need to "equate the second term then substitute it in the first term", concretely, equate the second term gives $ y=bx $, then
$tr_{2} ( (f\otimes f)\cdot \psi_{g}  ) x  = b(bx) x^{-1} g^{-1} x = b^2 g^{-1} x= (f\cdot c) x  $, which implies $tr_{2} ( (f\otimes f)\cdot \psi_{g}  )=f\cdot c .$ Similarly, from the identity $(f\otimes f) I^{-1} ( x\otimes y )= by \otimes bxy^{-1} gy  $,  equate the second term gives $ y= bxy^{-1} gy  $, which implies $y=gbx $. So $ tr_{2} ( (f\otimes f) I^{-1} ) (x) = b(gbx) = fc^{-1} (x)$, so we have $ tr_{2} ( (f\otimes f) I^{-1} )= fc^{-1}   $.

\end{pf}

\begin{rem}
\label{rem:extendableinvariant}
So by Theorem \ref{thm:extendedgivelink} and Theorem \ref{thm:extendableGRmatrix}, for any finite group $G$ and any extended pair $(g,b) $ for $G$  we have constructed a integer
invariant for oriented links. We denote this invariant as $\Lambda _{G;(g,b)} (L ) $, where $L$ is a oriented link. It is not hard to prove that if
  $(g_1 , b_1) $  conjugate to $ (g_2 ,b_2) $  then   $\Lambda _{G;(g_1 ,b_1 )} (L ) = \Lambda _{G;(g_2 ,b_2 )} (L ) $ for any $L$.  So it is suitable to  denote this invariant as $\Lambda _{G;[(g,b)]} (L ) $.
\end{rem}

\section{Invariants for oriented links colored by conjugacy classes }

Let $G$ be a finite group. In this section we define invariants of oriented links (and framed links ) colored by conjugacy classes of $G$ or conjugacy classes of extended pairs, as generalization of results of the last section, and as a basis for construction of  three manifold invariants. This construction is based on Theory of tangle category due to V.G.Turaev \cite{Tu2} and D.S.Freed, D.Yetter \cite{FY}. For a setting of colored tangles we consulted the paper by Y.Akutsu,T.Deguchi and T.Ohtsuki \cite{ADO}.   For such a approach, it is necessary to associate morphisms to the basic tangles
with shape "$\cup$ " and "$\cap $ ", luckily such morphisms exist naturally.

\begin{defi}
\label{defi:Gtangle} By a oriented, $G-$colored tame tangle, we mean a oriented tangle $T$ satisfying the following conditions.

(1) $T\subset \mathbb{C} \times [0,1]$;

(2) Every component of $T$ is associated with an element of $G$;

(3) The upper boundary  of $T$ is  $\{ 1,2,..., k \} \times \{ 1 \}  $ for some positive integer $k$ or empty. The lower boundary
is $\{ 1,2 ..., l \} \times \{ 0 \}  $ for some positive integer $l$ or empty.

If $T$ is a oriented $G-$ colored tame tangle, suppose the upper boundary and the lower boundary of $T$ are   $\{ 1,2,..., k \} \times \{ 1 \}  $ and
$\{ 1,2 ..., l \} \times \{ 0 \}  $ respectively, define
$\partial ^{+} T = [(g_1 , \epsilon _1  ), (g_2 , \epsilon _2 ),..., (g_k , \epsilon _k )   ] $ and
$\partial ^{-} T = [(h_1 , \mu _1 ), (h_2 , \mu _2 ),..., (h_l , \mu _l ) ] $  (  where $g_i , h_j \in G$ and $ \epsilon _i , \mu _j \in \{ +1 ,-1 \} $) according to the following rules:

(a) The arc of $T$ passing the boundary point $\{ i \} \times \{ 1 \} $ is colored by $g_i$ for $1\leq i\leq k $; the arc of $T$ passing the boundary
point $\{ j \} \times \{ 0 \}  $ is colored by $h_i $ for $1\leq j \leq l $;

(b) Taking into consideration of the orientation of arcs, if an arc is leaving (entering ) a upper boundary point $\{ i \} \times \{ 1 \} $, then $\epsilon _i = +1 (-1 )$;
 if an arc is entering ( leaving ) a lower boundary point $\{ j \} \times \{ 0 \} $, then $\mu _j = +1 (-1 )  $.

\end{defi}

\begin{rem}
\label{rem:colored tangle} (1) Eventually we will show the invariants to be constructed are essentially invariants for links colored by conjugacy classes
of $G$. Comparing with the theory of WRT  link invariants, where the invariants are defined for links colored by irreducible representations of quantum
groups. They two form a kind of "conceptional duality".

(2) $n-$string braid is a kind of tangle of special importance. We orient all its strings naturally from upward to downward. We call the string passing the
point $\{ i \} \times \{ 1 \}$ (as a end point) as the $i-$th string.
\end{rem}

Now we define a  "$G-$tangle category" $\mathcal {T}_{G}$ as follows.  The objects of $\mathcal{T}_{G}$ are sequences:

$[ (g_1 ,\epsilon_1   ), (g_2 , \epsilon_2 ),..., (g_k , \epsilon_k ) ] $, where $k=0,1,2,...$, $g_i \in G$, $\epsilon_i \in \{ +1,  -1  \} $. There is
a special object $\emptyset$ corresponding to the case $k=0$. For two objects  $X_1 =[(g_1 , \epsilon _1  ), (g_2 , \epsilon _2 ),..., (g_k , \epsilon _k )   ] $ and $X_2 = [(h_1 , \mu _1 ), (h_2 , \mu _2 ),..., (h_l , \mu _l ) ]$ , we set
$$ Hom (X_1 , X_2 )= \{  \ oriented\  G-tangle\  T\  |\  \partial ^{+} T = X_1 ,\  \partial ^{-} T = X_2 \  \} / isotopy .$$
Especially, the set
$Hom (\emptyset , \emptyset )$ is nothing but the set of isotopy classes of oriented $G-$colored links.

The composition $Hom (X_2 , X_3) \times Hom (X_1 , X_2  ) \rightarrow Hom (X_1 , X_3 ) $:$ (T_2  , T_1  ) \mapsto T_2 \circ T_1 $ is defined by the
operation of connecting two related tangles as shown in the following Figure 7.1.

 The category $\mathcal {T}_{G} $ has a natural tensor category strucure. For two objects

  $X_1 =[(g_1 , \epsilon _1  ), (g_2 , \epsilon _2 ),..., (g_k , \epsilon _k )   ] $ ,
 $X_2 = [(h_1 , \mu _1 ), (h_2 , \mu _2 ),..., (h_l , \mu _l ) ]$,the tensor product $X_1 \circ X_2$ is defined as
 $[g_1 , \epsilon _1  ), (g_2 , \epsilon _2 ),..., (g_k , \epsilon _k ),(h_1 , \mu _1 ), (h_2 , \mu _2 ),..., (h_l , \mu _l )  ] $.
 The tensor product of two morphisms (tangles ) $(T_1 , T_2  )\mapsto T_1 \otimes T_2 $ is explained in the following Figure 4.

 \begin{figure}[htbp]

  \centering
  \includegraphics[height=2.5cm]{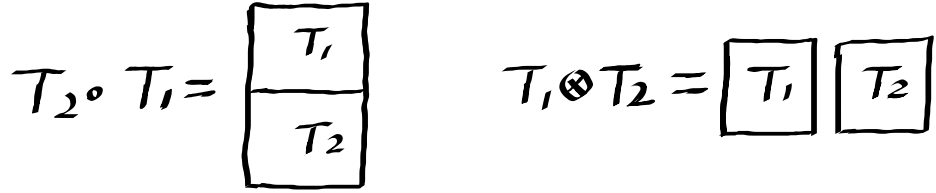}
  \caption{Composition of tangles }
\end{figure}

To present  the conception of tangle category we need to introduce another category of sliced oriented $G-$colored tangle diagrams
$\mathcal{D}\mathcal{T}_{G} $ as follows. First, a sliced oriented, $G-$colored tangle diagram means a  tangle diagram $T$ lying in $\mathbb{R} \times
[0,1] $ together with a set of horizontal lines such that

(1)each arc of $T$ is oriented and associated with an element of $G$;

(2)the upper boundary of $T$ is $ \{  1,2,..., k\} \times \{ 1 \} $, the lower boundary of $T$ is  $\{ 1,2,..., l \} \times \{ 0 \} $ for some natural numbers $k,l$;( whence we call such a tangle as "a $(k,l)$ tangle diagram  ")

(3)the part of the tangle diagram  $T$ between two adjacent horizontal lines (in above mentioned set of lines )consists of a disjoint union
of vertical lines and one of the elementary tangles shown in the Figure 5.

For a sliced oriented $G-$colored $(k,l)$ tangle diagram $T$ , we define:

$\partial ^{+} T =[(g_1 , \epsilon _1  ), (g_2 , \epsilon _2 ),..., (g_k , \epsilon _k ) ] $ and

$\partial ^{-} T =[(h_1 , \mu _1 ), (h_2 , \mu _2 ),..., (h_l , \mu _l ) ] $,

where the elements $g_1, ..., g_k , h_1 ,..., h_l \in G$, and the signatures $\epsilon_1 , ..., \epsilon_k , \mu_1 ,..., \mu_l \in \{ +1 , -1 \} $ are determined as the same way shown in $(3)$ of Definition \ref{defi:Gtangle}.

Now we define the category $\mathcal{D}\mathcal{T}_{G}$. The objects of $\mathcal{D}\mathcal{T}_{G}$ are the same as the objects of the category $\mathcal{T}_{G} $. For two objects

$X_1 =[(g_1 , \epsilon _1  ), (g_2 , \epsilon _2 ),..., (g_k , \epsilon _k )   ] $ and $X_2 = [(h_1 , \mu _1 ), (h_2 , \mu _2 ),..., (h_l , \mu _l ) ]$,

 a morphism from $X_1 $ to $X_2 $ is an equivalent class of a sliced, oriented , $G-$colored $(k,l)$ tangle $T$ such that
$\partial ^{+} T = X_1 $and $ \partial ^{-} T = X_2 $, with equivalence relations generated by relations represented by the Turaev moves $T1-T6$ as in Figure 6.
There is a natural functor $\mathcal{R}: \mathcal{D}\mathcal{T}_{G} \rightarrow \mathcal{T}_{G} $ by mapping a tangle diagram to the tangle represented by it.

\begin{thm}
\label{turaevtangle}
The functor  $\mathcal{R}$ is isomorphic. That is to say, $\mathcal{R}$ induces a bijection between the morphisms of the two categories.
\end{thm}

 \begin{figure}[htbp]
 \label{graph:basictangle}

  \centering
  \includegraphics[height=2.3cm]{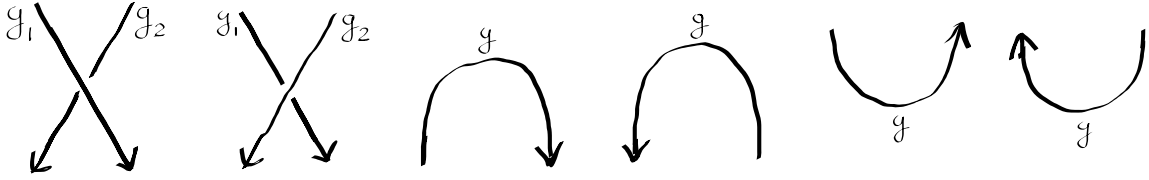}
  \caption{Basic tangles}
\end{figure}

\begin{figure}[htbp]

  \centering
    \includegraphics[height=13cm]{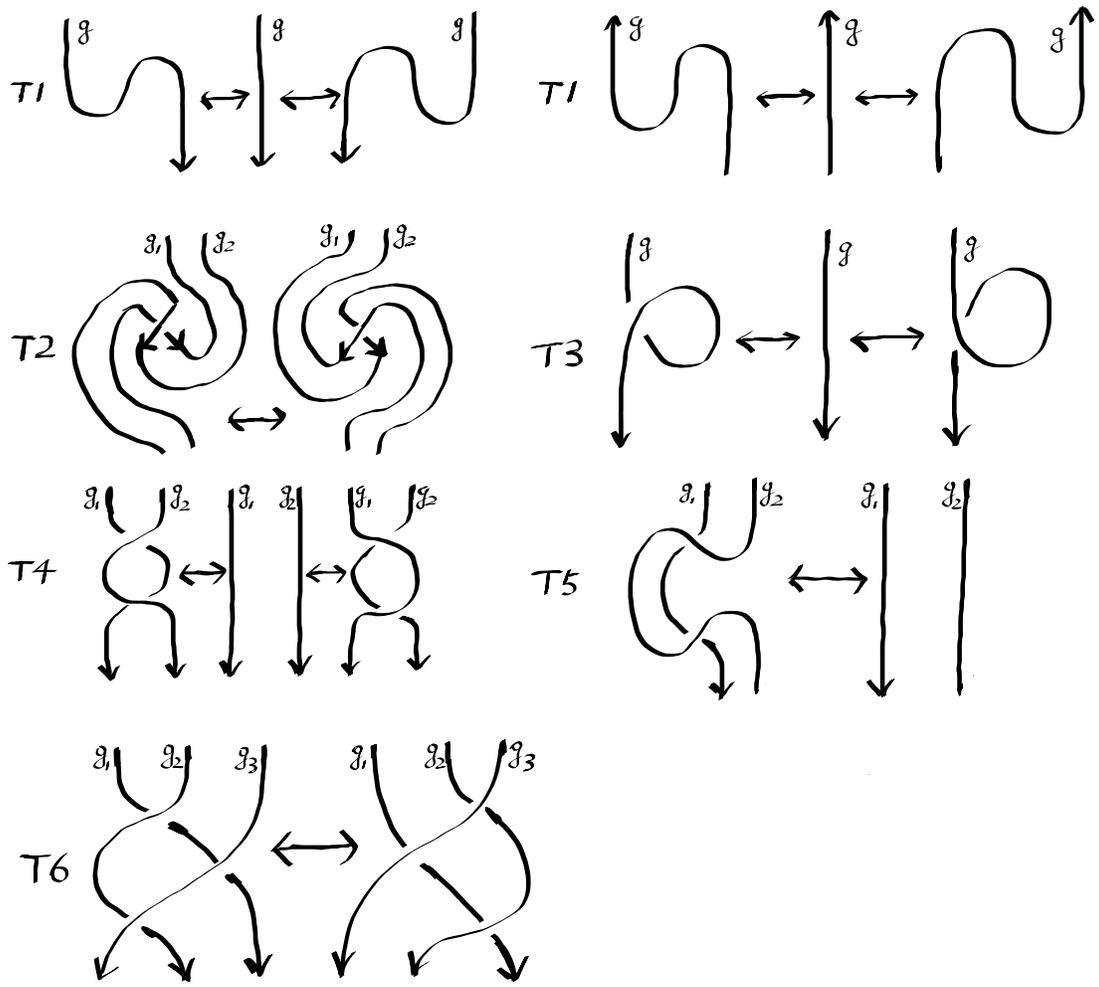}
  \caption{Turaev moves }
\end{figure}

So if one want to construct a representation of the tangle category $\mathcal{T}_{G}$, or, a functor $F: \mathcal{T} _{G } \rightarrow \mathscr{C}$ for some tensor category $\mathscr{C}$, one only need to do it according to the following procedures:

\begin{enumerate}
\item For every object $[(g, \epsilon ) ] $ of $\mathcal {T}_{G} $,  associate an object $F( (g , \epsilon ) )$ of $\mathscr{C}$;

\item For $X^{+}_{g_1 , g_2 }$ , associate  $c^{+} _{g_1 , g_2} \in Hom ( F( [(g_1 , +1), (g_2 , +1 ) ]  ) ,  F([(g_2 , +1 ), (g_1 , +1) ]) ) $;

 For $X^{-} _{g_1 ,g_2}$ , associate  $c^{-} _{g_1 , g_2} \in Hom ( F( [(g_1 , +1), (g_2 , +1 ) ]  ) ,  F([(g_2 , +1 ), (g_1 , +1) ]) ) $;

For $\cap ^{r} _{g} $, associate $b^{r} _{g} \in Hom ( F( \emptyset  ), F( [( g,-1  )  ]    ) \otimes F( [ ( g, +1  )  ]  )    ) $ ;

For $\cap ^{l} _{g} $ , associate  $b^{l} _{g} \in  Hom ( F( \emptyset  ), F( [( g,+1  )  ]    ) \otimes F( [ ( g, -1  )  ]  )    ) $;

For $\cup ^{r} _{g} $  , associate  $d^{r} _{g} \in Hom ( F( [( g,+1  )  ]    ) \otimes F( [ ( g, -1  )  ]  ) , F(\emptyset )  ) $ ;

For $\cup ^{l}_{g} $, associate a  $d^{l} _{g} \in Hom ( F( [( g,-1  )  ]    ) \otimes F( [ ( g, +1  )  ]  ) , F(\emptyset )  ) $ ;

\item Make sure that above morphisms $c^{+} _{g_1 , g_2  } , c^{-} _{g_1 , g_2 } , b^{r} _{g} , b^{l} _{g} , d^{r} _{g} , d^{l} _g   $   satisfy all relations represented by the Turaev moves $T1-T6$ .
\end{enumerate}

 Then by Theorem \ref{turaevtangle}, above association would extend to a functor $F : \mathcal{T}_{G} \rightarrow \mathscr{C} $.

 Now we construct a invariant for $G-$colored tangles, or a functor $ F: \mathcal{T}_{G} \rightarrow \mathcal{V}$ where
 $\mathcal{V}$ is the tensor category of complex linear spaces and linear maps, which is one of the main results of this paper. First recall the maps
\begin{align*}
 &\Psi _{g_1 , g_2 }: G\times G \rightarrow G\times G : (x,y ) \mapsto (yx^{-1} g_1 ^{-1} x, x  ), \\
&\Psi ^{-} _{g_1 , g_2} : G \times G \rightarrow G \times G  : (x,y) \mapsto (y, xy^{-1} g_2 y  ) .
\end{align*}

 Set $ N_{G} : \mathbb{C} \rightarrow \mathbb{C} G \otimes \mathbb{C} G$ to be the linear map $N_{G} (\lambda )= \lambda \sum_{g\in G} g\otimes g $;
 set $U_{G} : \mathbb{C} G \otimes \mathbb{C} G \rightarrow \mathbb{C} $ to be the linear map $U_{G} (g\otimes h ) = \delta _{g,h} , \forall g,h\in G$.

 \begin{thm}
 \label{thm:colored link invariant}
 The following association defines a functor $F: \mathcal{T}_{G} \rightarrow \mathcal{V} $. Where $\mathcal{V}$ is the category of vector spaces and linear maps.
\begin{enumerate}
\item[ (1)]  $F ( [(g , \epsilon )  ]  ) = \mathbb{C} G  $ for any $g\in G $ and any $\epsilon \in \{ +1 , -1 \} $; $F (\emptyset ) = \mathbb{C} $;

\item[ (2)] $F( X^{+} _{g_1 ,g_2 } ) = c^{+} _{g_1 ,g_2 }=( id_{\mathbb{C} G}  \otimes (g_1 )_{L} ) \Psi _{g_1 , g_2 } $;

\item[ (3)] $F( X^{-} _{g_1 ,g_2 } ) = c^{-} _{g_1 , g_2} = ( ( g_2 ^{-1} )_{L} \otimes id_{\mathbb{C} G}  ) \Psi ^{-} _{g_1 , g_2 }  $.

\item[ (4)] $F (\cap ^{r} _{g} ) = b^{r} _{g} = F( \cap ^{l} _{g} ) = b^{l} _{g} = N_{G} $;

 \item[(5)] $F(\cup ^{r} _{g} )= d^{r} _{g} = F(\cup ^{l} _{g} ) = d^{r} _{g} = U_{G} $.
\end{enumerate}
 \end{thm}

 \begin{pf} The proof is to check all the Turaev moves. Since the computations go straightforward without any difficulty we omit the details.
 Notice the key relation, that is ,the Yang-baxter relation is just content of Theorem \ref{thm:basicbraidrelation}.

 \end{pf}

 \begin{rem}
 Thus we obtain a invariant for oriented , $G-$colored tangles. In (3), (4) of Theorem \ref{thm:colored link invariant}, we see
 $F( \cap ^{l} )= F(\cap ^{r} ) $ and $ F (  \cup ^{l}) = F( \cup ^{r} ) $. It might imply that this tangle invariant is essentially a invariant
 for unoriented tangles.  It isn't hard to see this invariant generalize the link invariant defined in Section 6.
 \end{rem}

In very similar way we define a invariant for oriented extended pairs-colored links as follows. It is a natural extension of above construction, since
$G$ can be viewed as a subset of $\mathscr{E}G$ naturally. First, we define a category $\mathcal{T}_{\mathscr{E}G}$ of oriented $\mathscr{E}G$-colored tangles ,and  basic $\mathscr{E}G$-colored tangles $X^{+}_{(g_1 , b_1) ,(g_2 , b_2)} $, $X^{-}_{(g_1 , b_1) ,(g_2 , b_2)} $, $\cap^{r} _{(g,b)}$, $ \cap^{l} _{(g,b)} $,$\cup^{r} _{(g,b)}  $ and $\cup^{l} _{(g,b)}  $,  similar with the construction of the category $\mathcal{T}_{G}$ and those basic tangles.

\begin{thm}
 \label{thm:extended colored link invariant}
 The following association defines a functor $F: \mathcal{T}_{\mathscr{E}G} \rightarrow \mathcal{V} $. Where $\mathcal{V}$ is the category of vector spaces and linear maps.
\begin{enumerate}
\item[ (1)]  $F ( [((g,b) , \epsilon )  ]  ) = \mathbb{C} G  $ for any $(g,b)\in \mathscr{E}G $ and any $\epsilon \in \{ +1 , -1 \} $; $F (\emptyset ) = \mathbb{C} $;

\item[ (2)] $F( X^{+} _{(g_1 ,b_1 ) ,(g_2 ,b_2 ) } ) = c^{+} _{(g_1 ,b_1 ) , (g_2 ,b_2 ) }= ( id_{\mathbb{C}G} \otimes (b_1 g_1 )_{L}   ) \Psi _{g_1 , g_2 }  $;

\item[ (3)] $F(  X^{-} _{(g_1 ,b_1 ) ,(g_2 ,b_2 ) } ) = c^{-} _{(g_1 ,b_1 ) , (g_2 ,b_2 )} = (   (b_2 g_2 ^{-1} )_{L} \otimes id_{\mathbb{C}G} ) \Psi ^{-} _{g_1 , g_2 } $.

\item[ (4)] $F (\cap ^{r} _{(g,b)} )  = b^{r} _{g} = ( (b)_{L} \otimes id_{\mathbb{C}G}   ) N_{G}    $ ;

\item[ (5)] $F( \cap ^{l} _{(g,b))} ) = b^{l} _{g} = N_{G} $;

 \item[(6)] $F(\cup ^{r} _{(g,b)} )= d^{r} _{g} = U_{G} ( (b)_{L} \otimes id_{\mathbb{C}G}  ) $;

 \item[ (7)] $   F(\cup ^{l} _{g} ) = d^{r} _{g} = U_{G} $.
\end{enumerate}
 \end{thm}

 \begin{pf}
 The proof is also by checking all those Turaev moves.
 \end{pf}
 \\

 If restrict above tangle invariant to the set of $\mathscr{E}G$ colored links such that all components are colored with the same $(g,b)$, then we obtain the uncolored link invariant $\Lambda_{G,[(g,b)]}(L)$ as in Remark \ref{rem:extendableinvariant}.\\

 Next we consider the case for framed tangles. For a precise definition of framed tangles and examples please consult Chapter 3 of \cite{Oh}. Every framed tangle
 can also be presented by the same tangle diagrams as above, and from every tangle diagram  we recover a framed tangle in a canonical way.

\begin{defi}
\label{defi:framedGtangle} By a oriented, $G-$colored tame framed tangle, we mean a framed tangle( or a ribbon tangle) $T$ satisfying the following
conditions.

(1) $T\subset \mathbb{C} \times [0,1]$;

(2) Every component of $T$ is associated with an element of $G$;

(3) The upper boundary  of $T$ is  $([\frac{2}{3},1]\cup [\frac{5}{3},2 ]\cup ...\cup [k-\frac{1}{3},k ])\times \{ 1 \}  $ for some positive integer $k$ or empty. The lower boundary
is $ ([\frac{2}{3},1]\cup [\frac{5}{3},2 ]\cup ...\cup [l-\frac{1}{3},l ]) \times \{ 0 \}  $ for some positive integer $l$ or empty.

(4) Every component of $T$ has $2n \pi $ total twist for $n\in \mathbb{Z}$.

If $T$ is a oriented $G-$ colored tame framed tangle, suppose the upper boundary and the lower boundary of $T$ are
 $([\frac{2}{3},1]\cup [\frac{5}{3},2 ]\cup ...\cup [k-\frac{1}{3} ,k ])\times \{ 1 \} \times \{ 1 \}  $ and
$ ([\frac{2}{3},1]\cup [\frac{5}{3},1 ]\cup ...\cup [l-\frac{1}{3},l ])\times \{ 1 \} \times \{ 0 \}  $ respectively, define
$\partial ^{+} T = [(g_1 , \epsilon _1  ), (g_2 , \epsilon _2 ),..., (g_k , \epsilon _k )   ] $ and
$\partial ^{-} T = [(h_1 , \mu _1 ), (h_2 , \mu _2 ),..., (h_l , \mu _l ) ] $  (  where $g_i , h_j \in G$ and $ \epsilon _i , \mu _j \in \{ +1 ,-1 \} $) according to the following rules:

(a) The component of $T$ passing the boundary point $\{ i \} \times \{ 1 \} $ is colored by $g_i$ for $1\leq i\leq k $; the component of $T$ passing the boundary point $\{ j \} \times \{ 0 \}  $ is colored by $h_i $ for $1\leq j \leq l $;

(b) Taking into consideration of the orientation of components, if a component is leaving (entering ) a upper boundary point $\{ i \} \times \{ 1 \} $, then $\epsilon _i = +1 (-1 )$; if a component is entering ( leaving ) a lower boundary point $\{ j \} \times \{ 0 \} $, then $\mu _j = +1 (-1 )  $.

\end{defi}

Now one can define a category $\mathbb{T}_{G} $ of framed tangles as follows. The objects of $\mathbb{T}_{G}$ is the same as the objects of the category
$\mathcal{T}_{G}$. For any two objects $X_1 $ and $X_2$, a morphism from $X_1 $ to $X_2$ is a isotopy class of oriented, $G-$colored, tame framed tangle
$T$ such that $\partial ^{+} T = X_1  $ and $\partial ^{-} T = X_2 $.

Another category $\mathcal{D} \mathbb{T}_{G}$ is defined as follows.  The objects of $\mathcal{D} \mathbb{T} _{G}$ are as same as the objects of the category
$\mathbb{T}_{G}$. For two objects $X_1$ and $X_2$, a morphism in the category $\mathcal{D}\mathbb{T}_{G}$ is an equivalent class of  a sliced, oriented , $G-$colored  framed  tangle $T$ such that $\partial ^{+} T = X_1 $and $ \partial ^{-} T = X_2 $, with equivalence relations generated by relations represented by the Turaev moves $T1, T2 , T4 ,T5 ,T6$  in Figure 3 and $T3^{'}$ in Figure 4.

There is a natural functor $\mathcal{S}: \mathcal{D}\mathbb{T}_{G} \rightarrow \mathbb{T}_{G} $ by mapping a tangle diagram to the framed tangle represented by it.

\begin{thm}
The functor $\mathcal{S}$ is a isomorphism. Equivalently, the functor $\mathcal{S}$ induces a bijection between sets of morphisms of these two category.
\end{thm}

\begin{figure}[htbp]

  \centering
    \includegraphics[height=3.5cm]{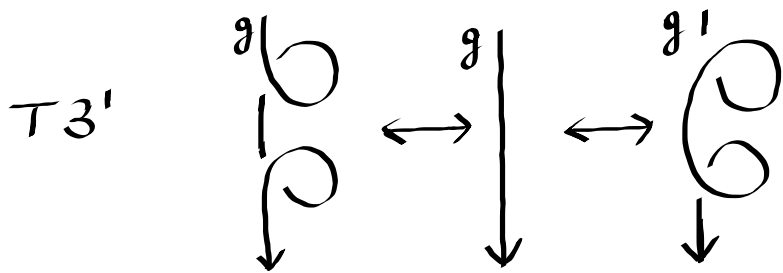}
  \caption{Ribbon Tureav moves }
\end{figure}

Let $\Psi _{g_1 ,g_2}$, $\Psi ^{-} _{g_1 , g_2 }$, $N_{G}$ and $U_{G}$ be defined as before Theorem \ref{thm:colored link invariant}.

 \begin{thm}
 \label{thm:colored framed link invariant}
 The following association defines a functor $\mathbb{F}: \mathbb{T}_{G} \rightarrow \mathcal{V} $. Where $\mathcal{V}$ is the category of vector spaces and linear maps.
\begin{enumerate}
\item[ (1)]  $\mathbb{F} ( [(g , \epsilon )  ]  ) = \mathbb{C} G  $ for any $g\in G $ and any $\epsilon \in \{ +1 , -1 \} $; $F (\emptyset ) = \mathbb{C} $;

\item[ (2)] $\mathbb{F}( X^{+} _{g_1 ,g_2 } ) = c^{+} _{g_1 ,g_2 }= \Psi _{g_1 , g_2 } $;

\item[ (3)] $\mathbb{F}( X^{-} _{g_1 ,g_2 } ) = c^{-} _{g_1 , g_2} = \Psi ^{-} _{g_1 , g_2 }$.

\item[ (4)] $\mathbb{F} (\cap ^{r} _{g} ) = b^{r} _{g} = \mathbb{F}( \cap ^{l}_{g} ) = b^{l} _{g} = N_{G} $;

\item[ (5)] $\mathbb{F}(\cup ^{r} _{g} )= d^{r} _{g} = \mathbb{F}(\cup ^{l}_{g} ) = d^{r} _{g} = U_{G} $.
\end{enumerate}
 \end{thm}

 \begin{pf}
  We only need to show the relations represented by the Turaev moves $T1,T2,T4,T5,T6$ and $T3^{'}$ are satisfied. Since the computation is straightforward we omit.
 \end{pf}\\

Now we define a category $\mathbb{T}_{\mathscr{E}G}$ by replacing $G$ with $\mathscr{E}G$, that is ,the category of oriented ribbon tangles colored with extended pairs. We introduce similar basic $\mathscr{E}G$-colored tangles $X^{+}_{(g_1 , b_1) ,(g_2 , b_2)} $, $X^{-}_{(g_1 , b_1) ,(g_2 , b_2)} $, $\cap^{r} _{(g,b)}$, $ \cap^{l} _{(g,b)} $,$\cup^{r} _{(g,b)}  $ and $\cup^{l} _{(g,b)}  $.  We define a invariant for this kind of tangles in the following theorem.

\begin{thm}
 \label{thm:extended colored framed link invariant}
 The following association defines a functor $\mathbb{F}: \mathbb{T}_{\mathscr{E}G} \rightarrow \mathcal{V} $. Where $\mathcal{V}$ is the category of vector spaces and linear maps.
\begin{enumerate}
\item[ (1)]  $\mathbb{F} ( [((g,b) , \epsilon )  ]  ) = \mathbb{C} G  $ for any $(g,b)\in \mathscr{E}G $ and any $\epsilon \in \{ +1 , -1 \} $; $F (\emptyset ) = \mathbb{C} $;

\item[(2)] $\mathbb{F}( X^{+} _{(g_1 ,b_1 ) ,( g_2 ,b_2 ) } ) = c^{+} _{g_1 ,g_2 }= \Psi _{g_1 , g_2 } $;

\item[ (3)] $\mathbb{F}( X^{-} _{( g_1 ,b_1 ) ,(g_2 ,b_2 ) } ) = c^{-} _{g_1 , g_2} = \Psi ^{-} _{g_1 , g_2 }$.

\item[ (4)] $\mathbb{F} (\cap ^{r} _{(g,b)} ) = b^{r} _{(g,b)} =( (b)_{L} \otimes id_{\mathbb{C}G}  )N_{G}  $ ;

\item[  (5)] $ \mathbb{F}( \cap ^{l}_{(g,b)} ) = b^{l} _{(g,b)} = N_{G} $;

\item[ (6)] $\mathbb{F}(\cup ^{r} _{(g,b)} )= d^{r} _{g} =U_{G} ( (b)_{L} \otimes id_{\mathbb{C}G}  )$;

\item[  (7)] $\mathbb{F}(\cup ^{l}_{(g,b)} ) = d^{r} _{(g,b)} = U_{G} $.
\end{enumerate}
 \end{thm}

\begin{pf}
 We only need to show the relations represented by the Turaev moves $T1,T2,T4,T5,T6$ and $T3^{'}$ are satisfied.
\end{pf}\\

Next we show above two type of invariants for tangles and framed tangles are essentially invariants for tangles colored by conjugacy classes of $\mathscr{E}G$. We
only show the unframed case. The other case can be shown similarly. Let $T$ be a
 morphism in the category $\mathcal{T}_{\mathscr{E}G}$, equivalently, an oriented, $G-$colored, tame tangle. Suppose $T=C_1 \sqcup C_2 \sqcup ...\sqcup C_m $ where
 $C_1 ,C_2 ,..., C_m $ are components of $T$. Suppose the component $C_i$ is colored by $(h_i ,b_i ) \in \mathscr{E}G$. Now chose randomly elements
 $\alpha_1 , \alpha_2 ,..., \alpha_m \ in\  G$, let $T^{'}$ be the same oriented tangle as $T$, but the component $C_i$ of $T^{'}$ are colored by
 $( \alpha _i  h_i \alpha _i ^{-1} , \alpha _i  b_i \alpha _i ^{-1} ) $.  Suppose
 \begin{align*}
 &\partial^{+} T =[( g_1 ,c_1 ) , \epsilon _1  ), ( ( g_2 ,c_2 ) , \epsilon _2 ),..., ( (g_k ,c_k ) , \epsilon _k )   ]  ,\  and\\
 &\partial^{-} T =  [( ( g^{'} _1 ,c^{'} _{1} ) , \mu _1 ), ( (g^{'} _2 ,c^{'} _{2}) , \mu _2 ),..., ( (g^{'} _l , c^{'} _{l} , \mu _l ) ] .
\end{align*}
  Certainly we have
 $\{ ( g_1 , c_1 ) ,..., (g_k ,c_k ) , (g^{'} _1 ,c^{'} _{1}) ,..., (g^{'} _{l} ,c^{'} _{l} ) \} \subset \{ ( h_1 ,b_1 ) ,( h_2 ,b_2 ) ,..., ( h_m ,b_m ) \} $. Now let $F$ be the functor defined in
 Theorem\ref{thm:colored link invariant}.

 \begin{thm}
 The following diagram is commutative

\[ \begin{CD}
  \underbrace{G\times G\times \cdots \times G }_{k}   @>f_1>>  \underbrace{G\times G\times \cdots \times G }_{k} \\
 @V\mathbb{F}( T )VV        @VV\mathbb{F}(T^{'})V \\
   \underbrace{G\times G\times \cdots \times G }_{l}    @>f_2>>   \underbrace{G\times G\times \cdots \times G }_{l}
 \end{CD}
 \]

  Suppose the upper boundary point $\{ v \} \times \{ 1 \} $ belongs to the component $T_{i_v }$ for $1\leq v\leq k$, and the lower boundary point
   $\{ w \} \times \{ 0 \} $ belongs to the component $T_{j_w } $ for $1\leq w\leq l$ .  We set\begin{align*}
   &f_1 ( (x_1 , x_2 ,..., x_k  )  )= ( \alpha_{i_1 } x_1 , \alpha_{i_2 } x_2 , ..., \alpha_{i_{k}} x_k  ) ,\\
   &f_2 ( (y_1 , y_2 ,..., y_l ) ) = (\alpha_{j_1 } y_1 , \alpha_{j_2 } y_2 ,..., \alpha_{j_{l}} y_l ) .
   \end{align*}
   Besides, if $k=0$, set $f_1 = id_{\mathbb{C}} : \mathbb{C} \rightarrow \mathbb{C}$. If $l=0$, set $f_2 = id_{\mathbb{C}}$ likewise.

\end{thm}
\begin{pf}
Since $\mathbb{F}$ is a functor, and every tangle is composed by the elementary tangles in Figure 2, we only need to prove this theorem for the cases when $T$ is one of
the six elementary tangles. Taking the first one as a example,

$\mathbb{F}(T^{'} ) \circ f_1  ((x_1 , x_2) ) = \mathbb{F}(T^{'} ) ( (\alpha_1 x_1 , \alpha_2 x_2    )  ) =
(\alpha_2 x_2 , x_1 \alpha_1 ^{-1} \alpha_1 g_1 ^{-1} \alpha_1 ^{-1} \alpha_1 x_1 , \alpha_1 x_1  )$

$= (\alpha_2 x_2  x_1 ^{-1} g_1 ^{-1} x_1 , \alpha_1 x_1 ).  $

$f_2 \circ \mathbb{F}(T) ((x_1 , x_2 ) )= f_2 ( x_2 x_1 g_1 ^{-1} x_1 , x_1  )= ( \alpha_2 x_2 x_1 g_1 ^{-1} x_1 , \alpha_1 x_1 ) .$

When $T$ is the elementary tangle $U^{r} _{(g,b)}$, for any $x,y\in G$,

$ \mathbb{F} (T^{'} ) \circ f_1 ( x\otimes y  )= \mathbb{F} ( T^{'} ) ( \alpha x \otimes \alpha y  )= \delta _{\alpha  b \alpha ^{-1} \alpha x , \alpha y     } =\delta _{\alpha  b  x , \alpha y     }  =  \delta _{bx,y} \in \mathbb{C}$,

$f_2 \circ \mathbb{T} ( x\otimes y )= id_{\mathbb{C}} ( \delta_{ bx,y  }   )= \delta_{bx,y}   $.

\end{pf}
$\\$
Apply above theorem to any oriented,  $G-$colored link $L$, which is a morphism from $\emptyset$ to $\emptyset$ (Recall $\mathbb{F}( \emptyset ) = \mathbb{C}$ ),
 we have the following corollary.

\begin{cor}
\label{cor:essentialconjugacy}

If $L$ and $L^{'}$ are the same oriented link but with possibly different colors ( in $\mathscr{E}G$ ) , such that for every component $L_i$ of $L$, the color of $L_i$ in  $T$  is conjugate to the color of $L_i $ in $T^{'}$, then $F( T )= F(T^{' } ) \in \mathbb{Z} ^{\geq 0} $.

\end{cor}

We end this section with the following theorem. It holds only for $G-$colored tangles but not for $\mathscr{E}G-$colored tangles.

\begin{thm}
\label{thm:reversing}
Let $T$ be $G-$colored oriented $(k,l)$ tangle. Suppose $T= T_1 \sqcup T_2 \sqcup \cdots \sqcup T_m $ where $T_i $'s are components of $T$. Choose any component $T_{j}$, suppose the color of $T_{j}$ is $h$.  Suppose $T^{'}$ is the oriented $G-$colored tangle obtained from $T$ by changing the orientation of $T_{j}$, and change the color of $T_j $ to $h^{-1}$ at the same time. Then
$$\mathcal{F}(T) = \mathcal{F}(T^{'}) : \underbrace{G\times G\times ...\times G}_{k} \rightarrow \underbrace{G\times G\times ...\times G}_{l }   .$$

The same statement holds for framed tangles and the functor $\mathbb{F}$.
\end{thm}

\begin{pf}
Because every tangle is composed by the elementary tangles in Figure 2, we only need to prove the theorem for elementary tangles. When $T= X^{-} _{g_1 , g_2 }$, Suppose $X^{-} _{g_1 \downharpoonright , g_2 \upharpoonright  }$ and $X^{-} _{g_1 \upharpoonright , g_2 \downharpoonright } $ are the tangles obtained by reversing the orientation of one of the two arcs of $X^{-} _{g_1 , g_2 }$.  Beside $X^{-} _{g_1 \downharpoonright , g_2 \upharpoonright  }$  and
 $X^{-} _{g_1 \upharpoonright , g_2 \downharpoonright } $ are sliced tangles representing them. By using these sliced tangles, we have
 \begin{align*}
 &\mathcal{F}( X^{-} _{g_1 \downharpoonright , g_2 \upharpoonright  } ) (x_1 \otimes x_2 ) \\
 &= (id_{G} \otimes id_{G} \otimes U_{G}  )\circ ( id_{G} \otimes \mathcal{F} ( X^{+} _{g_2 , g_1 } )  \otimes id_{G} )\circ (N_{G} \otimes id_{G} \otimes id_{G} ) (1\otimes x_1 \otimes x_2 )\\
 &= (id_{G} \otimes id_{G} \otimes U_{G}  )\circ ( id_{G} \otimes \mathcal{F} ( X^{+} _{g_2 , g_1 } )  \otimes id_{G} )(\sum_{x\in G} x\otimes x \otimes x_1 \otimes x_2 ) \\
 &= (id_{G} \otimes id_{G} \otimes U_{G}  )( \sum_{x\in G} x\otimes x_1 x^{-1} g^{-1}_{2} x  \otimes x \otimes x_2  ) \\
 &= x_2  \otimes x_1 x^{-1} _{2} g^{-1} _{2} x_2   = \mathcal{F} (X^{-} _{g_1 , g_2 ^{-1} } ) (x_1 \otimes x_2 ) .
 \end{align*}
 \begin{align*}
 &\mathcal{F} ( X^{-} _{g_1 \upharpoonright , g_2 \downharpoonright } )(x_1 \otimes x_2 ) \\
 &=(U_{G} \otimes id_{G} \otimes id_{G})\circ (id_{G} \otimes \mathcal{F}(X^{+} _{g_2 , g_1 } ) \otimes id_{G} ) \circ (id_{G} \otimes id_{G} \otimes N_{G} )
 (x_1 \otimes x_2 \otimes 1 ) \\
 &= (U_{G} \otimes id_{G} \otimes id_{G})\circ (id_{G} \otimes \mathcal{F}(X^{+} _{g_2 , g_1 } ) \otimes id_{G} )( \sum_{x\in G} x_1 \otimes x_2 \otimes x\otimes x) \\
 &=(U_{G} \otimes id_{G} \otimes id_{G}) (\sum_{x\in G } x_1 \otimes xx^{-1} _{2} g_2 x_2   \otimes x_2  \otimes x  ) = x_2 \otimes x_1 x^{-1} _2 g_2 x_2 \\
 &=\mathcal{F} ( X^{-} _{g_1 ^{-1} , g_2 }  ) ( x_1 \otimes x_2  ) .
 \end{align*}

 \begin{figure}[htbp]

  \centering
    \includegraphics[height=5.5cm]{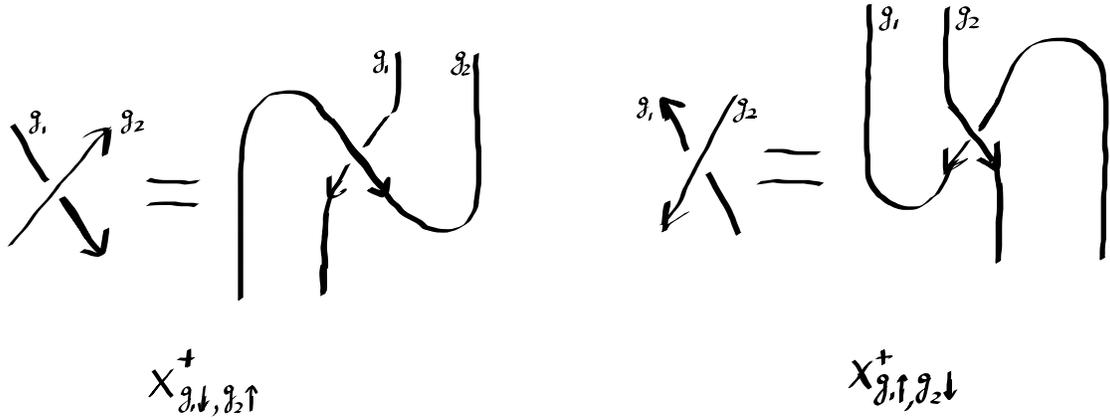}
  \caption{One string reversed crossings}
\end{figure}

The cases for other elementary tangles can be proved similarly.

\end{pf}

\section{Group invariants}

Now we restrict to the cases of links, and show these invariants for all kinds of finite groups $G$ are actually dominated by certain group, which can be
proved to be a invariant also. Let $L= L_1 \sqcup L_2 \sqcup \cdots \sqcup L_m  $ be a oriented link with $m$ components.  Let $b \in B_{n}$ be a braid whose closure $\hat{b}$ presents $L$. Let $G$ be a finite group, suppose $\bar{L}$ is a $G-$colored link based on $L$ by assign a color $g_i \in G$ to the component $L_i$ of $L$, for $1\leq i\leq m$. Suppose $\bar{b}$ is a $G-$colored braid based on $b$, whose closure is $\bar{L}$.

\begin{lem}
$\mathcal{F} ( \bar{L} ) = tr ( \mathcal{F} ( \bar{b})   ) .$
\end{lem}
\begin{pf}
It is proved by slicing the closure $\hat{ \bar{b} }$, and using the fact that $\cap ^{\pm } _{g } $ and $\cup ^{\pm} _{g}$ are represented by such simple maps
$N_G$ and $U_G$ respectively.
\end{pf}\\

Now, the isomorphism $\mathcal{F} (b) : \underbrace{\mathbb{C} G \otimes \mathbb{C} G \otimes \cdots \otimes \mathbb{C} G  }_{n} \rightarrow
\underbrace{\mathbb{C} G \otimes \mathbb{C} G \otimes \cdots  \otimes \mathbb{C} G  }_{n} $ is actually the linearization of the bijective map
$\mathcal{F}(b) : \underbrace{ G \times G \times \cdots \times G  }_{n} \rightarrow \underbrace{G \times G \times \cdots \times G  }_{n}   $.

Now we associate a pair of groups to $b $. We introduce variables $X_1 , X_2 ,..., X_n , C_1 , C_2 ,..., C_m $, and let $ F_{n+m}$ be the free group generated by them. Then we color the
$i-$th string of $b$ with $C_{n_i }$, if it (as a subset of $L= \hat{b}$) belongs to the component $L_{n_i}$. Denote the resulted $F_{n+m}-$colored braid as $b^{''}$. Recall the  map $\Psi_{g_1 , g_2}$ and $\Psi^{-} _{g_1 , g_2} $
before Theorem \ref{thm:colored link invariant} and (2),(3) of Theorem \ref{thm:colored link invariant}. For a $F_{n+m}-$colored $n-$string simple braid
$X^{+} _{i; C_{j^{'} } , C_{j^{''} } }$ or $X^{-} _{i;  C_{j^{'}} , C_{j^{''} } }  $, whose $i-$th string and $(i+1)-$th string are colored by $C_{j^{'}}$ and $C_{j^{''}}$ respectively, and $(Z_1 , Z_2 ,..., Z_n )\in \underbrace{ F_{n+m} \times F_{n+m} \times ... \times F_{n+m} }_{n}$,    we set
\begin{align*}
&\mathcal{F}^{'} ( X^{+} _{ i; C_{j^{'} } , C_{j^{''}} }  )(\cdots , Z_i , Z_{i+1} , \cdots    )=  (\cdots , Z_{i+1} Z_i ^{-1} C_{j^{'} }^{-1} Z_i , C_{j^{'} } Z_i ,\cdots   ) , \\
&\mathcal{F} ^{'} ( X^{-} _{i; C_{j^{'} } , C_{j^{''} }} )(\cdots ,Z_i ,Z_{i+1} , \cdots )= (\cdots , C_{j^{''} }^{-1} Z_{i+1} , Z_i Z_{i+1} ^{-1} C_{j^{''} }Z_{i+1} ,\cdots )  .\\
&\mathcal{F}^{''} ( X^{+} _{ i; C_{j^{'} } , C_{j^{''}} }  )(\cdots , Z_i , Z_{i+1}, \cdots    )=  (\cdots , Z_{i+1} Z_{i} ^{-1} C_{j^{'} }^{-1} Z_{i}, Z_{i} ,\cdots   ) , \\
&\mathcal{F} ^{''} ( X^{-} _{i; C_{j^{'} } , C_{j^{''} }} )(\cdots ,Z_i ,Z_{i+1} , \cdots )= (\cdots ,  Z_{i+1} , Z_{i} Z_{i+1} ^{-1} C_{j^{''} }Z_{i} ,\cdots )  .
\end{align*}
If $b^{''} = \Pi _{t=1} ^{N} X^{\epsilon_{t}} _{i_t ; C_{j^{'}_{t}} , C_{t, j^{''}_{t} }  }$, we set
\begin{align*}
&\mathcal{F}^{'} (b^{''} )= \mathcal{F}^{'} ( X^{\epsilon_{N}} _{i_N ; C_{j^{'} _{N}} , C_{ j^{''} _{N} }  } ) \circ \cdots \circ
\mathcal{F } ^{'} (  X^{\epsilon_{1}}   _{i_1 ; C_{j^{'} _{1}} ,C_{j^{''} _1 }  } ).\\
&\mathcal{F}^{''} (b^{''} )= \mathcal{F}^{''} ( X^{\epsilon_{N}} _{i_N ; C_{j^{'} _{N}} , C_{ j^{''} _{N} }  } ) \circ \cdots \circ
\mathcal{F } ^{''} (  X^{\epsilon_{1}}   _{i_1 ; C_{j^{'} _{1}} ,C_{j^{''} _1 }  } ) .
\end{align*}
 We define set of words $W^{b} _1 ,...,W^{b} _n \in F_{n+m }$ and $U^{b} _1 ,...,U^{b} _n \in F_{n+m }   $ by
\begin{align*}
&\mathcal{F}^{'} (b^{''} ) ( (X_1 ,X_2 ,..., X_n   )  )= (W^{b} _1 , W^{b} _2 ,..., W^{b} _n   ) ,\\
&\mathcal{F}^{''} (b^{''} ) ( (X_1 ,X_2 ,..., X_n   )  )= (U^{b} _1 , U^{b} _2 ,..., U^{b} _n   ) .
\end{align*}
We can view $\mathcal{F}^{'}(b^{''})$ and $\mathcal{F}^{''}(b^{''})$ as group isomorphisms from $F_{n+m}$ to itself, by sending generators $(X_1 ,..., X_n , C_1 , ..., C_m)$ to $( (W^{b} _1 ,..., W^{b} _n   , C_1 , ..., C_m  ) $ and $(U^{b} _1 ,..., U^{b} _n  , C_1 ,..., C_m )$ respectively.
 To show the reliance of $W^{b} _i $ ($U^{b}_i $ ) on the generators, we denote them as $W^{b} _i ( X_1 , ..., X_n , C_1 , ..., C_n  ) $ (
 $U^{b} _i ( X_1 , ..., X_n , C_1 , ..., C_n  ) $ ).  If a group morphism $\phi $ map $F_{n+m} $ to some group $G$ such that $\phi (X_i ) = x_i \in G $ for $1\leq i\leq n$ and $\phi (C_j ) = b_j \in G$ for $1\leq j \leq n $, then naturally we write the image
$\phi (W^{b} _i )$ as $W^{b} _i (x_1 ,..., x_n , b_1 ,..., b_n  ) \in G  $ for $1\leq i \leq n$. It is easy to see
\begin{align*}
&\mathcal{F}(b) (( x_1 , ..., x_n ) )= (W^{b} _1 (x_1 ,..., x_n , b_1 ,..., b_n ),..., W^{b} _n (x_1 ,..., x_n , b_1 ,..., b_n  )   ) \  and\\
&\mathbb{F}(b) (( x_1 , ..., x_n ) )= (U^{b} _1 (x_1 ,..., x_n , b_1 ,..., b_n ),..., U^{b} _n (x_1 ,..., x_n , b_1 ,..., b_n  )  )   .
\end{align*}
This shows the meaning of these words in $F_{n+m}$. Return to the original $G-$colored braid $b$.   Since $\hat{b}=L$, every string of $b$ is colored by some element in $\{ g_1 ,..., g_m \} $. Suppose the $t-$th string of $b$ is colored by $ g_{i_t }$.

\begin{defi}
\label{defi:groupinvariantbraidversion}
Assuming above notations,  we define two groups $\mathcal{G}_{L, b}$ and $\mathbb{G}_{L,b}  $ (as quotient groups of $F_{n+m}$)  by the following presentations:
$$ \mathcal{G}_{L,b} :<X_1 ,X_2 ,..., X_n , C_1 ,C_2 ,..., C_m | X_i = W^{b}_i (X_1 ,..., X_n , C_{i_1},..., C_{i_m } )\ for\ 1\leq i\leq n  > .$$
$$ \mathbb{G}_{L,b} :<X_1 ,X_2 ,..., X_n , C_1 ,C_2 ,..., C_m | X_i = U^{b}_i (X_1 ,..., X_n , C_{i_1},..., C_{i_m } )\ for\ 1\leq i\leq n  > . $$
\end{defi}

\begin{prop}
\label{prop:groupinandnumberin}
We have:

(1) $\mathcal{F} (\bar{L}) = tr( \mathcal{F}( \bar{b}) )= |\{  f\in Hom (\mathcal{G}_{L,b} , G   )   | f(C_1 )= g_1 , \ f(C_2 )= g_2 ,..., \ f(C_m) = g_m \} |.$

(2) $\mathbb{F} (\bar{L}) = tr( \mathbb{F}(\bar{b}) )= |\{  f\in Hom (\mathbb{G}_{L,b} , G   )   | f(C_1 )= g_1 , \ f(C_2 )= g_2 ,..., \ f(C_m) = g_m \} |.$
\end{prop}
\begin{pf}
Since the linear map $\mathcal{F}(\bar{b})$ comes from the bijective map( permutation) $\mathcal{F}(\bar{b}) : \underbrace{ G \times G \times \cdots \times G  }_{n} \rightarrow \underbrace{G \times G \times \cdots \times G  }_{n}  $, so $tr (\mathcal{F}(\bar{b}) )$ equals the number of fixed points
of the map $\mathcal{F}(\bar{b})$. But now there is a bijection from the set of fixed points of $\mathcal{F}(\bar{b})$ to the set of group morphisms
$\{  f\in Hom (\mathcal{G}_{L,b} , G   )   | f(C_1 )= g_1 , \ f(C_2 )= g_2 ,..., \ f(C_m) = g_m \}$, so we finish the proof.
\end{pf}

So all our link invariants $\mathcal{F}(\bar{L} )$ constructed in last section from various finite groups $G$ are dominated by the data
$( \mathcal{G}_{L,b } ; C_1 , ..., C_m  ) $, that is , a group with a set of prescribed elements. Next we show this data is a link invariant itself.

\begin{thm}
\label{thm:groupinvariant1}
Suppose $d\in B_{N} $ is another braid such that $\hat{d} = L$. Suppose the set of generators of the group $\mathcal{G}_{L,d} $ in above construction is
$\{ X_{1} ^{'} ,..., X_{N} ^{'} , Z_1 , ..., Z_m \} $. Then there exists a group isomorphism $ \Phi : \mathcal{G}_{L, b} \rightarrow \mathcal{G}_{L,d} $ such that
$\Phi (Y_i ) = Z_{ \epsilon (i ) } $ for $1\leq i\leq m$ and some permutation $\epsilon \in S_{m}$.
\end{thm}
\begin{pf}
Since $\hat{d}=\hat{b}$, by Markov's theorem,  we only need to consider the following two cases: Case 1. The braid $b$ is transformed to $d$ by performing  Markov move I once; Case 2. The braid $b$ is transformed to $d$ by performing Markov move II once.

Case 1. Suppose $b^{'}$ is the $F_{n,m}-$colored braid obtained from $b$ by replacing the colour $g_i \in G$ with $Y_{i}$ as before.
Recall the notations in the beginning of section 7.   Suppose $\partial ^{+}(b^{'} ) = \partial ^{-} (b^{'})=  ( (Y_{i_1 } ,+ ), ..., ( Y_{i_n } ,+  )   ) . $
Since $b^{'}$ is a braid whose strings are oriented from upward to downward, the signals in the brackets are all $+$. Since $b^{'}$ can be closed to give a
$F_{n,m}- $coloured link, we must have $\partial ^{+}(b^{'} ) = \partial ^{-} (b^{'}) $. Now suppose $\alpha $ is a $n-$string $F_{n,m}$ coloured braid such that $\partial^{+} ( \alpha  )= ((Y_{i_1 } ,+ ), ..., ( Y_{i_n } ,+  )    ) $.  Denote the natural epimorphism $B_n \rightarrow S_n$ from the braid group to the permutation group as $\pi$, and suppose $\pi (\alpha )=s \in S_n $. Denote the mirror image of $\alpha $, which is also a $F_{n,m}-$colored braid, as $\alpha^{-1}$. Then we have

$\partial^{+} (\alpha^{-1} )= \partial ^{-} (\alpha )= ((Y_{i_{s(1)}} ,+) , ..., (Y_{i_{s(n)}},+)   ) $ and

$\partial^{-} (\alpha^{-1} )= ((Y_{i_1}, + ) ,..., (Y_{i_n } ,+)   ) $.

Denote $d^{'} = \alpha \cdot b^{'} \cdot \alpha^{-1} $ as a $F_{n,m}-$colored braid, then $d^{'}$ is a braid obtained from $b^{'}$ by performing Markov move I once. Any braid  $d^{'}$ as that can be written in this way for suitable braid $\alpha$. What we need to show is $\mathcal{G}_{L,b^{'}}\cong \mathcal{G}_{L, d^{'}}. $

We denote a term like $W^{\alpha} _i (W^{\beta} _1 ,...,W^{\beta} _n ; Y_{i_1} ,..., Y_{i_n }   ) $ simply as $ W^{\alpha} _i (W^{\beta}  ; Y_{i_1} ,..., Y_{i_n }   )  $ for convenience. Since $\alpha \cdot \alpha^{-1} =id_{B_n } $ and $\alpha^{-1} \cdot \alpha = id_{B_n }  $, we have

(A) $W^{\alpha} _i (W^{\alpha^{-1}} (X_1 , ..., X_n ;Y_{i_{s(1)}},...,Y_{i_{s(n)}}  ); Y_{i_1} ,..., Y_{i_n }  ) = X_i $ for $1\leq i\leq n$.

(B) $W^{\alpha^{-1}} _i (W^{\alpha} (X_1 ,..., X_n ; Y_{i_1} ,..., Y_{i_n }   ) ; Y_{i_{s(1)}},..., Y_{i_{s(n)}}  )= X_i $ for $1\leq i\leq n$.

Now by definition, the group $\mathcal{G}_{L, d^{'}}$ has the following presentation:
$$<X_1 ,...,X_n ,Y_1 ,..., Y_m | X_i = W^{\alpha \cdot b^{'} \cdot \alpha^{-1}} _i (X_1 ,...,X_n ;Y_{i_{s(1)}} ,...,Y_{i_{s(n)}} ) , (1\leq i\leq n) >. $$
Here we add new generators $V_i$ (corresponding to $W^{\alpha^{-1}} _i $ ), by Tietze theorem, above presentation is equivalent to the following
$$<X_i , V_i , Y_j  |  V_i = W^{\alpha^{-1}} _i (X_1 ,..., X_n ;Y_{i_{s(1)}},...,Y_{i_{s(n)}} ),
X_i = W^{\alpha \cdot d }(V_1 ,..., V_n ; Y_{i_1} ,...,Y_{i_n } ) > ,$$
where $1\leq i\leq n$ and $1\leq j\leq m $ in above presentation for all $i,j$.  Then by above identity (A) and (B), we see the set of relations
$ V_i = W^{\alpha^{-1}} _i (X_1 ,..., X_n ;Y_{i_{s(1)}},...,Y_{i_{s(n)}} ) (1\leq i\leq n) $ are equivalent to the set of relations
$X_i = W^{\alpha} _i (V_1 ,..., V_n ; Y_{i_1 }, ..., Y_{i_n }) (1\leq i\leq n ) $. So above presentation is equivalent to the following
$$< X_i , V_i , Y_j  | X_i = W^{\alpha} _i (V_1 ,..., V_n ; Y_{i_1 }, ..., Y_{i_n })  ,
X_i = W^{\alpha \cdot b^{'} }(V_1 ,..., V_n ; Y_{i_1} ,...,Y_{i_n } )   > ,$$
where $1\leq i\leq n$ and $1\leq j\leq m$.

By using Tietze theorem again, above presentation is equivalent to the following
$$<V_1 ,..., V_n , Y_1 ,..., Y_m | W^{\alpha} _i (V_1 ,..., V_n ; Y_{i_1 }, ..., Y_{i_n } ) = W^{\alpha} _i (W^{b^{'}} ;Y_{i_1} ,..., Y_{i_n} )
(1\leq i\leq n) >.$$
By using above identity (A) and (B) again, we see above presentation is equivalent to the following
$$< V_1 ,..., V_n , Y_1 ,..., Y_m | V_i = W^{b^{'}} _i ( V_1 ,..., V_n ; Y_{i_1 } ,..., Y_{i_n } )  (1\leq i\leq n )> ,$$
which is exactly a presentation of the group $\mathcal{G} _{L, b}$.

Case 2. Suppose the $n-$string braid $b , b^{'}$ are as in the case 1. Consider the braid $\sigma_{n} b^{'} $ , for which we color the $(n+1)$-th string also
by $Y_{i_n }$. Then the closure of $\sigma_{n} b^{'} $ is also $L$ and it is a braid obtained from $b^{'}$ by performing the  Markov move II once. We need to show $\mathcal{G}_{L , \sigma_n b^{'}} \cong \mathcal{G}_{L, b^{' }} $. By definition, the map $\mathcal{F}(\sigma_n b^{'})$ maps
$(X_1 , ..., X_{n-1} , X_n , X_{n+1}) $ to $(W^{b^{'}} _1 ,..., W^{b^{'} } _{n-1} , X_{n+1} (W^{b^{'} } _n  )^{-1} Y_{i_n } W^{b^{'}} _{n} , Y_{i_n } W^{b^{'}} _n    )   $, so the group $\mathcal{G}_{L, \sigma_n b^{'}}$ has the following presentation
$$<X_i ,Y_j | X_k = W^{b^{'}} _{k} , X_n = X_{n+1} (W^{b^{'}} )^{-1} Y_{i_n }W^{b^{'}} _n , X_{n+1} =Y_{i_n }W^{b^{'}} _n     >,$$
where $1\leq i\leq n+1 , 1\leq j\leq m , 1\leq k\leq n-1$. It is easy to see by canceling the generator $X_{n+1}$ using the last relation, above presentation is equivalent to the following presentation
$$<X_1 ,..., X_{n} , Y_1 ,..., Y_m | X_i = W^{b^{'}} _i ( 1\leq i\leq n)  > ,$$
which is just a presentation of the group $\mathcal{G}_{L, b}.$
\end{pf}$\\$

Interestingly the group invariant can be built on a link diagram directly without using braids. Still suppose $L= L_1 \sqcup L_2 \sqcup \cdots \sqcup L_m $ be a oriented $m-$components link. Suppose $D= D_1 \cup D_2 \cup \cdots \cup D_m $ is a link diagram representing $L$. Suppose $D_i $ is the part of $D$ representing the component $L_i$. Every $D_i$ is divided by the crossing points of $D$ into a set of connected simple curves  which we call short arcs. So each short arc of $D_i$ is simply a path from one crossing point on $D_i$ to the next. Denote the set of short arcs on $D_i$  as $\{ C^{1} _i , C^{2} _i ,..., C^{N_i } _i \} .$ Denote the set of crossing points of $D$ as $\{  P_1 ,P_2 ,..., P_{M}   \}$.

\begin{defi}
We denote the following local moves in Figure 9 that changing one link diagram into another as $\mathcal{R} I   , \mathcal{R} II $ and $\mathcal{R} III$ from left to right. They are the famous Reidemeister moves. We denote another move shown in Figure 10  as $\mathcal{R} I^{'} $. Since the link diagrams we considering are oriented, we orient the arcs in the following diagrams arbitrarily.
 \begin{figure}[htbp]

  \centering
    \includegraphics[height=3.5cm]{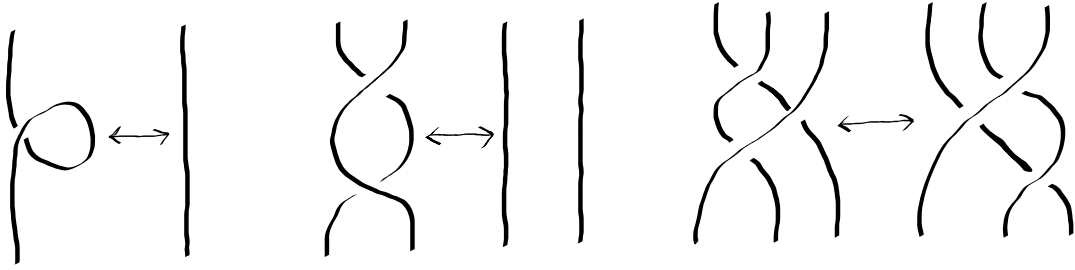}
  \caption{Deidemeister moves}
\end{figure}

 \begin{figure}[htbp]

  \centering
  \includegraphics[height=3.5cm]{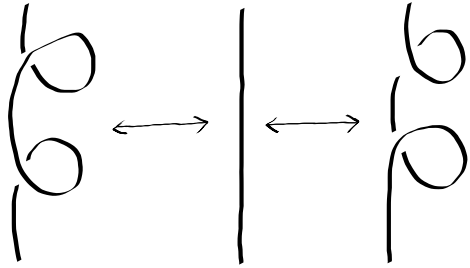}
  \caption{$\mathcal{R} I^{'} $ move }
\end{figure}
\end{defi}

\begin{thm}
Two oriented link diagrams $D_1 $ and $D_2 $ represent the same link if and only if $D_1$ can be transformed to $D_2$ by a sequence of Reidemeister moves
$\mathcal{R} I   , \mathcal{R} II $ and $\mathcal{R} III$.
\end{thm}

\begin{thm}
\label{thm:framedreidemeister}
Two oriented link diagram $D_1$ and $D_2$ represent the same framed link if and only if $D_1$ can be transformed to $D_2$ by a sequence of Reidemeister moves
$\mathcal{R} I^{'} , \mathcal{R} II $ and $\mathcal{R} III$.
\end{thm}

\begin{prop}
\label{prop:graphpresentation}
The following group $\mathcal{G}_{L, D}$ only depends on the link $L$, and $\mathcal{G}_{L, D} \cong \mathcal{G}_{L, b } $ where $b$ is any braid such that $L=\hat{b}$.

Generators of $\mathcal{G}_{L,D} $: $\{X^{j} _i \} _{1\leq i\leq m; 1\leq j\leq N_i } \sqcup \{ Y_k \} _{1\leq k \leq m }  $. Notice: we can understand as, the generator $X^{j} _i $ is assigned to the short arc $C^{j} _i$, and the generator $Y_k $ is assigned to the link component $L_k $ as something.  These assignments can be viewed as a kind of coloring.

Relations: $\{  R^{' }_k \} _{1\leq k\leq M } \cup \{ R^{''} _k \} _{1\leq k\leq M } $. Notice: every crossing point $P_k $ is associated with a pair of relations $R^{'} _k $ and $R^{''} _{k} $defined as follows.

If $P_k $ is a crossing like (1) of Figure 11, and $C^{j_1 } _j , C^{j_2 } _j , C^{i_1 } _i , C^{i_2 } _i $ are the segments surrounding $P_k$ as in Figure 11. Then  $R^{'} _k $ is: $X^{i_2 } _i = (Y_i  )^{-1} X^{i_1 } _{i} $,  and $R^{''} _{k}$ is: $X^{j_2 } _{j} = X^{j_1} _j (X^{i_1} _i )^{-1} Y_i X^{i_1 } _i $.

If $P_k $ is a crossing like (2) of Figure 11, and $C^{j_1 } _j , C^{j_2 } _j , C^{i_1 } _i , C^{i_2 } _i $ are the segments surrounding $P_k$ as in Figure 11. Then
$R^{'} _k $ is:  $ X^{j_2 }_j = Y_j X^{j_1 } _{j} $, and $R^{''} _k $ is:
$X^{i_2 } _i = X^{i_1 } _{i} ( X^{j_1 } _j  )^{-1} Y^{-1} _{j} X^{j_1 } _{j} $.

So we can denote the group $\mathcal{G}_{L,D}$ simply as $\mathcal{G}_{L}$.

 \begin{figure}[htbp]

  \centering
    \includegraphics[height=4.5cm]{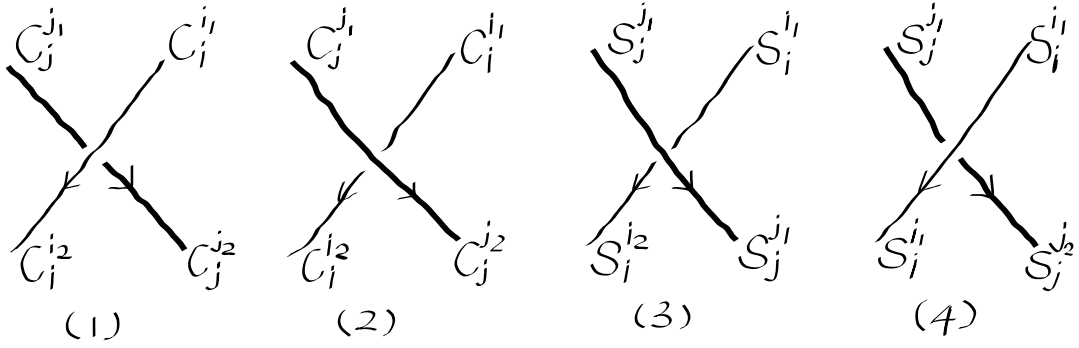}
  \caption{Crossing types}
\end{figure}

\end{prop}

\begin{rem}
The relations $R^{'} _{k} $ and $R^{''} _k$ are evidently derived from the maps in (2),(3) of Theorem \ref{thm:colored link invariant}.
 \end{rem}

\begin{pf}
We only need to prove that if $D^{'} $ is a link graph obtained by performing one of three types of the Reidmeister move on $D$, we have
$\mathcal{G}_{L, D^{'}} \cong \mathcal{G} _{L, D} .$ We can prove this through case by case check without any difficulty. Suppose $b$ is a braid such that
$L=\hat{b}$. Suppose $\Gamma $ is a "braid diagram " representing $b$.  Then the closure $\hat{ \Gamma }$ of the diagram $\Gamma$ is a link diagram representing $L$. We see evidently the group $\mathcal{G} _{L, \hat{ \Gamma }}$ is isomorphic to the group $\mathcal{G}_{L, b }$.
\end{pf}\\

We define a group invariant for oriented framed links similarly as follows. Still suppose $L= L_1 \sqcup L_2 \sqcup \cdots \sqcup L_m $ be a oriented framed $m-$components link. Suppose $D= D_1 \cup D_2 \cup \cdots \cup D_m $ is a link diagram representing $L$. Suppose $D_i $ is the part of $D$ representing the component $L$. Every $D_i$ is divided by the crossing points of $D$ into a set of curves which we call segments. So each segment of $D_i$ is a path on $D_i$ from one crossing point to the next crossing point along the orientation. Denote the set of segments on $D_i$  as $\{ D^{1} _i , D^{2} _i ,..., D^{N_i } _i \} .$ Denote the set of crossing points of $D$ as $\{  P_1 ,P_2 ,..., P_{M}   \}$.

\begin{prop}
\label{prop:framedgroupinvariant} The following group $\mathbb{G}_{L, D}$ only depends on the framed  link $L$.

Generators of $\mathbb{G}_{L,D} $: $\{X^{j} _i \} _{1\leq i\leq m; 1\leq j\leq N_i } \sqcup \{ Y_k \} _{1\leq k \leq m }  $. Notice: we can understand as, the generator $X^{j} _i $ is associated to the segment $D^{j} _i$, and the generator $Y_k $ is associated to the link component $L_k $ as something corresponding to the colour. We say the colour of the segment $X^{j} _i $ is $Y_i $.

Relations: $\{  R^{' }_k \} _{1\leq k\leq M } \cup \{ R^{''} _k \} _{1\leq k\leq M } $. Notice: every crossing point $P_k $ is associated with a pair of relations $R^{'} _k $ and $R^{''} _{k} $defined as follows.

If $P_k $ is a crossing like (1) of Figure 11, and $C^{j_1 } _j , C^{j_2 } _j , C^{i_1 } _i , C^{i_2 } _i $ are the segments surrounding $P_k$ as in Figure 11. Then

$R^{'} _k $ is: $X^{i_2 } _i =  X^{i_1 } _{i} $  ,  and $R^{''} _{k}$ is: $X^{j_2 } _{j} = X^{j_1} _j (X^{i_1} _i )^{-1} Y_i X^{i_1 } _i $.

If $P_k $ is a crossing like (2) of Figure 11, and $C^{j_1 } _j , C^{j_2 } _j , C^{i_1 } _i , C^{i_2 } _i $ are the segments surrounding $P_k$ as in Figure 11. Then  $R^{'} _k $ is:  $ X^{j_2 }_j = X^{j_1 } _{j} $ , and $R^{''} _k $ is:
$X^{i_2 } _i = X^{i_1 } _{i} ( X^{j_1 } _j  )^{-1} Y^{-1} _{j} X^{j_1 } _{j} $.

So we can denote the group $\mathbb{G}_{L, D}$ simply as $\mathbb{G}_{L}$.

\end{prop}

\begin{pf}
By Theorem \ref{thm:framedreidemeister}, we only need to show that if a link diagram $D^{'}$ is obtained from $D$ through performing the move $\mathscr{R}I$, $\mathcal{R}II$ or $\mathcal{R}III $ once, then $\mathbb{G}_{L, D^{'}} \cong \mathbb{G}_{L, D}$.
\end{pf}

Recall the group $\mathbb{G}_{L,b}$ introduced in Definition \ref{defi:groupinvariantbraidversion}, by a comparison of the presentation of the groups $\mathbb{G}_{L,D} $ with the group $\mathbb{G}_{L,b}$ we have the following lemma.
\begin{lem}
Suppose $L$ is a oriented framed link. Let $\Gamma$ is a $n$ string braid diagram such that the closure $\hat{\Gamma }$ represents $L$. Denote the $n$ string braid represented by $\Gamma$ as $b$. Then we have $\mathbb{G}_{L, b} \cong \mathbb{G}_{L, \hat{\Gamma }} $.
\end{lem}

Because of the relations $X^{j_2} _j = X^{j_1} _j $ and $X^{i_2} _i =X^{i_1} _i $ in Proposition \ref{prop:framedgroupinvariant}, the presentation of the group
$\mathbb{G}_{L, D}$ can be simplified as follows. On a link diagram $D=D_1 \sqcup \cdots \sqcup D_m $, we call a connected component as a "long arc". So a long arc is composed by some short arcs. For every crossing point on a long arc, the long arc pass over it except for the first one (the head ) and the last one (the tail). This conception is significant in the Wirtinger presentation of Link groups. Every component $D_i$ of $D$ is decomposed into a  union of long arcs: $D_i = S^{1}_i \cup \cdots \cup S^{M_i } _i$. Still denote the set of crossing points on $D$ are $\{ P_1 , P_2 ,..., P_M   \}$.

\begin{prop}
\label{prop:simplifiedpresentation}
The group defined as follows are isomorphic to the group $\mathbb{G}_{L, D}$.

Generators: $\{ X^{j} _{i} \} _{1\leq i\leq m , 1\leq j\leq M_i   }   \cup \{  Y_k \} _{1\leq k\leq m}$.

Relations: $ R_1 , ..., R_{M} $, where $R_k $ is a relation (associated to the crossing point $P_k $) as follows.

If $P_k $ is a crossing like (3) of Figure 11, then $R_k$ is: $ X^{i_2 } _{i} =X^{i_1} _{i} (X^{j_1 } _j )^{-1} Y^{-1} _j X^{j_1} _j  $,

if $P_k$ is a crossing like (4) of Figure 11, then $R_k$ is: $X^{j_2 } _{j}=X^{j_1} _{j} (X^{i_1} _i )^{-1} Y_i X^{i_1 } _{i} $.

\end{prop}

\begin{thm}
\label{thm:orientationindependent}
As a invariant of oriented links( framed links ), $\mathcal{G}_{L}$ ( $\mathbb{G}_{L} $ ) does not depend on the orientation of the components of $L$. That is, if $L^{'}$ is the oriented link (framed link ) obtained from $L$ by reversing the orientation of a component $L_i $ for any $1\leq i\leq m$, then
$\mathcal{G}_{L^{'} } \cong \mathcal{G}_{L}  $ (  $\mathbb{G}_{L^{'}} \cong \mathbb{G}_{L}  $   ).
\end{thm}

\begin{pf}
Suppose $L$ is  a oriented link diagram. Denote the link represented by $L$  still  as $L$ without making any confusion. Denote the oriented link diagram obtained by reversing the part of $L$ representing the component $L_i$ as $L^{'}$. So $L^{'}$ is a link diagram representing $L^{'}$.  As before, suppose the set of long arcs on the component $L_j $ is $\{ S^{1} _i , ..., S^{M_j } _{i} \} $ for $1\leq j \leq m$. Suppose the set of crossing points on $L$ is $\{ P_1 ,..., P_M \} . $  It is easy to see the set of long arcs on $L^{'}$ are in natural one to one correspondence with those on $L$. The difference is that the orientation of the arcs on $L_i $ are reversed. Denote the long arc on $L^{'}$ corresponding to $S^{d} _j$ still as $S^{d} _j $. Similarly, denote the crossing point on $L^{'}$ corresponding to $P_j $ ($1\leq j\leq M $) still as $P_j $. We compare the presentations in Proposition \ref{prop:simplifiedpresentation}  for $\mathbb{G}_{L}$ and $\mathbb{G}_{L^{'}}$.  Suppose the set of generators of $\mathbb{G}_{L}$ as in Proposition \ref{prop:simplifiedpresentation} is $\{ X^{d}  _j \} _{1\leq j\leq m, 1\leq d\leq M_j } \cup \{  Y_k \} _{1\leq k\leq M } $. Suppose the set of generators of $\mathbb{G}_{L^{'}}$ as in Proposition \ref{prop:simplifiedpresentation} is $\{ \bar{X}^{d}  _j \} _{1\leq j\leq m, 1\leq d\leq M_j } \cup \{ \bar{ Y}_k \} _{1\leq k\leq M } $. We claim that the map $\phi( X^{d} _j )= \bar{X}^{d} _j $ ($1\leq j\leq m, 1\leq d\leq M_j $  ); $\phi (Y_k ) = Y_k  $ ($k\neq i $); $\phi (Y_i ) = (Y_i )^{-1} $ extends to an isomorphism $\phi : \mathbb{G}_{L} \rightarrow \mathbb{G}_{L^{'}} $. We only need to consider those six different cases in Figure 12. Where on every column of Figure 12, the graph above means a crossing point in $L$, the graph below means the corresponding crossing point in $L^{'}$. We omit the detail since it is easy.

 \end{pf}

 \begin{figure}[htbp]

  \centering
    \includegraphics[height=5.4cm]{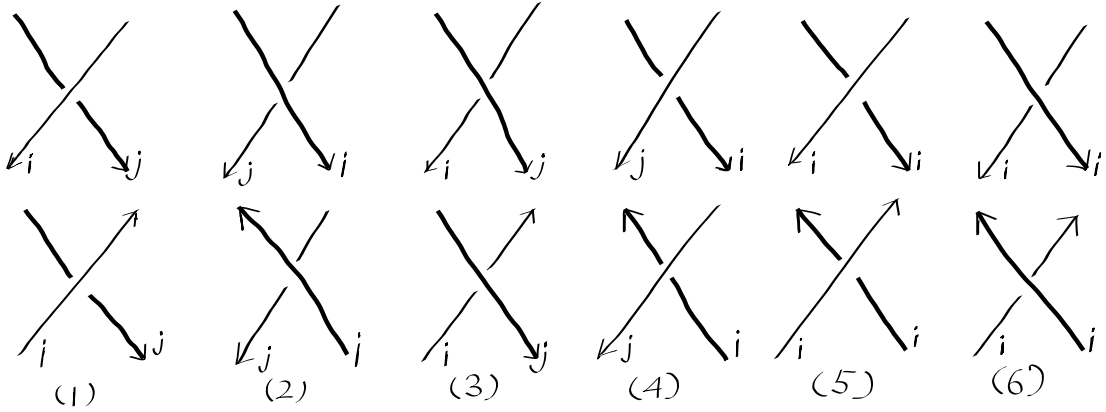}
  \caption{Different crossings }
\end{figure}

\section{A link group invariant}

In Section 6 we construct a link group invariant dominating the link invariant $\Lambda_{G;[(g,b)]}$ (Theorem \ref{thm:extendableGRmatrix}, Remark \ref{rem:extendableinvariant} ) by using extended pairs, and their colored version $F$ (Theorem \ref{thm:extended colored link invariant} ).  The construction is similar with the group invariant $\mathcal{G}_{L}$ in last section but a little more complicated than that.

Still let $L= L_1 \sqcup L_2 \sqcup \cdots \sqcup L_m  $ be a oriented link with $m$ components.  Let $b \in B_{n}$ be a braid whose closure $\hat{b}$ presents $L$. Let $G$ be a finite group, suppose $\bar{L}$ is a $\mathscr{E}G-$colored link based on $L$ by assign a color $(g_i ,b_i ) \in \mathscr{E} G$ to the component $L_i$ of $L$, for $1\leq i\leq m$. Suppose $\bar{b}$ is a $\mathscr{E}G-$colored braid based on $b$, whose closure is $\bar{L}$. Suppose
$\partial ^{-} (\bar{b})= \partial ^{-} (\bar{b}) = ( (g_{j_1} b_{j_1} ), ..., (g_{j_{n}} ,b_{j_{n}})  )$. Recall the functor $F(-)$ constructed in Theorem \ref{thm:extended colored link invariant}.

\begin{lem}
$F ( \bar{L} ) = tr ( ( (b_{j_1 })_{L} \otimes ...(b_{j_n })_{L}   )  F ( \bar{b})   ) .$
\end{lem}

Now we associate a group to $b $. We introduce variables
$$X_1 , X_2 ,..., X_n , C_1 ,B_1 , C_2 ,B_2,..., C_m ,B_m ,$$
 and let $ F_{n,m}$ be the group generated by these elements along with the following relations:
$$(1)\  C_i B_i = B_i C_i  , \  for\  1\leq i\leq m  ;\ \ \ \ \ \ (2)\  B_i ^{2} =e \  (the\  unit).$$
Then we color the
$i-$th string of $b$ with the extended pair $( C_{j_i } , B_{j_i } )$, if it (as a subset of $L= \hat{b}$) belongs to the component $L_{j_i}$. Denote the resulted $\mathscr{E}F_{n,m}-$colored braid as $b^{''}$. Recall the  map $\Psi_{g_1 , g_2}$ and $\Psi^{-} _{g_1 , g_2} $
before Theorem \ref{thm:colored link invariant} and (2),(3) of Theorem \ref{thm:extended colored link invariant}. For a $\mathscr{E}F_{n,m}-$colored $n-$string simple braid
$X^{+} _{i; (C_{j^{'} }, B_{j^{'}} ) ,( C_{j^{''} } ,B_{j^{''}} )}$ or $X^{-} _{i; (C_{j^{'} }, B_{j^{'}} ) ,( C_{j^{''} } ,B_{j^{''}} ) }  $, whose $i-$th string and $(i+1)-$th string are colored by $(C_{j^{'}} ,B_{j^{'}})$ and $(C_{j^{''}} ,B_{j^{''}})$ respectively, and $(Z_1 , Z_2 ,..., Z_n )\in \underbrace{ F_{n,m} \times F_{n,m} \times ... \times F_{n,m} }_{n}$, we set
\begin{align*}
&\mathcal{F}^{'} ( X^{+} _{ i; (C_{j^{'} } ,B_{j^{'}}) ,( C_{j^{''}}  ,B_{j^{''}})}  )(\cdots , Z_i , Z_{i+1} , \cdots    )=  (\cdots , (B_{j^{''}} C_{j^{''}})Z_{i+1} Z_i ^{-1} C_{j^{'} }^{-1} Z_i ,  Z_i ,\cdots   ) , \\
&\mathcal{F} ^{'} ( X^{-} _{i; (C_{j^{'} } ,B_{j^{'}}) ,( C_{j^{''}}  ,B_{j^{''}})} )(\cdots ,Z_i ,Z_{i+1} , \cdots )= (\cdots , Z_{i+1} , (B_{j^{'}} C_{j^{'}} ^{-1}   )Z_i Z_{i+1} ^{-1} C_{j^{''} }Z_{i+1} ,\cdots )  .
\end{align*}

If $b^{''} = \Pi _{t=1} ^{N} X^{\epsilon_{t}} _{i_t ; ( C_{j^{'}_{t}} , B_{j^{'} _{t}}) , (C_{ j^{''}_{t} } ,B_{j^{''} _{t}} )  }$, we set
$$\mathcal{F}^{'} (b^{''} )= \mathcal{F}^{'} ( X^{\epsilon_{N}} _{i_N ; (C_{j^{'} _{N}} , B_{j^{'} _{N} } ) , (C_{ j^{''} _{N}  },B_{ j^{''} _{N}   }    ) }\circ \cdots \circ \mathcal{F } ^{'} (  X^{\epsilon_{1}}   _{i_1 ; (C_{j^{'} _{1}} , B_{j^{'} _{1} } ) , (C_{ j^{''} _{1}  },B_{ j^{''} _{1}   }    )  } ) .$$
 We define set of words $\bar{W}^{b} _1 ,...,\bar{W}^{b} _n \in F_{n,m }$  by
$$\mathcal{F}^{'} (b^{''} ) ( (X_1 ,X_2 ,..., X_n   )  )= (\bar{W}^{b} _1 , \bar{W}^{b} _2 ,..., \bar{W}^{b} _n   ).$$

 We can view $\mathcal{F}^{'}(b^{''})$ as a group isomorphisms from $F_{n,m}$ to itself, by sending generators $(X_1 ,..., X_n , C_1 ,B_1 , ..., C_m ,B_m )$ to $( \bar{W}^{b} _1 ,..., \bar{W}^{b} _n   , C_1 ,B_1 , ..., C_m ,B_m ) $ .
 To show the reliance of $\bar{W}^{b} _i $  on the generators and the order of colors of the strings of $b$, we denote them as $\bar{W}^{b} _i ( X_1 , ..., X_n ;C_{j_1} ,B_{j_1}, ..., C_{j_n} ,B_{j_n}  ) $ .

 If a group morphism $\phi $ map $F_{n,m} $ to some group $G$ such that $\phi (X_i ) = x_i \in G $ , $\phi (C_j ) = g_j \in G$ and $ \phi (B_j) = b_j $ for $1\leq j \leq n $, then naturally we write the image
$\phi (\bar{W}^{b} _i )$ as $\bar{W}^{b} _i (x_1 ,..., x_n ;g_{j_1}, b_{j_1} ,..., g_{j_n} , b_{j_n}  ) \in G  $ for $1\leq i \leq n$. It is easy to see
$$\mathcal{F}(\bar{b}) (( x_1 , ..., x_n ) )= (\bar{W}^{b} _1 (x_1 ,..., x_n ; g_{j_1}, b_{j_1} ,...,g_{j_n}, b_{j_n} ),..., \bar{W}^{b} _n (x_1 ,..., x_n ;g_{j_1}, b_{j_1} ,...,g_{j_n}, b_{j_n}  )   ) . $$
This shows the meaning of these elements in $F_{n,m}$. Return to the original $\mathscr{E}G-$colored braid $\bar{b}$.   Since $\hat{\bar{b}}=\bar{L}$, every string of $\bar{b}$ is colored by some element in $\{ (g_1,b_1) ,..., (g_m,b_m) \} $. Suppose the $t-$th string of $b$ is colored by $ (g_{j_t },b_{j_t})$.

\begin{defi}
\label{defi:groupinvariantbraidversion2}
Assuming above notations,  we define a group $\bar{\mathcal{G}}_{L, b}$ (as quotient groups of $F_{n,m}$) generated by letters $X_1 ,X_2 ,..., X_n , C_1 ,B_1,C_2 ,B_2,..., C_m ,B_m$,  with the following set of relations.
(1) $C_i B_i = B_i C_i $  for  $ 1\leq i\leq m $ ;\\
(2) $ B_i ^{2} =e $ ;\\
(3) $ X_t = B_{j_t} \bar{W}^{b} _t (X_1 ,..., X_n ;C_{j_1} , B_{j_1},..., C_{j_n}, B_{j_n}  )$ for  $  1\leq t \leq n $ .

\end{defi}

\begin{prop}
\label{prop:groupinandnumberin}
We have:

(1) $\mathcal{F} (\bar{L}) = tr(((b_{j_1})_{L}\otimes (b_{j_2})_{L}\otimes ...\otimes (b_{j_n})_{L}  ) \mathcal{F}( \bar{b}) )$

$= |\{  f\in Hom (\bar{\mathcal{G}}_{L,b} , G   )   | f(C_1 )= g_1 , \ f (B_1 )=b_1,\ ..., \ f(C_m) = g_m, \ f(B_m)=b_m \} |.$

\end{prop}
\begin{pf}
Since the linear map $((b_{j_1})_{L} \otimes ...\otimes (b_{j_n})_{L}   )\mathcal{F}(\bar{b})$ comes from the bijective map( permutation) $((b_{j_1})_{L} \otimes ...\otimes (b_{j_n})_{L}   )\mathcal{F}(\bar{b}) : \underbrace{ G \times G \times \cdots \times G  }_{n} \rightarrow \underbrace{G \times G \times \cdots \times G  }_{n}  $, so $tr (((b_{j_1})_{L} \otimes ...\otimes (b_{j_n})_{L}   )\mathcal{F}(\bar{b}) )$ equals the number of fixed points
of the permutation $((b_{j_1})_{L} \otimes ...\otimes (b_{j_n})_{L}   )\mathcal{F}(\bar{b})$. But now there is a bijection from the set of fixed points of $((b_{j_1})_{L} \otimes ...\otimes (b_{j_n})_{L}   )\mathcal{F}(\bar{b})$ to the set of group morphisms
$\{  f\in Hom (\bar{\mathcal{G}}_{L,b} , G   )   | f(C_1 )= g_1 , \ ,f(B_1)=b_1,\  f(C_2 )= g_2 ,\ f(B_2)=b_2,\ ..., \ f(C_m) = g_m,\ f(B_m)=b_m \}$, so we finish the proof.
\end{pf}

So all our link invariants $\mathcal{F}(\bar{L} )$ constructed in Theorem \ref{thm:extended colored link invariant} from various finite groups $G$  and the sets $\mathscr{E}G$ are dominated by the data
$( \bar{ \mathcal{G}}_{L,b } ; (C_1 ,B_1 ), ...,( C_m ,B_m )  ) $, that is , a group with a set of prescribed ordered set of extended pairs . As before in the next theorem we show this data is a link invariant itself.

\begin{thm}
\label{thm:groupinvariant2}
Suppose $d\in B_{N} $ is another braid such that $\hat{d} = L$. Suppose the set of generators of the group $\mathcal{G}_{L,d} $ in above construction is
$\{ X_{1} ^{'} ,..., X_{N} ^{'} , Z_1 , ..., Z_m \} $. Then there exists a group isomorphism $ \Phi : \bar{\mathcal{G}}_{L, b} \rightarrow \bar{\mathcal{G}}_{L,d} $ such that
$\Phi (Y_i ) = Z_{ \epsilon (i ) } $ for $1\leq i\leq m$ and some permutation $\epsilon \in S_{m}$.
\end{thm}
\begin{pf}
Since $\hat{d}=\hat{b}$, by Markov's theorem,  we only need to consider the following two cases: Case 1. The braid $b$ is transformed to $d$ by performing  Markov move I once; Case 2. The braid $b$ is transformed to $d$ by performing Markov move II once.

Case 1. Suppose $b^{'}$ is the $F_{n,m}-$colored braid obtained from $b$ by color the $t-$th string of $b$ with $(C_{j_t} ,B_{j_{t}})$ as before.
Recall the notations in the beginning of section 7.   Suppose $\partial ^{+}(b^{'} ) = \partial ^{-} (b^{'})=  ( ((C_{j_1 },B_{j_1}) ,+ ), ..., ( (C_{j_n },B_{j_n}) ,+  )   ) . $
Since $b^{'}$ is a braid whose strings are oriented from upward to downward, the signals in the brackets are all $+$. Since $b^{'}$ can be closed to give a
$F_{n,m}- $coloured link, we must have $\partial ^{+}(b^{'} ) = \partial ^{-} (b^{'}) $. Now suppose $\alpha $ is a $n-$string $\mathscr{E}F_{n,m}$ coloured braid such that $\partial^{+} ( \alpha  )= (((C_{j_1 } ,B_{j_1}) ,+ ), ..., ( (C_{j_n } ,B_{j_n} ),+  )    ) $.  Denote the natural epimorphism $B_n \rightarrow S_n$ from the braid group to the permutation group as $\pi$, and suppose $\pi (\alpha )=s \in S_n $. Denote the mirror image of $\alpha $, which is also a $F_{n,m}-$colored braid, as $\alpha^{-1}$. Then we have
\begin{align*}
&\partial^{+} (\alpha^{-1} )= \partial ^{-} (\alpha )= (((C_{j_{s(1)} } ,B_{j_{s(1)}}) ,+ ), ..., ( (C_{j_{s(n)} } ,B_{j_{s(n)}} ),+  )    ),\  and\\
&\partial^{-} (\alpha^{-1} )=(((C_{j_1 } ,B_{j_1}) ,+ ), ..., ( (C_{j_n } ,B_{j_n} ),+  )    ).
\end{align*}
Denote $d^{'} = \alpha \cdot b^{'} \cdot \alpha^{-1} $ as a $\mathscr{E}F_{n,m}-$colored braid, then $d^{'}$ is a braid obtained from $b^{'}$ by performing Markov move I once. Any braid  $d^{'}$ as that can be written in this way for suitable braid $\alpha$. What we need to show is $\bar{\mathcal{G}}_{L,b^{'}}\cong \bar{\mathcal{G}}_{L, d^{'}}. $

We denote a term like $\bar{W}^{\alpha} _i (\bar{W}^{\beta} _1 ,...,\bar{W}^{\beta} _n ; C_{j_1} ,B_{j_1},..., C_{j_n },B_{j_n}   ) $ simply as
$$ \bar{W}^{\alpha} _i (\bar{W}^{\beta}  ; C_{j_1} ,B_{j_1},..., C_{j_n },B_{j_n}   )  $$ for convenience. Since $\alpha \cdot \alpha^{-1} =id_{B_n } $ and $\alpha^{-1} \cdot \alpha = id_{B_n }  $, we have
\begin{align*}
&(A) \bar{W}^{\alpha} _i (\bar{W}^{\alpha^{-1}} (X_1 , ..., X_n ;C_{j_{s(1)}},B_{j_{s(1)}},...,C_{j_{s(n)}}, B_{j_{s(n)}}  ); C_{j_1}, B_{j_1} ,..., C_{j_n },B_{j_n}  ) = X_i \  for\ 1\leq i\leq n.\\
&(B) \bar{W}^{\alpha^{-1}} _i (\bar{W}^{\alpha} (X_1 ,..., X_n ; C_{j_1} ,B_{j_1},..., C_{j_n }, B_{j_n}   ) ; C_{j_{s(1)}}, B_{j_{s(1)}} ,..., C_{j_{s(n)} },B_{j_{s(n)}}  )= X_i \  for\  1\leq i\leq n .
\end{align*}
Now by definition, the group $\mathcal{G}_{L, d^{'}}$ has the following presentation:
\begin{align*}
&<X_1 ,...,X_n ,C_1 ,B_1,..., C_m ,B_m | X_i =\\
&B_{j_{s(i)}} \bar{W}^{\alpha \cdot b^{'} \cdot \alpha^{-1}} _i (X_1 ,...,X_n ;C_{j_{s(1)}},B_{j_{s(1)}} ,...,C_{j_{s(n)}},B_{j_{s(n)}} ) , (1\leq i\leq n) >
\end{align*}
Here we add new generators $V_i$ (corresponding to $\bar{W}^{\alpha^{-1}} _i $ ), by Tietze theorem, above presentation is equivalent to the following
\begin{align*}
&<X_i , V_i , C_j ,B_j |  V_i = \bar{W}^{\alpha^{-1}} _i (X_1 ,..., X_n ;C_{j_{s(1)}}, B_{j_{s(1)}},...,C_{j_{s(n)}}, B_{j_{s(n)}} ),\\
&X_i =B_{j_{s(i)} } \bar{W}^{\alpha \cdot b^{'} }(V_1 ,..., V_n ; C_{j_1} ,B_{j_1},...,C_{j_n } ,B_{j_n}) > ,
\end{align*}
where $1\leq i\leq n$ and $1\leq j\leq m $ in above presentation for all $i,j$.  Now by using (1) and (2) of Definition \ref{defi:groupinvariantbraidversion} and by  definition of the elements $\bar{W}^{b} $, we have
\begin{align*}
& B_{j_{s(i)}}   \bar{W}^{\alpha \cdot b^{'} } _{i}(V_1 ,..., V_n ; C_{j_1} ,B_{j_1},...,C_{j_n } ,B_{j_n}) = \\
&\bar{W} ^{\alpha}_{i} ( B_{j_1}  \bar{W} ^{b^{'}}_1 ( V_1 ,..., V_n ; C_{j_1} ,B_{j_1},...,C_{j_n } ,B_{j_n} ),...,B_{j_n}  \bar{W} ^{b^{'}}_n ( V_1 ,..., V_n ; C_{j_1} ,B_{j_1},...,C_{j_n } ,B_{j_n} ) ; \\
&C_{j_{s(1)}} ,B_{j_{s(1)}},...,C_{j_{s(n)} } ,B_{j_{s(n)}}   ) .
\end{align*}

 Then by above identity (A) and (B), we see the set of relations

$ \{ V_i = \bar{W}^{\alpha^{-1}} _i (X_1 ,..., X_n ;C_{j_{s(1)}}, B_{j_{s(1)}},...,C_{j_{s(n)}}, B_{j_{s(n)}} ) \} _{(1\leq i\leq n)} $  can be replaced by  the set of relations
$\{ X_i = \bar{W}^{\alpha} _i (V_1 ,..., V_n ; C_{j_{1}}, B_{j_{1}},...,C_{j_{n}}, B_{j_{n}} )\} _{ (1\leq i\leq n )} $. So above presentation is equivalent to the following
\begin{align*}
&<X_i , V_i , C_j ,B_j | X_i = \bar{W}^{\alpha} _i (V_1 ,..., V_n ; C_{j_{1}}, B_{j_{1}},...,C_{j_{n}}, B_{j_{n}} ),\\
&X_i =\bar{W} ^{\alpha} _{i} ( B_{j_1}  \bar{W} ^{b^{'}}_1 ( V_1 ,..., V_n ; C_{j_1} ,B_{j_1},...,C_{j_n } ,B_{j_n} ),...,B_{j_n}  \bar{W} ^{b^{'}}_n ( V_1 ,..., V_n ; C_{j_1} ,B_{j_1},...,C_{j_n } ,B_{j_n} ) > ,
\end{align*}
where $1\leq i\leq n$ and $1\leq j\leq m $.

By using Tietze theorem again, above presentation is equivalent to the following
\begin{align*}
&< V_i , C_j ,B_j | X_i = \bar{W}^{\alpha} _i (V_1 ,..., V_n ; C_{j_{1}}, B_{j_{1}},...,C_{j_{n}}, B_{j_{n}} )=\\
&\bar{W} ^{\alpha} _{i} ( B_{j_1}  \bar{W} ^{b^{'}}_1 ( V_1 ,..., V_n ; C_{j_1} ,B_{j_1},...,C_{j_n } ,B_{j_n} ),...,B_{j_n}  \bar{W} ^{b^{'}}_n ( V_1 ,..., V_n ; C_{j_1} ,B_{j_1},...,C_{j_n } ,B_{j_n} ) > ,
\end{align*}
where $1\leq i\leq n$ and $1\leq j\leq m $.

By using above identity (A) and (B) again, we see above presentation is equivalent to the following
$$<  V_i , C_j ,B_j | V_i = \bar{W}^{b^{'}} _i ( V_1 ,..., V_n ; C_{j_1} ,B_{j_1},...,C_{j_n } ,B_{j_n} ) )  (1\leq i\leq n )> ,$$
which is exactly a presentation of the group $\mathcal{G} _{L, b}$.

Case 2. Suppose the $n-$string braid $b , b^{'}$ are as in the case 1. Consider the braid $\sigma_{n} b^{'} $ , for which we color the $(n+1)$-th string also
by $Y_{i_n }$. Then the closure of $\sigma_{n} b^{'} $ is also $L$ and it is a braid obtained from $b^{'}$ by performing the  Markov move II once. Since the closure of $\sigma_{n} b^{'} $ is also $L$, we see the colour of the $(n+1)-$th string of $\sigma _{n} b^{'}$ should also be $( C_{j_n}  , B_{j_n } )$.

 We need to show $\mathcal{G}_{L , \sigma_n b^{'}} \cong \mathcal{G}_{L, b^{' }} $. By definition, the map $\mathcal{F}(\sigma_n b^{'})$ maps

$(X_1 , ..., X_{n-1} , X_n , X_{n+1}) $ to $(\bar{W}^{b^{'}} _1 ,..., \bar{W}^{b^{'} } _{n-1} , B_{j_n} C_{j_n} X_{n+1} (\bar{W}^{b^{'} } _n  )^{-1} C_{j_n } ^{-1} \bar{W}^{b^{'}} _{n} ,  \bar{W}^{b^{'}} _n    )   $, so the group $\mathcal{G}_{L, \sigma_n b^{'}}$ has the following presentation
$$<X_i ,C_j ,B_j | X_k = B_{j_k} \bar{W}^{b^{'}} _{k} , X_n = B_{j_n }  B_{j_n} C_{j_n} X_{n+1} (\bar{W}^{b^{'}} _n )^{-1} C_{j_n }^{-1}\bar{W}^{b^{'}} _n , X_{n+1} =B_{j_n}W^{b^{'}} _n     >,$$
where $1\leq i\leq n+1 , 1\leq j\leq m , 1\leq k\leq n-1$. It is easy to see by canceling the generator $X_{n+1}$ using the last relation, above presentation is equivalent to the following presentation
$<X_i, C_j ,B_j  | X_i = B_{j_i } \bar{W}^{b^{'}} _i ( 1\leq i\leq n)  > $, which is just a presentation of the group $\mathcal{G}_{L, b}$.

\end{pf}$\\$

\section{A further generalization of the extended $R-$matrix }

It turns out the conditions of the extended $R-$matrix can be relaxed further, and still give link invariants. We give a sketch of this construction in this section, whose proofs are similar.

\begin{defi}
\label{defi:generalizedextended}
A generalization of the extended $R-$matrix is a combination $(I,f,c_1 ,c_2 ,d )$, such that $I\in End(V\otimes V)$, $f, c_1 , c_2 ,d \in End (V) $ and:
\begin{enumerate}
\item[(1)] $I$ is a $R-$matrix;

\item[(2)] $tr_2 (I) =c_1 $, $tr_2 (I^{-1})=c_2 $;

\item[(3)] $c_1 c_2 =d^2 $;

\item[(4)] $c_1 , c_2 $ and $d$ are central with respect to $I$;

\item[(5)] $(f\otimes f) I = I(f\otimes f)$;

\item[(6)] $c_1 ,c_2 $ and $d$ commute with $I$.
\end{enumerate}
\end{defi}

A special case of Definition \ref{defi:generalizedextended} is as follows.

\begin{defi}
\label{defi:specialgeneralizedextended}
A special generalized extended $R-$matrix is a combination $( I, f, d ) $ such that  $r\in End(V\otimes V)$, $f,d \in End (V) $ and:
\begin{enumerate}
\item[(1)] $I$ is a $R-$matrix;

\item[(2)] $tr_2 (I) =tr_2 (I^{-1})=d $;

\item[(3)] $d$ is central with respect to $I$;

\item[(4)] $(f\otimes f) I=I(f\otimes f)$;

\item[(5)] $d f= fd$.
\end{enumerate}

\end{defi}

It is evident if $c_1 =c_2 =d$ in Definition \ref{defi:generalizedextended}, then that combination is a special generalized extended $R-$matrix.  Nevertheless, we have the following lemma.

\begin{lem}
Suppose $(I,f,  c_1 ,c_2 , d , )$ is a generalized extended $R-$matrix, if set $c= d^{-1} c_1 = d (c_2 )^{-1} $, and $\bar{I} = (id_{V} \otimes c^{-1} )I $, then, $(\bar{I} , f, d )$ is a special generalized extended $R-$matrix.
\end{lem}
\begin{pf}
Since $c$ is central with respect to $I$, then we can prove $\bar{I}$ is a $R-$matrix in the same way as Lemma \ref{lem:modifiedRmatrix}. And we have
$(\bar{I})^{-1} = (c\otimes id_{V})I^{-1}$. The main claim of the lemma follows from the following identities.
\begin{align*}
&tr_2 ( (id_{V} \otimes f)\bar{I} )=tr_2 ((id_{V} \otimes f) I (c^{-1} \otimes id_{V} ) )= tr_2 ((id_{V} \otimes f)I) c^{-1}
 = c_1 c^{-1} = d ;\\
 &tr_2 ( (id_{V} \otimes f)\bar{I}^{-1} )=tr_2 ( (c \otimes id_{V} ) (id_{V} \otimes f) I^{-1} )=c tr_2 ((id_{V} \otimes f)I^{-1})
 = c c_2 = d   .
 \end{align*}
\end{pf}

Now suppose $(I, f, d)$ is special generalized extended $R-$matrix as in Definition \ref{defi:specialgeneralizedextended}, we construct a link invariant.
First, for any $n$, denote the braid group representation induced by the $R-$matrix $I$ as $\rho_n :$ $B_n \rightarrow End (V^{\otimes n}) $. Denote the natural surjection from $B_n$ to the permutation group $S_n$ as $\pi $. Suppose $g\in B_n$. Choose any concrete braid $a$ representing $g$. Then the $\pi (g)^{-1} (i)-$th string connect the $\pi(g)^{-1} (i)- $th end point in $\partial ^{+} (g) $ to the $i-$th end point in $\partial ^{-} (g)$. For any $1\leq i\leq n$, we define the following integers $m^{+} _i (a) ,m^{-} _i (a), m_i (g) $.
\begin{itemize}
\item In the braid $a$, the $\pi (g)^{-1} (i)-$th string pass over another string as (1) of Figure 13 for $m^{+} _i (a) $ times;

\item In the braid $a$, the $\pi (g)^{-1} (i)-$th string pass over another string as (2) of Figure 13 for $m^{-} _i (a)$ times;

\item $m_i (g) = m^{+} _i (a) - m^{-} _i (a) $.
\end{itemize}

\begin{thm}
The number $m_i (g)$ is well defined, that is, does not depend on representatives of the braid $g\in B_n$.

\end{thm}
\begin{pf}
Suppose a concrete braid $a$ is a representative of the braid $g$. It is easy to see whenever another concrete braid $b$ is obtained from $a$ through a transformation (in some part of $a$) $\sigma_i  \sigma^{-1} _i  \longleftrightarrow 1 $, $ \sigma_i \sigma_j \longleftrightarrow \sigma_j \sigma_i $ for $|i-j|\geq 2$, or $\sigma_i \sigma_{i+1} \sigma_i \longleftrightarrow \sigma_{i+1} \sigma_i \sigma_{i+1}  $, we have $m_i (a) = m_i (b)$. So the theorem is proved.
\end{pf}

 \begin{figure}[htbp]

  \centering
    \includegraphics[height=4.2cm]{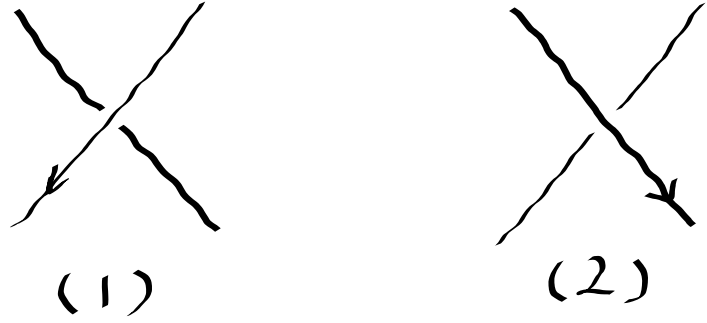}
  \caption{Overcrossing types }
\end{figure}

\begin{lem}
\label{lem:movingd}
Let $(I,f,d)$ be a generalized extended $R-$matrix, and let $\rho_n$ be the representation defined as above. Suppose $g\in B_n$. Let $\pi : B_n \rightarrow S_n$ be the natural surjective map. Denote the abelian subgroup of $S_n$ generated by the permutation $\pi(g)$ as $<\pi (g)>$.  The action of $< \pi(g)  >$ on the set $\{ 1,2,...,n \} $ decompose it into several equivalence classes ( corresponding to the "cycle decomposition of the permutation $\pi(g)$") : $\{  1,2,...,n\} = \sqcup_{l=1,...,M} I_l $. Let $d_1 , ..., d_n$ be a sequence of integers. Then we have:
\begin{enumerate}
\item[(1)] $tr ( (d^{k_1 } \otimes ...\otimes d^{k_l })f^{\otimes n} \rho_n (g )   ) = tr (d^{k_1} \otimes ...\otimes d^{k_i - 1} \otimes ...\otimes d^{k_{\pi(g)(i)}+ 1} \otimes ...\otimes d^{k_n }    ) f^{\otimes n} \rho_n (g) ).$
\item[(2)] For $1\leq l\leq M$, denote $d_{I_l } = \sum_{j\in I_l } d_j  $. Then the value
$tr ( (d^{k_1 } \otimes ...\otimes d^{k_l })f^{\otimes n} \rho_n (g )   ) $ only depends on the sequence $d_{I_1 } ,..., d_{I_M }$.
\end{enumerate}
\end{lem}
\begin{pf}
First, (1) follows from the following identities.
\begin{align*}
 &tr ( (d^{k_1 } \otimes ...\otimes d^{k_l })f^{\otimes n} \rho_n (g )   ) \\
 =&tr ((d^{k_1} \otimes ...\otimes d^{k_i - 1} \otimes ...\otimes d^{k_n }    ) f^{\otimes n} \rho_n (g) ( id^{i-1} _V \otimes d \otimes  id^{n-i} _V  )  \\
 =&tr ( (d^{k_1} \otimes ...\otimes d^{k_i - 1} \otimes ...\otimes d^{k_n }    )( id^{\pi(g)(i)-1} _V \otimes d \otimes  id^{n-\pi(g)(i)} _V  )   f^{\otimes n} \rho_n (g)   )\\
 =&tr (d^{k_1} \otimes ...\otimes d^{k_i - 1} \otimes ...\otimes d^{k_{\pi(g)(i)}+ 1} \otimes ...\otimes d^{k_n }    ) f^{\otimes n} \rho_n (g) ).
\end{align*}

By using (1) repeatedly we can prove (2).
\end{pf}

\begin{thm}
\label{thm:linkinvariantfromgeneralizedextended}
Let $(I, f, d)$ is special generalized extended $R-$matrix. Let $g\in B_n$ be a $n-$string braid . Let $\rho_n , \pi , m^{+} _i (g) , m^{-} _i (g) , m_i (g)$ be defined as in above setting.  Let
$$\Delta _{(I,f,d)  } (g) = tr (  f^{\otimes n} ( d^{m_1 (g)} \otimes d^{m_2 (g)} \otimes ...\otimes d^{m_n (g)}  )   \rho (g ) ) .$$
Then $\Delta _{(I,f,d)}$ is invariant under both two types of Markov moves, thus define a link invariant.
\end{thm}

\begin{pf}
Suppose $h\in B_n$, first we show $\Delta_{I,f,d} (h^{-1}gh)= \Delta_{I,f,d}(g)$.  Since $\pi (h) ( 1,2,...,n  )= ( \pi(h)(1),\pi(h)(2),..., \pi(h)(n)  )$, so we have
 $ \pi (h^{-1}) ( \pi(h)(1), \pi(h) (2), ..., \pi(h)(n)  ) =(1,2,...,n)$. It is easy to see the following identities.
 $$ m_{  i} (h^{-1}) = - m_{\pi(h) (i)} (h )   \  for\  1\leq i\leq n .$$

 And we can see
 \begin{align*}
   (*)\ \ \   m_{i} (h^{-1}gh ) = m_{i} (h^{-1}) + m_{\pi(h)(i)}(g) + m_{ \pi( g^{-1} h )(i) } (h)\\
 = m_{i} (h^{-1}) + m_{\pi(h)(i)}(g) - m_{ \pi( h^{-1} g^{-1} h )(i) } (h).
 \end{align*}

 So we have
 \begin{align*}
  \Delta_{I,f,d} (h^{-1}gh) &= tr( (f^{\otimes n}) ( d^ {m_{1} (h^{-1}gh)} \otimes ...\otimes d^{m_n (h^{-1}gh )}  )\rho(h^{-1}gh) )\\
  &= tr( (f^{\otimes n}) (d^{m_{\pi(h)(1)}(g)} \otimes ...\otimes d^{m_{\pi(h)(n)}(g) }) \rho(h^{-1}gh) )\\
   &= tr ((f^{\otimes n} ) \rho(h^{-1}) ( d^{m_1 }\otimes d^{m_2 }... \otimes d^{m_n }   ) \rho(g)\rho(h) )\\
   &=tr ( \rho(h^{-1})(f^{\otimes n} ) ( d^{m_1 }\otimes d^{m_2 }... \otimes d^{m_n }   ) \rho(g)\rho(h)  )\\
   &=tr((f^{\otimes n} ) ( d^{m_1 }\otimes d^{m_2 }... \otimes d^{m_n }   ) \rho(g)   )=\Delta_{I,f,d}( g).
 \end{align*}
Where the second identity sign is by using Lemma \ref{lem:movingd} and above identities $(*)$.   The third equality sign is by (3) of Definition \ref{defi:specialgeneralizedextended}, and the fourth equality sign is by (4) of Definition \ref{defi:specialgeneralizedextended}.  So we proved the value $\Delta_{I,f,d}(-)$ is invariant under the first type of Markov moves.

Now consider the braid $\sigma_{n} g \in B_{n+1} $ where $\sigma_{n}$ is the $n-$th generator of $B_{n+1}$. It is easy to see
$$m_{1} (\sigma_{n} g ) = m_1 (g) ,..., m_{n-1}(\sigma_{n}g) =m_{n-1}(g), m_{n} (\sigma_{n} g) =0, m_{n+1}(\sigma_{n} g) =m_{n} -1   .$$
So we have
\begin{align*}
\Delta_{I,f,d} &( \sigma_{n} g ) = tr (f^{\otimes (n+1)} (d^{m_1} \otimes ...\otimes d^{m_{n-1}} \otimes d^{0} \otimes d^{m_{n}-1} ) \rho(\sigma_{n} g) )\\
 &=tr ( tr_{n+1} (( f^{\otimes n}\otimes id_{V})( id^{\otimes n} _{V} \otimes f  ) \rho(\sigma_n) (d^{m_1 } \otimes d^{m_2 } \otimes ...\otimes d^{m_n -1}\otimes id_{V} )
 (\rho(g)\otimes id_{V}       )  )\\
 &=tr(  f^{\otimes n} (   id^{\otimes(n-1)} \otimes d ) (d^{m_1 } \otimes d^{m_2 } \otimes ...\otimes d^{m_n -1}  )\rho(g)   )\\
 &=\Delta_{I,f,d} ( g ).
\end{align*}
Where the second equality sign is by (3) of Definition \ref{defi:specialgeneralizedextended}, the third equality sign is by Lemma \ref{lem:partialtrace}. Similarly we can prove $\Delta_{I,f,d} ( \sigma^{-1} _{n} g) =\Delta _{I,f,d} (g) $. So the value $\Delta_{I,f,d} (-)$ is also invariant under the second type of Markov moves. So we complete the proof.

\end{pf}

An important example of generalized extended $R-$matrix still come from finite groups.

\begin{thm}
\label{thm:commutingpairsgiveinvariant}
Let $G$ be a finite group. Suppose $b,g\in G$ such that $gb=bg$. Let $V=\mathbb{C} G$, $\phi_{g} \in End (V\otimes V)$ as in Theorem \ref{thm:extendableGRmatrix}. Then the combination
$$(I= \phi_{g} , f=(b)_{L} , c_1 = (bg^{-1})_{L} , c_2 = (gb)_{L} , d= (b)_{L}   )$$
 is a generalized extended $R-$matrix. We denote the resulted invariant of oriented links as $\Delta_{G;(g,b)}(-)$.
\end{thm}

A slightly modification of the construction of section 9 give a stronger link group invariant as follows. Let $L= L_1 \sqcup L_2 \sqcup \cdots \sqcup L_m  $ be a oriented link with $m$ components.  Let $b \in B_{n}$ be a braid whose closure $\hat{b}$ presents $L$. First we introduce a simple definition generalizing extended pair.
\begin{defi}
Let $G$ be any group. A commuting pair of $G$ means a couple $(g,b)\in G\times G$ such that $gb=bg$. We denote the set of all commuting pairs of $G$ as $\mathscr{C}G$.
\end{defi}

Evidently there is $G\subset \mathscr{E}G\subset \mathscr{C}G$. We can similarly consider $\mathscr{C}G-$colored links (braids) , that is , a link (braid) with every component associated with an element of $\mathscr{C}G$.

Now we associate a group to $b $. We introduce variables
$$X_1 , X_2 ,..., X_n , C_1 ,B_1 , C_2 ,B_2,..., C_m ,B_m ,$$
 and let $ \bar{F}_{n,m}$ be the group generated by these elements along with the following relations:
$$  C_i B_i = B_i C_i  , \  for\  1\leq i\leq m  .$$
Then we color the
$i-$th string of $b$ with the commuting  pair $( C_{j_i } , B_{j_i } )$, if it (as a subset of $L= \hat{b}$) belongs to the component $L_{j_i}$. Denote the resulted $\mathscr{C}\bar{F}_{n,m}-$colored braid as $b^{''}$. Recall the  map $\Psi_{g_1 , g_2}$ and $\Psi^{-} _{g_1 , g_2} $
before Theorem \ref{thm:colored link invariant} and (2),(3) of Theorem \ref{thm:extended colored link invariant}. For a $\mathscr{C}\bar{F}_{n,m}-$colored $n-$string simple braid
$X^{+} _{i; (C_{j^{'} }, B_{j^{'}} ) ,( C_{j^{''} } ,B_{j^{''}} )}$ or $X^{-} _{i; (C_{j^{'} }, B_{j^{'}} ) ,( C_{j^{''} } ,B_{j^{''}} ) }  $, whose $i-$th string and $(i+1)-$th string are colored by $(C_{j^{'}} ,B_{j^{'}})$ and $(C_{j^{''}} ,B_{j^{''}})$ respectively, and $(Z_1 , Z_2 ,..., Z_n )\in \underbrace{ \bar{F}_{n,m} \times \bar{F}_{n,m} \times ... \times \bar{F}_{n,m} }_{n}$, we set
\begin{align*}
&\mathcal{F}^{'} ( X^{+} _{ i; (C_{j^{'} } ,B_{j^{'}}) ,( C_{j^{''}}  ,B_{j^{''}})}  )(\cdots , Z_i , Z_{i+1} , \cdots    )=  (\cdots , C_{j^{''}}Z_{i+1} Z_i ^{-1} C_{j^{'} }^{-1} Z_i ,  Z_i ,\cdots   ) , \\
&\mathcal{F} ^{'} ( X^{-} _{i; (C_{j^{'} } ,B_{j^{'}}) ,( C_{j^{''}}  ,B_{j^{''}})} )(\cdots ,Z_i ,Z_{i+1} , \cdots )= (\cdots , Z_{i+1} ,  C_{j^{'}} ^{-1}   Z_i Z_{i+1} ^{-1} C_{j^{''} }Z_{i+1} ,\cdots )  .
\end{align*}

\begin{rem}
 Here the maps are different with the corresponding ones in section 9, as $B_{j^{'}}$ or $B_{j^{"}}$ do not appear in the right side of the identities. Compare (2) of the following Definition \ref{defi:groupinvariantbraidversion3} with (2) of Definition \ref{defi:groupinvariantbraidversion2}, we can understand those difference as "putting all those $B_{j}$'s to the left end".
\end{rem}

If $b^{''} = \Pi _{t=1} ^{N} X^{\epsilon_{t}} _{i_t ; ( C_{j^{'}_{t}} , B_{j^{'} _{t}}) , (C_{ j^{''}_{t} } ,B_{j^{''} _{t}} )  }$, we set
$$\mathcal{F}^{'} (b^{''} )= \mathcal{F}^{'} ( X^{\epsilon_{N}} _{i_N ; (C_{j^{'} _{N}} , B_{j^{'} _{N} } ) , (C_{ j^{''} _{N}  },B_{ j^{''} _{N}   }    ) }\circ \cdots \circ \mathcal{F } ^{'} (  X^{\epsilon_{1}}   _{i_1 ; (C_{j^{'} _{1}} , B_{j^{'} _{1} } ) , (C_{ j^{''} _{1}  },B_{ j^{''} _{1}   }    )  } ) .$$
 We define set of words $\bar{W}^{b} _1 ,...,\bar{W}^{b} _n \in \bar{F}_{n,m }$  by
$$\mathcal{F}^{'} (b^{''} ) ( (X_1 ,X_2 ,..., X_n   )  )= (\bar{W}^{b} _1 , \bar{W}^{b} _2 ,..., \bar{W}^{b} _n   ).$$

 We can view $\mathcal{F}^{'}(b^{''})$ as a group isomorphisms from $\bar{F}_{n,m}$ to itself, by sending generators $(X_1 ,..., X_n , C_1 ,B_1 , ..., C_m ,B_m )$ to $( \bar{W}^{b} _1 ,..., \bar{W}^{b} _n   , C_1 ,B_1 , ..., C_m ,B_m ) $ .
 To show the reliance of $\bar{W}^{b} _i $  on the generators and the order of colors of the strings of $b$, we denote them as $\bar{W}^{b} _i ( X_1 , ..., X_n ;C_{j_1} ,B_{j_1}, ..., C_{j_n} ,B_{j_n}  ) $ .

\begin{defi}
\label{defi:groupinvariantbraidversion3}
Assuming above notations,  we define a group $\tilde{\mathcal{G}}_{L, b}$ (as quotient groups of $\bar{F}_{n,m}$) generated by letters $X_1 ,X_2 ,..., X_n , C_1 ,B_1,C_2 ,B_2,..., C_m ,B_m$,  with the following set of relations.
\begin{enumerate}
\item[(1)] $C_i B_i = B_i C_i $  for  $ 1\leq i\leq m $ ;
\item[(2)] $ X_t = B^{m_t (b)+1} _{j_t}  \bar{W}^{b} _t (X_1 ,..., X_n ;C_{j_1} , B_{j_1},..., C_{j_n}, B_{j_n}  )$ for  $  1\leq t \leq n $ .
\end{enumerate}
Where  the integers $m_t (b )$ in (2) are defined in the place before Theorem \ref{thm:linkinvariantfromgeneralizedextended}.
\end{defi}
\

In similar method  in the proof of Theorem \ref{thm:groupinvariant2}, we can prove the following theorem.
\begin{thm}
The group $\tilde{\mathcal{G}}_{L,b}$ only depends on the link $L$, so it is a group invariant of the link and we simply denote it as $\tilde{\mathcal{G}}_{L}$.
\end{thm}

There is the following proposition similar as Proposition \ref{prop:groupinandnumberin}.Suppose $K$ is a knot. Let $G$ be a finite group and $(g,b)$ a commuting pair of $G$. By definition of the group $\tilde{\mathcal{G}}_{K}$ we have

\begin{prop}
$\Delta_{G;(g,b)} (K) = |\{  f\in Hom( \tilde{\mathcal{G}}_{K} , G    )  | f(C)=g,\ f(B)=b    \} |.$
\end{prop}

Denote the link group of an oriented  link $L= L_1 \sqcup L_2 \sqcup \cdots \sqcup L_m\subset S^3 $ as $\pi_{L}$.  For $1\leq i \leq m$,  the counterclockwise meridian around the $i-$th component determine a conjugacy class of $\pi_{L}$, we chose any representative of this conjugacy class and denote it as $m_i$.

\begin{defi}
\label{defi:reflectivelinkgroup}
For an oriented link as above, we denote the normal group of $\pi_{L}$ generated by $\{ m^2 _i \} _{1\leq i\leq m}$ as $N_{L}$, and define a group $R \pi  _{L}$ by   $R\pi _{L} = \pi_{L}/N_{L}$. Denote the element in $R\pi_{L}$ represented by $m_i$ as $\bar{m}_i$.
\end{defi}

From definition, the relationship between the group $R\pi_{L}$ and the link group $\pi_{L}$ is quite similar with the relationship between Coxeter group and its related Artin group. But we do not know if there are deeper similarities between them.

    We hope the groups $\bar{\mathcal{G}}_{L}$ and $\tilde{\mathcal{G}}_{L}$ are new group invariants for links. But it seems the following facts are true.
\begin{enumerate}
\item[(1)]  There is an isomorphism from the group $\bar{\mathcal{G}}_{L}$  to a free product of $R \pi_{L}$ with certain group, which send the generators $C_i , B_i $ to $l_i ,\bar{m} _i  \in R\pi_{L}$ respectively, where $\bar{m}_i$ is defined in Definition \ref{defi:reflectivelinkgroup}, and $l_i \in R\pi _{L}$ is the equivalent class of suitable longitude element in $\pi_{L}$ along the $i-$th component $L_i$.
\item[(2)] There is an isomorphism from the group $\tilde{\mathcal{G}}_{L}$  to a free product of $\pi_{L}$ with a free group, which send the generators $C_i , B_i $ to $l_i ,m_i  \in \pi_{L}$ respectively, where  $m_i $ are defined as above, and $l_i$ represent suitable longitude cycle along the $i-$th component $L_i$.
\end{enumerate}

We will clarify those facts in our forthcoming papers.

\section{Invariants for three dimensional manifolds}

Based on the invariant of oriented, $G-$colored framed links we constructed in last section, we present a integer invariant for closed three dimensional manifolds. Still fix a finite group $G$.  Because of Corollary \ref{cor:essentialconjugacy}, we consider framed links colored by conjugacy classes of $G$.
As in the beginning of section 2, suppose the set of conjugacy classes of $G$ is $\{ [g_i ]  \}_{i=0,..., N} $. Let $C_{g_i}$ be the centralizer subgroup of $g_i$, and set $v_i = |C_{g_i} |$. Suppose $L= L_1 \sqcup L_2 \sqcup \cdots \sqcup L_m  $ is a $m-$component link where $L_i ( 1\leq i\leq m)$ are its components. Denote $L_{d_1 , d_2 ,..., d_m } $ as the conjugacy class colored link, where the component $L_i $ of $L$ is colored by the conjugacy class $[g_{d_i }]$ for $1\leq i\leq m $. Recall the functor $\mathbb{F}$ constructed in Theorem \ref{thm:colored framed link invariant}. Apply the functor $\mathbb{F}$ to present link $L_{d_1 ,d_2 ,..., d_m }$ we have an integer $\mathbb{F} (L_{d_1 ,d_2 ,..., d_m } ) $. After a comparison of our invariants with the WRT invariants for closed three dimensional manifolds, we can guess that combine these link invariants suitably we can get a manifold invariant. More explicitly, there exist suitable constants $\lambda_0 , \lambda_1 , ..., \lambda_N$ such that the number for $L$:
$$\sum_{d_1 =0 , d_2 =0 ,..., d_m =0 } ^{N,N,..., N} \lambda_{d_1 } \lambda_{d_2} \cdots \lambda_{d_m } \mathbb{F} (L_{d_1 , d_2 ,..., d_m  } ) $$
is invariant under the both two Kirby moves, thus gives a invariant for close three dimensional manifolds. This is certified in the following theorem.

\begin{thm}
\label{thm:3mfdintegerinvariant}
Under above setting, the following function on the set of oriented links is invariant under the Kirby move I and II.
$$\mathscr{F} (L) = \frac{1}{|G|^{m}}  \sum_{d_1 =0 , d_2 =0 ,..., d_m =0 } ^{N,N,..., N} |[g_{d_1}]|| [g_{d_2}]| \cdots |[g_{d_m }]| \mathbb{F} (L_{d_1 , d_2 ,..., d_m  } ) . $$

\end{thm}

\begin{pf}
Suppose $b$ is a $G-$colored braid such that $L= \hat{b}$. By $(2)$ of Proposition \ref{prop:groupinandnumberin}, we have

$  |[g_{d_1}]|| [g_{d_2}]| \cdots |[g_{d_m }]| \mathbb{F} (L_{d_1 , d_2 ,..., d_m  } ) =
 |\{  f\in Hom ( \mathbb{G}_{L,\bar{b}}, G ) | f(Y_1 )\in [g_{d_1} ], ..., f(Y_m ) \in [g_{d_m} ] \}    | $
  So we have
  $  \sum_{d_1 =0 , d_2 =0 ,..., d_m =0 } ^{N,N,..., N} |[g_{d_1}]|| [g_{d_2}]| \cdots |[g_{d_m }]| \mathbb{F} (L_{d_1 , d_2 ,..., d_m  } ) =
  | Hom ( \mathcal{G}_{L ,\bar{b}} , G )  | $.

  So this theorem is a corollary of the following Theorem \ref{thm:groupinvariantmaintheorem}. Because if $L^{'}$ is obtained from from $L$ by adding one trivial knot with $\pm$ framing (that is the Kirby move $\mathcal{K}I$), then $L^{'}$ also has $m+1$ components. By (1) of Theorem \ref{thm:groupinvariantmaintheorem}, $\mathbb{G}_{L^{'}} \cong \mathbb{G}_{L} \ast \mathbb{Z}$. So
  $$| Hom ( \mathbb{G}_{L^{'}} , G  ) | = |G| | Hom ( \mathbb{G}_{L} , G  ) |.$$
  So we have $\mathscr{F}(L^{'}) = \frac{1}{|G|^{m+1}} |G||Hom( \mathbb{G}_{L} , G ) |= \mathscr{F}(L) .$

  On the other hand, if $L^{'}$ is obtained from $L$ by a Kirby move $\mathcal{K}II$, then $L^{'}$ also has $m$ components. By (2) of Theorem \ref{thm:groupinvariantmaintheorem}, we have $\mathbb{G}_{L^{'}} \cong \mathbb{G}_{l} $. Which implies
  $$\mathscr{F}(L^{'}) = \frac{1}{|G|^{m}} |Hom(\mathbb{G}_{L^{'}} , G  )  | =  \frac{1}{|G|^{m}} |Hom(\mathbb{G}_{L} , G  )  |=\mathscr{F}(L). $$
\end{pf}

 \begin{figure}[htbp]

  \centering
    \includegraphics[height=4.2cm]{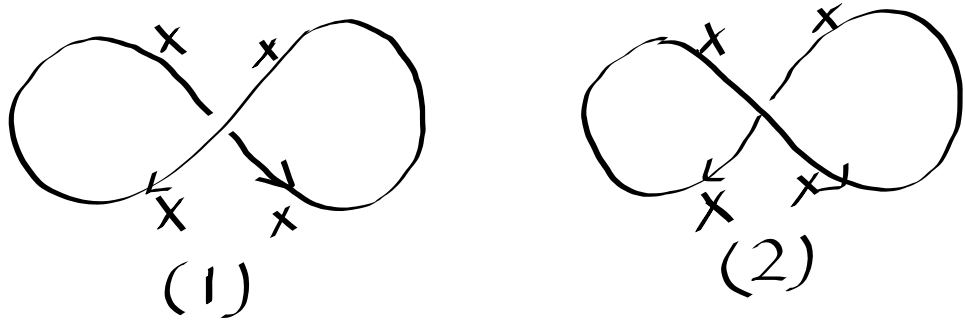}
  \caption{Different twists }
\end{figure}

\begin{thm}
\label{thm:groupinvariantmaintheorem}
Let $L$ be a framed link. Then,

(1) If $L^{'}$ is a framed link obtained from $L$ through a Kirby move $\mathcal{K}I$  by adding a trivial knot with $\pm$ framing, then $\mathbb{G}_{L^{'}} \cong \mathbb{G}_{L} \ast \mathbb{Z} $.

(2) If $L^{'}$ is a framed link obtained from $L$ through a  Kirby move $\mathcal{K}II$, then $\mathbb{G}_{L^{'}} \cong \mathbb{G}_{L} $.
\end{thm}
\begin{pf} First, suppose $L^{'}$ and $L$ are framed links as in (1). Since the group $\mathbb{G}_{L}$ is a invariant for oriented framed links, we can suppose
on the diagram, the graph of the added trivial knot is separate from the graph for $L$. Now if the added trivial knot has framing $+1$ as (1) of Figure 14, since the trivial knot has one long under arc and one crossing point, so according to Proposition \ref{prop:simplifiedpresentation}, the group $\mathbb{G}_{L^{'}}$ has two more generators $X, Y$ where $Y $ is the generator related with colour, and one more relation $X= X X^{-1} Y^{-1} X $ than the group $\mathbb{G}_{L}$. So we have $\mathbb{G}_{L^{'}} \cong \mathbb{G}_{L} \ast \mathbb{Z} $.  If the added trivial knot has framing $-1$ as (2) of Figure 14, then similarly the group
$\mathbb{G}_{L^{'}}$ has two more generators $X,Y$ and one more relation $X=XX^{-1} YX $. So we also have $\mathbb{G}_{L^{'}} \cong \mathbb{G}_{L} \ast \mathbb{Z} $.

Next we suppose $L^{'}$ and $L$ are framed links as in (2). Suppose the components of $L$ are $L_1 ,L_2 ,...,L_{J} $. Suppose the Kirby move from $L$ to
$L^{'}$ is to replace a short section of arc $C$  in $L_i $ with a string $S$ parallel the component $L_j $ (always lying on the right hand side of $L_j$
)as (1) of Figure 15. Denote the components of $L^{'} $ corresponding (in natural sense ) to  the component $L_i$ as $L^{'} _i$ $(1\leq i\leq J)$.

For $1\leq k\leq J$, suppose the long arcs on the component $L_k$ are $S^{1} _{k} ,S^{2} _k ,..., S^{N_k } _{k} $. Suppose the short section of arc
replaced in the Kirby move is on $S^{N_i } _i $. And suppose the long arc on $L_j $ paralleling the head and tail of the new string $S$ is $S^{1} _j $ as shown in (1) of Figure 15.

 \begin{figure}[htbp]

  \centering
  \includegraphics[height=6cm]{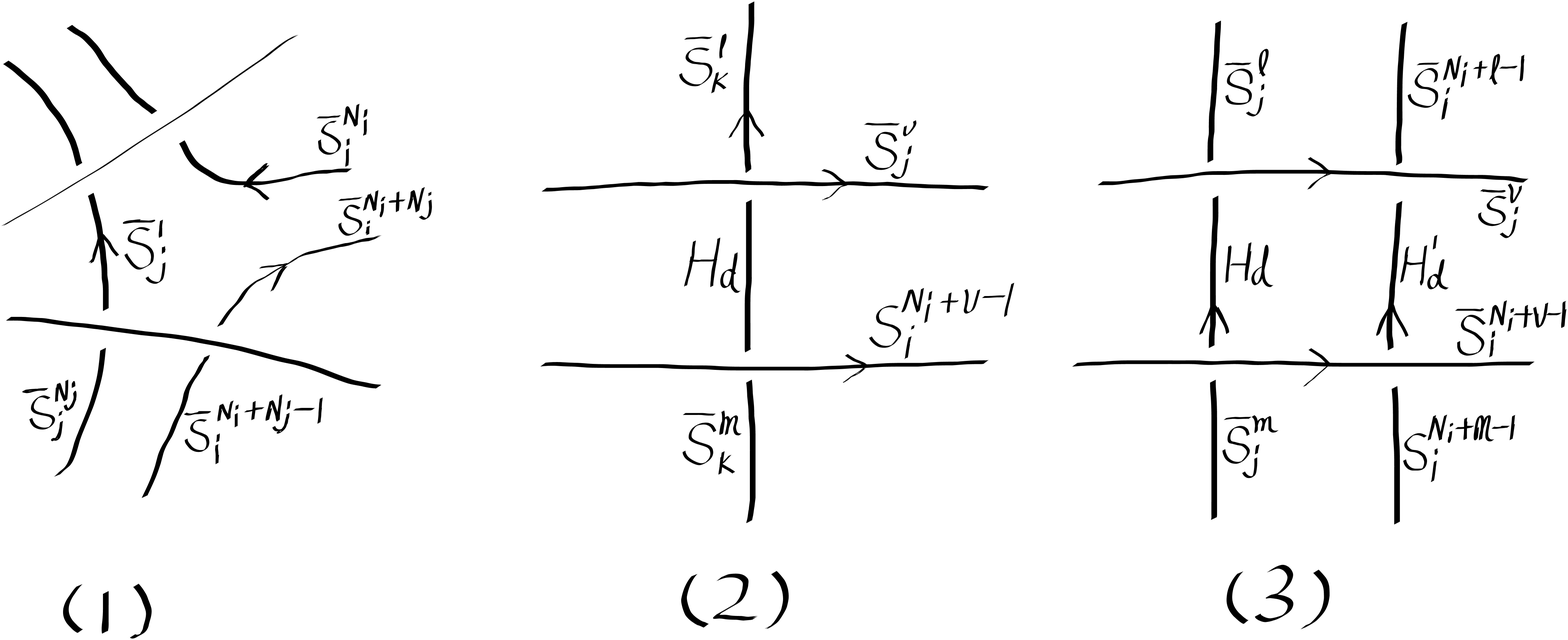}
  \caption{Fig }
\end{figure}

Denote the crossing points on $L$ as $P_1 , P_2 ,..., P_M $. We need to name the long arcs of $L^{'}$ suitably. Intuitively, $L^{'}$ has more long arcs
located on the new string $S$ or cut out by the new string $S$ near the crossing points on which $L_j $ passing over.

case1. For a long arc $S^a _b $ (on $L_b $ ) $b\neq i$, we denote the corresponding long arc on $L^{'}_b $ as $\bar{S}^a _b $.

case2. Denote the long arc on the new string $S$ paralleling $S^{a} _j$ $(1<a\leq N_j )$ as $\bar{S} ^{N_i +a-1} _i $; denote the long arc containing the
head of the new string $S$ as $\bar{S} ^{N_i} _i$, the long arc containing the tail of $S$ as $\bar{S} ^{N_i +N_j } _i $. Then the long arcs of $L^{'} _i$
are $ \bar{S} ^{1} _{i} , ... \bar{S} ^{N_i -1} _i , \bar{S} ^{N_i } _i  ,..., \bar{S} ^{N_i +N_j } _i $.

case3. Now if on a crossing point $P_d $ the component $L_j$ pass over another component, denote the (short) long arc cut by the new string $S$ near $P_d $
as $H_d $,as shown in (2) of Figure 15.

case4. Now if on a crossing point $P_d $ the component $L_j$ pass over itself, then we have two new (short) long arc near $P_d$ parallel to each other .
Denote the one on $\bar{L} _j$ as $ H_d $ , and the one on $\bar{L} _{i} $ as $H^{'} _d$, as shown in (3) of Figure 15.  Then we have nominated all long arcs of $L^{'}$.

Denote the generator of $\mathbb{G} _L$ associated to the long arc $S^{d} _{k} $  as $X^{d} _{k}$, the generator associated with the component $L_k$ as $
Y_k$, for $1\leq k\leq J$, $1\leq d\leq N_k $.

As for the group $\mathbb{G}_{L^{'}}$, for the long arcs $\bar{S} ^a _b $ in case 1 and case 2, we associate a generator $\bar{X} ^{a} _b $. For a long arc
$H_d$ as in case 3, we associate a generator $Z_d$. For long arcs $H_d , H^{'} _d $ appeared in case 4, we associate generator $Z_d , Z^{'} _d$ to them
respectively.

Now we construct a isomorphism from $\mathbb{G}_{L}$ to $\mathbb{G}_{L^{'}} $.

First, for $1\leq d\leq N_j $, we set $T^{d} =\bar{X}^{d} _{j} (\bar{X}^{N_i +d-1} _{i} )^{-1} $, and set $T^{N_j +1} = \bar{X}^{1} _{j} (\bar{X}^{N_i +N_j
} )^{-1} $.

By the following Lemma \ref{lem:aboutT} , we have $T^{l} = T^{m} $ for any $1\leq l,m\leq N_j +1$. So we set $T= T^{1}$ for convenience. Especially, from $T^{1}= T^{N_j
+1}$ we have $\bar{X}^{N_i} _{i} = \bar{X} ^{N_i +N_j } _{i}$.

Now construct a map $\phi$ from the set of generators of $\mathbb{G} _L$ to $\mathbb{G}_{L^{'}}$ as follows.
\begin{enumerate}
\item[(1)] $\phi (X^b _a) = \bar{X}^{b} _a $ for any $a$ and any $1\leq b\leq N_a $;
\item[(2)] $\phi (Y_k ) =\bar{Y}_k  $ for any $k\neq j $;
\item[(3)] $\phi (Y_j )= T \bar{Y}_i T^{-1} \bar{Y}_j  $.
\end{enumerate}

We claim the map $\phi$ extends to a morphism from $\mathbb{G} _L$ to $\mathbb{G}_{L^{'}}$, which is still denoted by $\phi$. It would suffice to check
that the right hand sides of above identities satisfy all relations for the group $\mathbb{G} _L$.

  \begin{figure}[htbp]

  \centering
  \includegraphics[height=5.8cm]{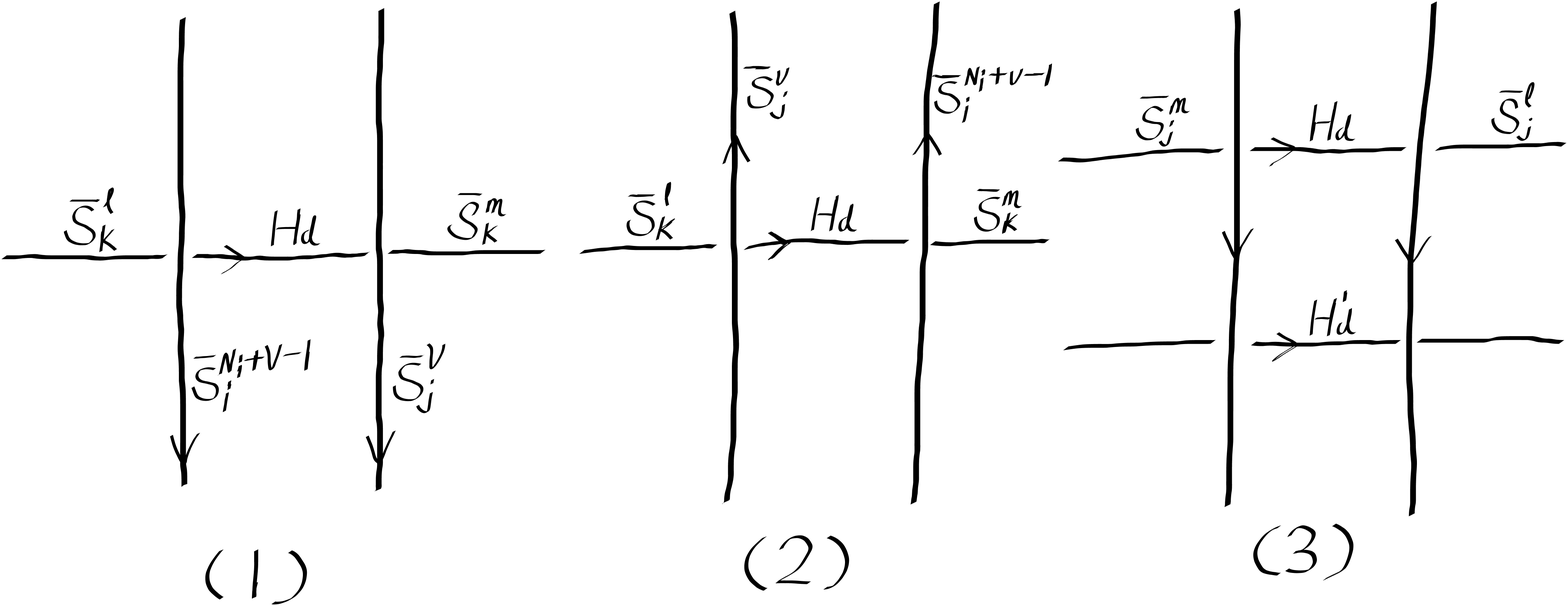}
  \caption{Fig }
\end{figure}

The morphism $\phi$ is actually a isomorphism. We prove this simply by construct the a morphism $\psi$ from $\mathbb{G}_{L^{'}}$ to $\mathbb{G}_{L}$ and
show this morphism is the inverse of $\phi$. The restriction of $\psi$ on  generators of $\mathbb{G}_{L^{'}}$ is as follows. Let $T=X^{1}_j ( X^{N_i -1} _i
)^{-1}$.
\begin{enumerate}
\item[(1)] $ \psi ( \bar{X} ^b _a ) = X^b _a$, except for $a=i$ and $b\leq N_i $;
\item[(2)] $\psi ( \bar{X} ^{N_i +d} _{i} ) = T^{-1} X^{d+1} _{j} $ for $0\leq d\leq N_j -1 $;  $\psi (\bar{X} ^{N_i + N_j} _{i} )= T^{-1} X^{1} _j $;
\item[(3)] $\psi (\bar{Y}_a  )= Y_a $ for $a\neq j$;
\item[(4)] $\psi (\bar{Y}_j )= T Y^{-1} _i T^{-1} Y_j $;
\item[(5)] Suppose $P_d$ is a crossing point where the component $L_j$ pass over another component $L_k$. In cases shown in (1) of Figure 16, we set $ \psi (Z_d ) = X^{m} _{k} (X^{v} _j
)^{-1} Y ^{-1} _{j} X ^{v} _{j}$, and we set $\psi ( Z_d) = X^{l} _k (X^{v} _j )^{-1} Y^{-1} _j X^{v} _j $ in cases shown in (2) of Figure 16 .
\item[(6)] Suppose $P_d$ is a crossing point where the component $L_j$ pass over itself. Suppose the related notations are shown in (3) of Figure 16,  we set $\psi (Z_d) = X^{l} _j (X^{v} _j)^{-1} Y^{-1} _j X^{v} _j $ and
 $\psi (Z^{'} _d  )= T^{-1} X^{l} _j ( X^{v} _j  )^{-1} Y^{-1} _j X^{v} _j $.
\end{enumerate}

 To prove above $\psi$ extends to a group morphism, it suffice to show the right hand sides satisfies the relations of $\mathbb{G} _L$. It is readily to
 prove $\psi \circ \phi = Id_{\mathbb{G} _L}$, and $\phi \circ \psi = Id_{\mathbb{G}_{L^{'}}  } $.  Because certification of these identities are not hard but tedious, let's explain where the elements $"T", "T Y^{-1} _i T^{-1} Y_j  "$ come from instead of giving a complete case by case proof.

  \begin{figure}[htbp]

  \centering
  \includegraphics[height=6cm]{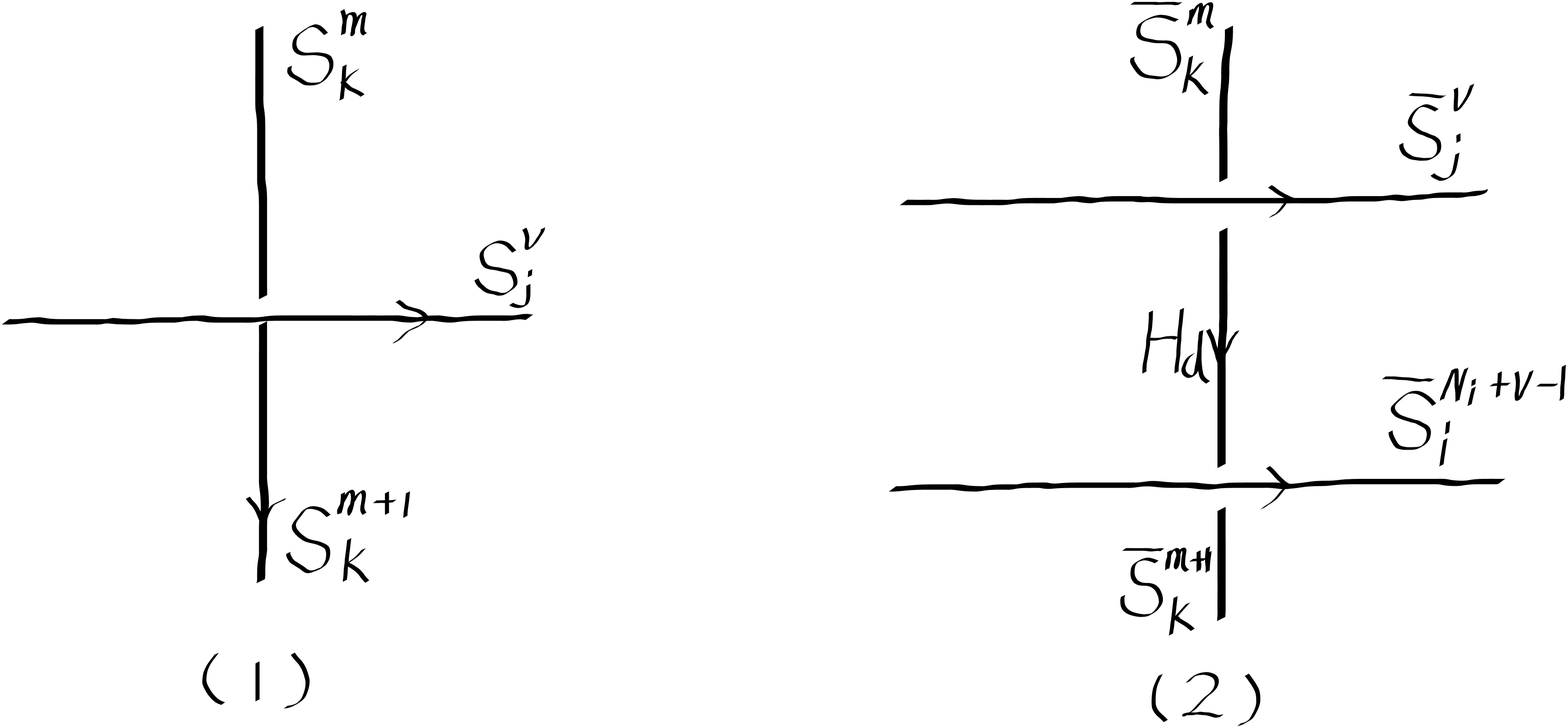}
  \caption{Fig }
\end{figure}

 Suppose a crossing point $P_d$ on the graph of $L$ is as (1) of Figure 17, and suppose (2) of Figure 17 describes the related region of the graph of $L^{'}$. Then from the graph (1) we have a relation for $\mathbb{G}_L $:   $X^{m+1} _{k} = X^{m} _{k} (X^{v} _{j} )^{-1} Y^{-1} _{j} X^{v} _{j} $.

 From the graph (2) we have two relations for the group $\mathbb{G}_{L^{'}}$:

 $Z_d = \bar{X}^{m} _{k} (\bar{X} ^{v} _{j}  )^{-1} \bar{Y} ^{-1} _{j} \bar{X}^{v} _{j} $,  $\bar{X} ^{m+1} _{k} = Z_d ( \bar{X} ^{N_i +v -1 }  )^{-1} \bar{Y}^{-1} _{i} \bar{X} ^{N_i +v-1} _{i} $.

 So we have $X^{m+1} _{k} =  \bar{X}^{m} _{k} (\bar{X} ^{v} _{j}  )^{-1} \bar{Y} ^{-1} _{j} \bar{X}^{v} _{j} ( \bar{X} ^{N_i +v -1 }  )^{-1} \bar{Y}^{-1} _{i} \bar{X} ^{N_i +v-1} _{i} $

 $= \bar{X}^{m} _{k} (\bar{X} ^{v} _{j}  )^{-1} \bar{Y} ^{-1} _{j} \bar{X}^{v} _{j} ( \bar{X} ^{N_i +v -1 }  )^{-1} \bar{Y}^{-1} _{i} \bar{X} ^{N_i +v-1} _{i}
 ( \bar{X} ^{v} _{j} )^{-1} \bar{X}^{v} _{j} $

 $=  \bar{X}^{m} _{k} (\bar{X} ^{v} _{j}  )^{-1} ( T \bar{Y}_i T^{-1} \bar{Y} _{j}   )^{-1}  \bar{X}^{v} _{j}  $.

This identity shows that being passed over by the parallel strings $\bar{S}^{v} _{j} $ and $\bar{S}^{N_i +v-1}_i $ ( colored by $\bar{Y}_j $ and $ \bar{Y}_i $ respectively ) has the same effect of being passed over by a single string $\bar{S}^{v} _{j} $ colored by the element $T \bar{Y}_i T^{-1} \bar{Y} _{j}  $. It explains the morphism $\phi$.

 \begin{lem}
 \label{lem:aboutT}
  Suppose elements $T^v $ are defined as in above proof of Theorem \ref{thm:groupinvariantmaintheorem}, then $T^l = T^m $ for any $m,l$.
 \end{lem}
\begin{pf}
When the related local graph is as (1) of Figure 18, we have

$ \bar{X}^{N_i +v } _i  = \bar{X}^{N_i +v-1 } _{i} ( \bar{X}^{m} _{k} )^{-1} \bar{Y}^{-1} _{k} \bar{X}^{m} _{k} $,
and $\bar{X}^{v+1 } _j  = \bar{X}^{v } _{j} ( \bar{X}^{m} _{k} )^{-1} \bar{Y}^{-1} _{k} \bar{X}^{m} _{k} $.

So we have $T^{v+1} = \bar{X}^{v+1} _j ( \bar{X}^{N_i +v } _i  )^{-1} = T^{v} $.  When the related local graph is as in (2) of Figure 18, we have

$\bar{Y}^{'} _d = \bar{X}^{N_i +v-1 } _{i} ( \bar{X}^{m} _{j} )^{-1} \bar{Y}^{-1} _{j} \bar{X}^{m} _{j}  $,
$\bar{Y} _d = \bar{X}^{v } _{i} ( \bar{X}^{m} _{j} )^{-1} \bar{Y}^{-1} _{j} \bar{X}^{m} _{j}   $,

$\bar{X}^{N_i +v} _{i} = \bar{Y}^{'} _{d} ( \bar{X}^{N_i +m -1} _{i} )^{-1} \bar{Y}^{-1} _{i} \bar{X} ^{N_i +m-1} _{i} $,
$ \bar{X}^{v+1} _{j} = \bar{Y} _{d} ( \bar{X}^{N_i +m -1} _{i} )^{-1} \bar{Y}^{-1} _{i} \bar{X} ^{N_i +m-1} _{i}  $.

So we have $T^{v+1} = \bar{X}^{v+1} _j ( \bar{X}^{N_i +v } _i  )^{-1} = \bar{Y}_{d} (\bar{Y}^{'} _{d} )^{-1} = \bar{X}^{v} _j ( \bar{X}^{N_i +v-1 } _i  )^{-1}  =T^{v} $.

There are other two cases obtained by reversing the horizontal orientations of the arcs in (1) and (2) of Figure 18, these cases can be discussed similarly.

   \begin{figure}[htbp]

  \centering
  \includegraphics[height=5.8cm]{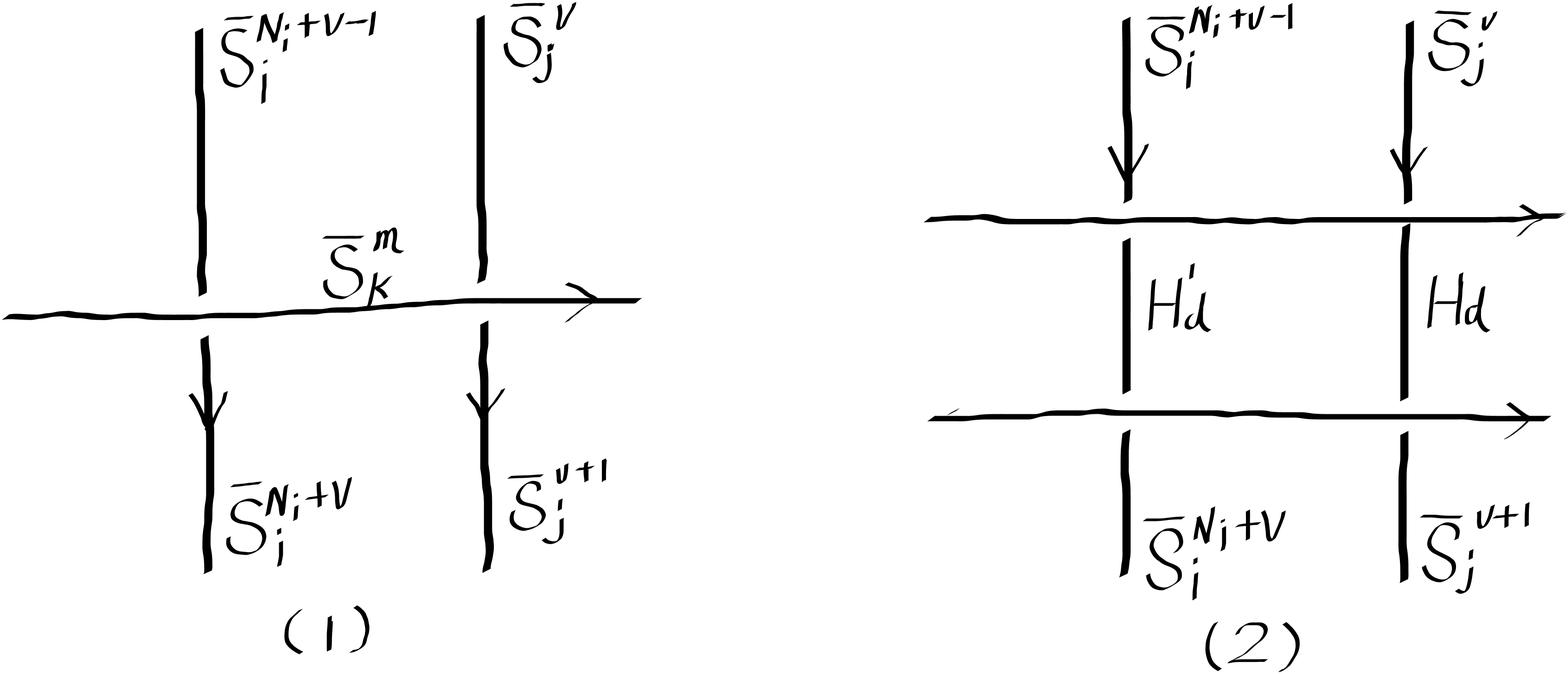}
  \caption{Fig }
\end{figure}

\end{pf}

 Since the group  $\mathbb{G} _L$ is only invariant under the second type of Kirby moves, it is only a invariant of the so called framed three dimensional
 manifolds. We have a very simple method to get a invariant for three dimensional manifolds.

 \begin{thm} Let $L$ be a oriented link diagram. Suppose $\mathbb{G}_L$ is the group defined in Proposition \ref{prop:simplifiedpresentation},suppose $\hat{ \mathbb{G}
 }_{L}$ is a subgroup of $\mathbb{G}_L$ such that $\mathbb{G}_L \cong \hat{ \mathbb{G} }_{L} \ast  F_{N} $ where $F_N$ is a free subgroup
 with maximal possible $N$, and the symbol "$ \ast $" means free product. Then the group $\hat{ \mathbb{G} }_L $ is well defined and invariant
 under both two types of Kirby moves.

 \end{thm}

 \begin{pf}
By the well known Grushko decomposition theorem \cite{Gru}, the group $\hat{ \mathbb{G} }_L $ is uniquely determined thus well defined. It is evident that this group
is invariant under both Kirby moves.
 \end{pf}$\\$

 We hoped the group invariant $\hat{\mathbb{G}}_{L}$ could be a new noncommutative group invariant of three manifolds, but some computations imply that $\mathbb{G}_{L}\cong \pi_1 (M_{L}) \ast F_{N}  $, where $N$ is the number of components in $L$.  We will clarify this in some forthcoming papers. So if this is true, instead of producing a new group invariant of three manifold, we give a "reconstruction" of the fundamental group (free product with a free group),   with a new interesting presentation similar to the Wirtinger presentation of link groups  .

\hspace{-0.70cm} {\sc School of Mathematics}

\nd{\sc Hefei university of technology}

\nd {\sc Hefei 230009 China}

\nd {\sc E-mail addresses}: {\sc Zhi Chen} ({\tt
zzzchen@ustc.edu.cn}).

\end{pf}


\begin{thebibliography}{99}

\bibitem[ADO]{ADO}{\sc Y.Akutsu, T.Deguchi, T.Ohtsuki},Invariants of colored links,{\it Journ. of Knot Theory and Its Rami. }{\bf 1}No.2(1992),161-184.

\bibitem[BK]{BK}{\sc B.Bakalov, J.A.Kirillov }, Lectures on tensor categories and modular functor, {Amer. Math. Socie. }(2000)

\bibitem[Bir]{Bir}{\sc J.S.Birman},Braids, links and mapping class groups, {\it Ann. of Math, Studies.}{\bf 82}, Princeton University Press,(1974)

\bibitem[DW]{DW}{\sc R.Dijkgraaf, E.Witten}, Topological gauge theories and group cohomology, {\it Commun. Math. Phys.  }{\bf 129}(1990),393-429.

\bibitem[Gru]{Gru}{\sc I.A.Grushko},On the bases of a free product of groups, {\it Matematicheskii Sbornik}{\bf 8}(1940),169-182.

\bibitem[Ku]{Ku}{\sc A.G.Kurosh },The theory of groups. Vol 1, Translated and edited by K.A.Hirsch. Chelsea Publ. Co. New York,N.Y (1955).

\bibitem[Oh]{Oh}{\sc T.Ohtsuki  }, Quantum invariants, A study of knots, 3-manifolds, and their sets { World Scientific Publishing Co }(2002).

\bibitem[RT]{RT}{\sc N.Reshetikhin, V.G.Turaev},Invariants of 3-manifolds via link polynomials and quantum groups,{\it Invent.Math. }{\bf 34}(1991),547-597.

\bibitem[Se]{Se}{\sc J.P.Serre}, Linear representations of finite groups, {\it Graduate Texts in Mathematics}{\bf 42}(1977).

\bibitem[Tu1]{Tu1}{\sc V.G.Turaev},The Yang-Baxter equation and invariants of links,{\it Invent. Math.}{\bf 92}(1988),527-553.

\bibitem[Tu2]{Tu2}{\sc V.G.Tureav}, Operator invariants of tangles, and R-matrices,{\it Izv.Akad.Nauk SSSR Ser.Mat.}{\bf 53}(1989),411-444.

\bibitem[FY]{FY}{\sc P.Freed, D.Yetter}, Braided compact closed categories with applications to low dimensional topology,{Adv. in Math.}{\bf77(2)}(1989),79-117.

\bibitem[Kir]{Kir}{\sc R.C.Kirby}, A calculus for framed links in $S^3$, {Invent. Math.}{\bf 45}(1978),35-56.

\bibitem[Wi]{Wi}{\sc E.Witten}, Quantum field theory and the Jones polynomial,{Comm.Math.Phys.} {\bf 121}(1989),351-399.




\end{thebibliography}
   \end{document}